\documentclass[reqno, 11pt]{smfart}

\usepackage{srcltx}
\usepackage[utf8]{inputenc}

\usepackage{a4wide}
\usepackage[author={Max Schlepzig}]{pdfcomment}
\usepackage{url}
\usepackage{amsfonts, amsthm, amsmath, amssymb}
\usepackage{amsmath, amsfonts, amssymb, amsthm}
\usepackage{graphicx, subfigure}
\usepackage{fancyhdr}
\usepackage{stmaryrd}
\usepackage{enumerate}
\usepackage{pstricks}
\usepackage{pstricks-add}
\usepackage{shorttoc}
\usepackage[colorinlistoftodos,french,bordercolor=black,backgroundcolor=red,linecolor=red,textsize=footnotesize,textwidth=0.75in,shadow]{todonotes}
\usepackage{etex}
\reserveinserts{28}
\usepackage{xypic}
\usepackage{dsfont}
\usepackage{cases}

\newcommand{\longhookrightarrow}{}
\DeclareRobustCommand{\longhookrightarrow}{\lhook\joinrel\relbar\joinrel\rightarrow}

\input cyracc.def

\theoremstyle{plain}
\newtheorem{theor}{Théorème}[section]
\newtheorem{lemme}{Lemme}[section]
\newtheorem{prop}{Proposition}[section]
\newcommand{\Mod}[1]{\left({\rm{mod}}\ #1\right)}
\renewcommand{\theequation}{\thesection.\arabic{equation}}

\author{Kevin Destagnol}
\address{Institut de Mathématiques de Jussieu-Paris Rive Gauche\\ UMR 7586\\ Université Paris Diderot-Paris 7\\ Case postale 6052\\ Bâtiment Sophie Germain\\75205 Paris Cedex 13\\ France} 
\email{kevin.destagnol@imj-prg.fr}
\urladdr{webusers.imj-prg.fr/~kevin.destagnol/}
\date{19 décembre 2016}

\begin{document}
\renewcommand{\contentsname}{Sommaire}
\renewcommand{\abstractname}{Abstract}
\setcounter{section}{0}
\rhead{}
\chead{}
\lhead{}
\fancyhead[CO]{\footnotesize{ \textsc{La conjecture de Manin pour une famille de variétés en dimension supérieure}}}
\fancyhead[CE]{\footnotesize{\textsc{K. Destagnol}}}

\title{\textbf{La conjecture de Manin pour\\ une famille de variétés en dimension supérieure}}
\maketitle
\begin{abstract} 
Inspired by a method of La Bretèche relying on some unique factorisation, we generalize work of Blomer, Brüdern, and Salberger to prove Manin's conjecture in its strong form conjectured by Peyre for some infinite family of varieties of higher dimension. The varieties under consideration in this paper correspond to the projective varieties defined by the following equation
$$
x_1 y_2y_3\cdots y_n+x_2y_1y_3 \cdots y_n+ \cdots+x_n y_1 y_2 \cdots y_{n-1}=0.
$$
in $\mathbb{P}^{2n-1}_{\mathbb{Q}}$ for all $n \geqslant 3$. This paper comes with an Appendix by Per Salberger.
\end{abstract}

\setcounter{tocdepth}{3}

\begin{center}
\tableofcontents
\end{center}

\vspace{-7.8mm}
\section{Introduction et principaux résultats}

\subsection{Introduction}
Soit $n \geqslant 2$ un entier fixé et $W_n$ l'hypersurface normale projective de $\mathbb{P}^{2n-1}_{\mathbb{Q}}$ définie par l'équation
\begin{eqnarray}
x_1 y_2y_3\cdots y_n+x_2y_1y_3 \cdots y_n+ \cdots+x_n y_1 y_2 \cdots y_{n-1}=0
\label{eq}
\end{eqnarray}
où $(x_1,\dots,x_n,y_1,\dots,y_n)$ désignent des coordonnées homogènes de $\mathbb{P}^{2n-1}_{\mathbb{Q}}$. Lorsque $n=2$, la variété $W_2$ est lisse tandis que pour $n \geqslant 3$, la variété $W_n$ est singulière, le lieu singulier étant donné par la réunion des sous-espaces fermés définis par les équations
$$
y_i=y_j=y_k=0 \quad \mbox{pour} \quad 1 \leqslant i<j<k\leqslant n
\quad \mbox{ou} \quad
y_i=y_j=x_i=x_j=0 \quad \mbox{pour} \quad 1 \leqslant i<j\leqslant n.
$$
On pose $U_n$ l'ouvert de Zariski de $W_n$ défini par la condition $y_1y_2\cdots y_n \neq 0$. Sur cet ouvert, on peut réécrire l'équation (\ref{eq}) sous la forme
$$
\sum_{i=1}^n \frac{x_i}{y_i}=0.
 $$
 On constate alors que toutes les sous-variétés accumulatrices et les points singuliers de $W_n$ sont inclus dans $W_n \smallsetminus U_n$ (voir \cite{Blomer2014}). La variété $W_n$ possède la structure algébrique suivante qui s'avérera très utile dans la suite. On considère le groupe algébrique
$$
H=\left\{ \begin{pmatrix}
b&a\\
0&b\\
\end{pmatrix} \hspace{0.5mm} : \hspace{0.5mm} (a,b) \in \mathbb{Q}\times \mathbb{Q}^{\ast}
\right\}
$$
ainsi que
$$
\Psi_n: \left\{
\begin{array}{ccc}
H^n &\longrightarrow& \mathbb{G}_a\\
\left( \begin{pmatrix}
b_1&a_1\\
0&b_1\\
\end{pmatrix},\dots,\begin{pmatrix}
b_n&a_n\\
0&b_n\\
\end{pmatrix} \right) & \longmapsto& \overset{n}{\underset{i=1}{\sum}} \frac{a_i}{b_i}\\
\end{array}
\right.
$$
et $G_n=\mbox{Ker}(\Psi_n)/ \mathbb{G}_m$ où l'on utilise le plongement
$$
\begin{array}{ccc}
\mathbb{G}_m &\longhookrightarrow& H^n\\
b&\longmapsto& \left( \begin{pmatrix}
b&0\\
0&b\\
\end{pmatrix},\dots,\begin{pmatrix}
b&0\\
0&b\\
\end{pmatrix}  \right)\\
\end{array}.
$$
Il vient alors que $G_n \cong U_n$. En effet, tout point de $U_n$ admet un unique représentant de la forme $[x_1:\cdots:x_n:1:y_2:\cdots:y_n]$ et de même tout point de $G_n$ possède une unique représentation du type $ \left( \begin{pmatrix}
1&a_1\\
0&1\\
\end{pmatrix},\dots,\begin{pmatrix}
b_n&a_n\\
0&b_n\\
\end{pmatrix}  \right)$ si bien que
$$
\begin{array}{ccc}
U_n &\longrightarrow& G_n\\
{[x_1:\cdots:x_n:1:y_2:\cdots:y_n]}&\longmapsto&
\left( \begin{pmatrix}
1&x_1\\
0&1\\
\end{pmatrix},\dots,\begin{pmatrix}
y_n&x_n\\
0&y_n\\
\end{pmatrix}  \right)\\
\end{array}
$$
est un isomorphisme qui munit $U_n$ d'une structure de groupe commutatif lorsque la multiplication de deux points $(x_1,\dots,x_n,y_1,\dots,y_n)$ et $(x'_1,\dots,x'_n,y'_1,\dots,y'_n)$ est donnée par
$$
(x_1y'_1+x'_1y_1,\dots,x_ny'_n+x'_ny_n,y_1y'_1,\dots,y_ny'_n).
$$
On obtient alors une immersion ouverte naturelle
$
j:G_n \hookrightarrow W_n
$
et une action naturelle $\alpha:G_n \times W_n \rightarrow W_n$ de~$G_n$ sur $W_n$. Ce groupe 
$$
G_n \cong \left( \mathbb{G}_a\right)^{n-1} \times \left( \mathbb{G}_m\right)^{n-1}
$$
est un produit de groupes additifs et multiplicatifs et les techniques d'analyse harmonique utilisées par exemple dans \cite{BT} et \cite{CLT} ne s'adaptent pas \textit{directement} à ce cas même s'il est possible qu'elles puissent s'y généraliser. La méthode utilisée dans cet article repose sur une approche différente de la conjecture de Manin, à savoir une descente sur le torseur versel et avec \cite{Blomer2014}, il s'agit vraisemblablement du seul exemple de groupe de ce type pour lequel les conjectures de Manin-Peyre sont établies. Dans \cite{Blomer2014}, Blomer, Brüdern et Salberger établissent la conjecture de Manin sous sa forme forte conjecturée par Peyre pour la variété $W_3$ pour la hauteur anticanonique
$$
H\left([x_1:x_2:x_3:y_1:y_2:y_3]\right)=\max_{1\leqslant i \leqslant 3}\max(|x_i|,|y_i|)^3
$$
lorsque $(x_1,x_2,x_3,y_1,y_2,y_3)$ sont des entiers premiers entre eux. Leur méthode repose sur la construction d'une résolution crépante $X \rightarrow W_3$ de $W_3$ puis sur une descente sur le torseur versel de~$X$. Le problème se réduit alors à un problème de comptage de points d'un réseau de dimension $10$ dans une région dont la frontière rend le traitement difficile. Ce problème de géométrie des nombres est alors traité à l'aide de séries de Dirichlet multiples et de transformations de Mellin multidimensionnelles. La résolution crépante construite n'est plus valable lorsque $n \geqslant 4$ et leur méthode de comptage se complique considérablement dès que $n \geqslant 4$. Il est cependant à noter que, pour $n=3$, leur formule asymptotique (\cite[theorem 1]{Blomer2014}) est plus précise que la conjecture de Manin à proprement parler puisqu'ils obtiennent un terme principal de la forme $BQ(\log(B))$ avec $Q$ un polynôme de degré 4 dont le coefficient dominant correspond à la constante de Peyre. Enfin, pour $n \geqslant 4$, les auteurs indiquent dans \cite{Blomer2014} sans donner de détails que leur méthode de comptage peut se généraliser afin de donner l'équivalent prédit par la conjecture de Manin sans terme d'erreur. On donne dans cet article une démonstration des conjectures de Manin-Peyre pour tout $n \geqslant 3$ en utilisant une méthode de comptage différente adaptée à la combinatoire du problème.\\
\indent
Comme remarqué dans ce même article \cite[section 1.3]{Blomer2014}, l'équation (\ref{eq}) définit également pour $n \geqslant 3$ une variété singulière $\widehat{W_n}  \subset \left(\mathbb{P}^{1}_{\mathbb{Q}}\right)^n$ pour laquelle il est intéressant d'étudier les conjectures de Manin et Peyre pour la hauteur anticanonique
$$
\widehat{H}\left([x_1:y_1],\dots,[x_n:y_n]\right)=\prod_{i=1}^n\max (|x_i|,|y_i|)
$$
lorsque $(x_i,y_i)$ sont deux entiers premiers entre eux pour tout $i \in \llbracket 1,n\rrbracket$. \`A la connaissance de l'auteur, ces deux conjectures ne sont connues pour aucune valeur de $n \geqslant 3$ pour les variétés $\widehat{W}_n$. Ces dernières ont pourtant un intérêt arithmétique puisqu'elles interviennent dans la détermination du cardinal des matrices stochastiques à coefficients rationnels tous de hauteur inférieure à une certaine borne $B \geqslant 1$. Dans \cite{Sh}, Shparlinski obtient une borne supérieure du bon ordre de grandeur pour ce cardinal et dans \cite{Br}, La Bretèche parvient à en obtenir une formule asymptotique en déterminant le nombre de points rationnels de hauteur inférieure à $B$ sur $\widehat{W}_n$ lorsque $n \geqslant 3$ pour la hauteur
$$
\widehat{H}_n\left([x_1:y_1],\dots,[x_n:y_n]\right)=\max_{1 \leqslant i \leqslant n} (|x_i|,|y_i|)
$$
lorsque $(x_i,y_i)$ sont deux entiers premiers entre eux pour tout $i \in \llbracket 1,n\rrbracket$. Il est cependant important de noter que cette hauteur n'est pas anticanonique et qu'ainsi les conjectures de Manin et Peyre ne s'appliquent pas dans ce cas.\\
\indent
Enfin, l'équation (\ref{eq}) définit également pour $n \geqslant 3$ une variété singulière biprojective $\widetilde{W}_n \subset \left(\mathbb{P}^{n-1}_{\mathbb{Q}}\right)^2$ pour laquelle il peut également être intéressant d'étudier les conjectures de Manin et Peyre par rapport à la hauteur anticanonique
$$
\widetilde{H}\left([x_1:\dots:x_n],[y_1:\dots:y_n]\right)=\max_{1 \leqslant i \leqslant n} |x_i|^{n-1}\max_{1 \leqslant i \leqslant n} |y_i|
$$
lorsque $(x_1,\dots,x_n)$ et $(y_1,\dots,y_n)$ sont deux $n$-uplets d'entiers premiers entre eux. Il est à noter que, contrairement au cas de $W_n$, la variété $\widetilde{W}_n$ ne peut pas être écrite de façon naturelle comme la compactification équivariante d'un groupe. Dans le cas de $W_n$, c'est un aspect essentiel de la preuve des conjectures de Manin-Peyre de \cite{Blomer2014} et du présent travail. Le seul cas pour lequel un résultat est connu est le cas $n=3$. La variété~$\widetilde{W}_3$ est alors une cubique de dimension 3 et Blomer et Brüdern obtiennent dans un premier temps dans \cite{BB} le bon ordre de grandeur pour le problème de comptage associé avant de parvenir récemment avec Salberger à obtenir une formule asymptotique et les conjectures de Manin-Peyre dans \cite{BBS}. Les méthodes utilisées présentent des similarités avec celles de \cite{Blomer2014} et reposent là encore sur une descente sur le torseur versel associée à une résolution crépante de $\widetilde{W}_n$. Blomer, Brüdern et Salberger font en revanche appel à des techniques d'analyse de Fourier pour compter les points entiers sur ce torseur. La différence principale en termes de géométrie de la résolution crépante provient de l'absence de structure naturelle de groupe.\\
\indent

\subsection{Résultats}
L'objet de cet article est de démontrer les conjectures de Manin et Peyre dans le cas de $W_n$ pour tout $n \geqslant 2$. Puisque $W_2$ est une variété torique lisse, le cas $n=2$ est inclus dans les travaux généraux de Batyrev et Tschinkel sur les variétés toriques lisses \cite{BT} et le cas $n=3$ est couvert par le résultat de \cite{Blomer2014}. Cet article est accompagné d'une annexe de Per Salberger qui explicite une résolution crépante pour la variété $W_n$ pour tout $n \geqslant 3$.\\
\par
Pour tout $n \geqslant 3$ et pour un point $(\mathbf{x},\mathbf{y}) \in \mathbb{P}^{2n-1}_{\mathbb{Q}}$ représenté par $(\mathbf{x},\mathbf{y}) \in \mathbb{Z}^{2n}$ premiers entre eux, on considère la hauteur
$$
H(\mathbf{x},\mathbf{y})=\max_{1\leqslant i \leqslant n}(|x_i|,|y_i|)^n
 $$
 qui est une hauteur anticanonique naturelle sur $W_n$ comme remarqué dans \cite{Blomer2014}.
 De plus, pour tout $B \geqslant 1$, on pose
$$
N(B;U_n)=\#\{ (\mathbf{x},\mathbf{y}) \in U_n(\mathbb{Q}) :  H(\mathbf{x},\mathbf{y}) \leqslant B\}.
$$

\noindent
Lorsque $n \geqslant 3$, notre résultat est alors le suivant.
\begin{theor}
Soient $B \geqslant 2$ et  $n \geqslant 3$. La variété $W_n$ est "presque de Fano" au sens de \cite[Définition 3.1]{P03} et il existe une constante $c_n>0$ telle que
$$
N(B;U_n)=c_nB(\log B)^{2^n-n-1}+O\left(B(\log B)^{2^n-n-2}\log(\log B)\right).
$$
De plus, l'expression de la constante $c_n$ est en accord avec la conjecture de Peyre.
\label{theor1}
\end{theor}
\noindent
 \textbf{Remarques.} 
\hspace{-0.5cm}
\begin{minipage}[t]{14.2cm}
 \begin{itemize}
\item Une expression explicite de la constante $c_n$ est donnée par la formule (\ref{cn}) en section 3.2.
\item Le Théorème \ref{theor1} permet de retrouver l'équivalent asymptotique qui découle de \cite[Theorem 1]{Blomer2014} pour $n=3$. La démonstration donnée ici est en revanche différente de celle de \cite{Blomer2014}. 
\item Une adaptation simple de la démonstration de ce théorème conduirait à l'étude des conjectures de Manin-Peyre sur l'hypersurface $W_{\mathbf{a},n}$ de $\mathbb{P}^{2n-1}$ définie par l'équation
$$
a_1x_1 y_2y_3\cdots y_n+a_2x_2y_1y_3 \cdots y_n+ \cdots+a_nx_n y_1 y_2 \cdots y_{n-1}=0,
$$
pour un $n$-uplet $(a_1,\dots,a_n)$ d'entiers tous non nuls. On constate clairement que l'hypersurface $W_{\mathbf{a},n}$ est isomorphe à $W_n$ et par conséquent on ne détaillera pas cet aspect.

\end{itemize}   
  \end{minipage}

$$
{}
$$
\noindent
\textbf{{Remerciements.--}} L'auteur tient à exprimer ici toute sa gratitude à son directeur de thèse, Régis de la Bretèche, pour ses conseils, son soutien et ses relectures tout au long de ce travail, ainsi qu'à Marc Hindry et Tim Browning pour quelques discussions éclairantes. L'auteur tient également à remercier chaleureusement Per Salberger pour lui avoir communiqué ses résultats concernant la résolution crépante de $W_n$ qui figurent en annexe et pour avoir été à l'origine de nombreux éclaircissements ainsi que Valentin Blomer et Jörg Brüdern.

 \section{Notations}
On introduit dans cette section des notations qui seront utilisées tout au long de cet article. On notera $\mathbb{N}$ l'ensemble des entiers positifs non nuls, $\mbox{pgcd}(m,n)$ le pgcd de deux entiers $m$ et $n$, $[m,n]$ leur ppcm et, lorsqu'il existe, $\overline{n}^a$ l'inverse de $n$ modulo un entier $a$. On omettra la dépendance en $a$ et on notera plutôt $\overline{n}$ dans les cas où aucune ambiguïté n'est possible sur l'entier $a$. Pour $n \geqslant 2$, on considère l'entier $N=2^n-1$ ainsi que pour tout $h \in \llbracket 1,N\rrbracket$, son développement en binaire
$$
h=\sum_{1 \leqslant j \leqslant n} \varepsilon_j(h)2^{j-1},
$$ 
avec $\varepsilon_j(h) \in \{0,1\}$. On notera $s(h)=\sum_{j \geqslant 1} \varepsilon_j(h)$ la somme des chiffres en base 2 de $h$. On utilise dans la suite de cet article la variante de la définition des ensembles de décomposition unique de \cite{5} et \cite{4} introduite dans \cite{Br}. On dit qu'un entier $h$ est dominé par $\ell$ (respectivement strictement dominé) si pour tout $j \in \mathbb{N}$, on a $\varepsilon_j(h) \leqslant \varepsilon_j(\ell)$ (respectivement $\varepsilon_j(h) < \varepsilon_j(\ell)$). On notera $h \preceq \ell$ (respectivement $h \prec \ell$). On dira qu'un $N-$uplet $(z_{1},\dots,z_N)$ de $\llbracket 1,n \rrbracket$ est réduit si $\mbox{pgcd}(z_h,z_{\ell})=1$ lorsque $h \not\preceq \ell$ et $\ell \not\preceq h$. On a alors le lemme fondamental suivant.

\begin{lemme}
Il existe une bijection entre l'ensemble des $n-$uplets $(y_1,\cdots,y_n)$ de $\llbracket 1,n \rrbracket$ et les $N-$uplets réduits de $\llbracket 1,n \rrbracket $ tels que
$$
\forall j \in \llbracket 1,n \rrbracket, \quad y_j=\prod_{1 \leqslant h \leqslant N} z_h^{\varepsilon_j(h)} \quad \mbox{et} \quad [y_1,\dots,y_n]=\prod_{1 \leqslant h \leqslant N} z_h.
$$
\label{lemmefacto}
\end{lemme}
\noindent
\textit{Démonstration--}
Il suffit de définir les $z_h$ à $s(h)$ décroissant. Soit $k \in \llbracket 1,n-1\rrbracket$. On pose $z_{N}=\mbox{pgcd}(y_1,\dots,y_n)$ et supposons les $z_h$ construits pour $s(h)\geqslant k+1$. On pose alors pour $h$ tel que $s(h)=k$
$$
z_h={\rm{pgcd}}\left(\frac{y_j}{\underset{1\leqslant \ell \leqslant N, s(\ell) \geqslant k+1 \atop \varepsilon_j(\ell) =1}{\prod}z_{\ell}} \quad : \quad \varepsilon_j(h) =1 \right).
$$
Enfin, lorsque $h=2^{i-1}$ pour $i \in \llbracket 1,n \rrbracket$, on pose
$$
z_h=\frac{y_i}{\underset{1\leqslant \ell \leqslant N, s(\ell) \geqslant 2 \atop \varepsilon_i(\ell) =1}{\prod}z_{\ell}}.
$$
Il est alors facile de vérifier que le $N$-uplet $(z_{1},\dots,z_N)$ de $\llbracket 1,n \rrbracket$ est réduit.
\hfill
$\square$\\
\newline

\section{Démonstration du Théorème \ref{theor1}}


Cette section est consacrée à la preuve de la formule asymptotique du Théorème~\ref{theor1}. On montrera ensuite que la variété $W_n$ est "presque de Fano" au sens de \cite[Définition 3.1]{P03} et que la constante $c_n$ obtenue est en accord avec la conjecture de Manin dans sa forme forte conjecturée par Peyre en section 4. La méthode employée, qui consiste à compter d'abord les $\mathbf{x}$ à $\mathbf{y}$ fixés n'est valable que lorsque la borne sur les $x_i$ est plus grande que celle sur les $y_i$ pour $i \in \llbracket 1,n \rrbracket$. C'est en particulier pourquoi elle ne s'adapte pas, en tout cas directement, aux variétés définies lors de l'introduction $\widehat{W}_n$ et $\widetilde{W}_n$.
\subsubsection{Réduction au cas $y_i \geqslant 1$ pour tout $i \in \llbracket 1,n\rrbracket$}
Quitte à changer le signe des $x_i$ pour $i \in \llbracket 1,n\rrbracket$, on peut supposer les $y_i$ positifs pour tout $i \in \llbracket 1,n\rrbracket$ et obtenir l'égalité
$$
N(B;U_n)=2^{n-1}\#\left\{ (\mathbf{x},\mathbf{y}) \in\mathbb{Z}^{n} \times \mathbb{N}^n \quad :  
\begin{array}{c}
\underset{1\leqslant i \leqslant n}{\max}(|x_i|,y_i) \leqslant B^{1/n},\\[1mm]
 (x_1,\dots,x_n,y_1,\dots,y_n)=1,\\[1mm]
(\mathbf{x},\mathbf{y}) \hspace{2mm} {\rm{ v\acute{e}rifie }} \hspace{2mm} (\ref{eq})
\end{array}
\right\}.
$$
Une inversion de Möbius fournit alors
$$
N(B;U_n)=2^{n-1}\sum_{k=1}^{+\infty}\mu(k) N\left(\frac{B}{k^n}\right)
$$
où
$$
\begin{aligned}
N(B)&=\#\left\{ (\mathbf{x},\mathbf{y}) \in\mathbb{Z}^{n} \times \mathbb{N}^n \quad :  
\begin{array}{c}
\underset{1\leqslant i \leqslant n}{\max}(|x_i|,y_i) \leqslant B^{1/n},\\[1mm]
(\mathbf{x},\mathbf{y}) \hspace{2mm} {\rm{ v\acute{e}rifie }} \hspace{2mm} (\ref{eq})
\end{array}
\right\}\\[2mm]
&=\sum_{\mathbf{y} \in \mathbb{N}^n \atop 1 \leqslant y_i\leqslant \sqrt[n]{B}} N_{\mathbf{y}}\left(B^{1/n}\right),\\
\end{aligned}
$$
avec, pour tous $\mathbf{y} \in \mathbb{N}^n$ et $X \geqslant 1$,
\begin{eqnarray}
N_{\mathbf{y}}(X)=\#\left\{\mathbf{x} \in\mathbb{Z}^{n} \quad :
\begin{array}{c}
\underset{1\leqslant i \leqslant n}{\max}|x_i| \leqslant X,\\[1mm]
(\mathbf{x},\mathbf{y}) \hspace{2mm} {\rm{ v\acute{e}rifie }} \hspace{2mm} (\ref{eq})
\end{array}
\right\}.
\label{Ny}
\end{eqnarray}

\subsubsection{Utilisation du Lemme \ref{lemmefacto}}
On se place ici dans le cas où $\mathbf{y} \in \mathbb{N}^n$ est fixé. L'équation (\ref{eq}) peut se réécrire en utilisant le Lemme \ref{lemmefacto} sous la forme
\begin{eqnarray}
\sum_{j=1}^n d_jx_j=0
\quad
\mbox{avec}
\quad 
\forall i \in \llbracket 1,n \rrbracket, \quad d_i=\prod_{1 \leqslant h \leqslant N} z_h^{1-\varepsilon_i(h)}.
\label{torsor}
\end{eqnarray}
On obtient ainsi la relation de divisibilité $z_{2^{j-1}} \mid x_j$ pour tout $j \in \llbracket 1,n \rrbracket$. Mais contrairement au cas de \cite{Br}, on ne peut pas ici en déduire que $z_{2^{j-1}}=1$ puisqu'on n'a pas les conditions $\mbox{pgcd}(x_j,y_j)=1$.\\
On considère alors, pour $X \geqslant 1$, l'ensemble
$$
\mathcal{A}(\mathbf{y};X)=\left\{ \boldsymbol{\alpha} \in \mathbb{Z}^n \quad :  \quad
\max_{1\leqslant i \leqslant n}|\alpha_i|\leqslant X,\quad
\sum_{i=1}^n d_i\alpha_i=0
\right\}.
$$
On définit pour tout $r\geqslant 1$, l'ensemble $\mathcal{A}_r(\mathbf{y},X)$ par
\begin{eqnarray}
\mathcal{A}_r(\mathbf{y},X)=\left\{
(\alpha_{r+1},\dots,\alpha_n) \in \mathbb{Z}^{n-r}  \quad :  \quad
\max_{r+1\leqslant i \leqslant n}|\alpha_i|\leqslant X,\quad
\sum_{i=r+1}^n d_i\alpha_i \equiv 0 \Mod{d_{1,r}}
\right\}
\label{ar}
\end{eqnarray}
où l'on a posé
\begin{eqnarray}
\forall r \in \llbracket 2,n \rrbracket, \quad d_{1,r}=\prod_{\varepsilon_1(h)=\cdots=\varepsilon_r(h)=0} z_h.
\label{dr}
\end{eqnarray}
De plus, on introduit les notations
\begin{eqnarray}
d'_1=\prod_{\varepsilon_1(h)=0 \atop \varepsilon_2(h)=1} z_h \quad \mbox{et} \quad \forall r \in \llbracket 2,n \rrbracket, \quad d_1^{(r-1)}=\prod_{\varepsilon_1(h)=\cdots=\varepsilon_{r-1}(h)=0 \atop \varepsilon_r(h)=1} z_h
\label{dprime}
\end{eqnarray}
de sorte que $d_{1,r-1}=d_{1,r} d_1^{(r-1)}$. Enfin, on pose
\begin{eqnarray}
\forall r \in \llbracket 2,n \rrbracket, \quad d'_r=\prod_{\varepsilon_r(h)=0 \atop \varepsilon_1(h)+ \cdots+\varepsilon_{r-1}(h) \neq 0 } z_h.
\label{dpprime}
\end{eqnarray}
On obtient ainsi que $d_r=d_{1,r}d'_r$ et on constate en utilisant la section 2 que pour tout $2\leqslant r \leqslant n$, on a $\mbox{pgcd}\left(d_1^{(r-1)},d'_r\right)=1$.\\
\newline
\indent
On démontre alors les lemmes suivants qui permettent d'estimer le cardinal de l'ensemble $\mathcal{A}(\mathbf{y};X)$ pour $X \geqslant 1$.

\begin{lemme}
Soient $1 \leqslant r \leqslant n-1$ et $X \geqslant 1$. Avec les notations (\ref{ar}) et (\ref{dprime}), on a l'estimation
$$
\#\mathcal{A}_r(\mathbf{y},X)=2^{n-r}\overset{n}{\underset{j=r+1}{\prod}}\frac{X}{ d_1^{(j-1)}}+O\left(X^{n-r-1}\right).
$$
\label{lemme1}
\end{lemme}
\noindent
\textit{Démonstration--} 
Soit $r \geqslant 1$. L'idée principale de la preuve, inspirée de \cite{Br}, est de relier le cardinal de $\mathcal{A}_{r+1}(\mathbf{y},X)$ à celui de $\mathcal{A}_r(\mathbf{y},X)$. Démontrons pour ce faire que pour tout $1\leqslant r \leqslant n-2$, on a
\begin{eqnarray}
\#\mathcal{A}_{r}(\mathbf{y},X)=\left(\frac{2X}{d_1^{(r)}}+O(1)\right)\#\mathcal{A}_{r+1}(\mathbf{y},X).
\label{rec}
\end{eqnarray}
En effet, pour un point $(\alpha_{r+1},\dots,\alpha_n)$ de $\mathcal{A}_r(\mathbf{y},X)$, $(\alpha_{r+2},\dots,\alpha_n)$ appartient à $\mathcal{A}_{r+1}(\mathbf{y},X)$. D'autre part, considérons $(\alpha_{r+2},\dots,\alpha_n)$ appartenant à $\mathcal{A}_{r+1}(\mathbf{y},X)$. Il existe donc $k \in \mathbb{Z}$ tel que
\begin{eqnarray}
\sum_{i=r+2}^n d_i\alpha_i =d_{1,r+1}k.
\label{rela}
\end{eqnarray}
Le $(n-r-1)$-uplet $(\alpha_{r+2},\dots,\alpha_n)$ provient de la projection sur les $n-r-1$ dernières coordonnées d'un $(n-r)$-uplet $(\alpha_{r+1},\dots,\alpha_n)$ de $\mathcal{A}_r(\mathbf{y},X)$ s'il existe un entier $|\alpha_{r+1}| \leqslant X$ et $\ell \in \mathbb{Z}$ tel que
$$
\sum_{i=r+1}^n d_i\alpha_i =d_{1,r}\ell.
$$
Grâce à (\ref{rela}), on en déduit la relation
$$
d_{r+1}\alpha_{r+1}+d_{1,r+1}k=d_{1,r}\ell,
$$
qui, au vu de (\ref{dprime}) et (\ref{dpprime}), se réécrit sous la forme
$$
\ell d_1^{(r)}-\alpha_{r+1}d'_{r+1}=k.
$$
On en déduit que $d'_{r+1} \mid k-\ell d_1^{(r)}$, autrement dit que $\ell \equiv \overline{d_1^{(r)}}k \Mod{d'_{r+1}}$ et puisque
$\alpha_{r+1}=\frac{\ell d_1^{(r)}-k}{d'_{r+1}}$, on a également
$$
\frac{k-Xd'_{r+1}}{d_1^{(r)}} \leqslant \ell \leqslant \frac{k+Xd'_{r+1}}{d_1^{(r)}}.
$$
Réciproquement, les relations (\ref{rela}) et (\ref{dprime}) fournissent les égalités
$$
\forall \ell \in \mathbb{Z}, \quad  \sum_{i=r+2}^n d_i\alpha_i =d_{1,r} \ell +d_{1,r+1}\left(k-\ell d_1^{(r)}\right),
$$
si bien que pour $\ell$ vérifiant
\begin{eqnarray}
\ell \equiv \overline{d_1^{(r)}}k \Mod{d'_{r+1}} \quad \mbox{et} \quad\frac{k-Xd'_{r+1}}{d_1^{(r)}} \leqslant \ell \leqslant \frac{k+Xd'_{r+1}}{d_1^{(r)}},
\label{auxx}
\end{eqnarray}
on obtient que $(\alpha_{r+2},\dots,\alpha_n)$ provient de $(\alpha_{r+1},\dots,\alpha_n) \in \mathcal{A}_r(\mathbf{y},X)$ avec $\alpha_{r+1}=\frac{\ell d_1^{(r)}-k}{d'_{r+1}}.$
Ainsi, on en déduit bien la formule (\ref{rec}). Or, par définition
$$
\mathcal{A}_{n-1}(\mathbf{y},X)=\left\{
\alpha_n \in \mathbb{Z}  \quad :  \quad
|\alpha_n|\leqslant X,\quad
d_n\alpha_n \equiv 0 \Mod{d_{1,n-1}}
\right\}.
$$
On constate que $d_{1,n-1}=d_1^{(n-1)}$ et que $\mbox{pgcd}\left(d_{1,n-1},d_n\right)=1$ de sorte que
$$
\#\mathcal{A}_{n-1}(\mathbf{y},X)=\frac{2X}{d_1^{(n-1)}}+O(1).
$$
Le lemme suit alors par récurrence.
\\
\hfill
$\square$



\begin{lemme}
Soit $n \geqslant 3$. Pour tout $X \geqslant 1$ et $\mathbf{y} \in \mathbb{N}^n$, on a
$$
\#\mathcal{A}(\mathbf{y};X)=\frac{X^{n-1}}{d_1}\left( b(\mathbf{y})+O\left( \frac{1}{X}\overset{n-1}{\underset{j=2}{\sum}}d_1^{(j-1)} \right)  \right)
$$
avec les notations (\ref{dprime}) et où
\begin{eqnarray}
 b(\mathbf{y})={\rm{vol}}\left\{(\alpha_2,\dots,\alpha_n) \in [-1,1]^{n-1} \quad : \quad  \left|\sum_{i=2}^n d_i\alpha_i\right| \leqslant d_1\right\}.
\label{by}
\end{eqnarray}
\label{Lemme}
\end{lemme}
\noindent
\textit{Démonstration--} La preuve proposée ne repose pas sur des arguments de géométrie des nombres mais suit les grandes lignes de la preuve du lemme 2.3 de \cite{Br}. On a l'égalité
$$
\#\mathcal{A}(\mathbf{y};X)=\#\left\{
(\alpha_2,\dots,\alpha_n) \in \mathbb{Z}^{n-1}  \quad :  \quad
\max_{2\leqslant i \leqslant n}|\alpha_i|\leqslant X,\quad
\begin{array}{c}
\overset{n}{\underset{i=2}{\sum}} d_i\alpha_i \equiv 0 \Mod{d_1},\\[3mm]
\left|\overset{n}{\underset{i=2}{\sum}}  d_i\alpha_i\right| \leqslant  d_1 X
\end{array}
\right\}.
$$
En raisonnant comme pour établir la formule (\ref{auxx}) de la preuve du Lemme \ref{lemme1}, on obtient alors
\begin{eqnarray}
\#\mathcal{A}(\mathbf{y};X)=\!\!\!\!\!\!
\sum_{\boldsymbol{\alpha} \in \mathcal{A}_2(\mathbf{y},X)} \!\!\!\#\left\{ \frac{k_2-Xd'_2}{d'_1}\leqslant \ell \leqslant \frac{k_2+Xd'_2}{d'_1}  \hspace{1mm} : \hspace{-1mm} 
\begin{array}{c}
\ell \equiv -\overline{d'_1}k_2 \Mod{d'_2},\\[2mm]
\left|d_{1,2}\left(\ell d'_1-k_2\right)+\overset{n}{\underset{i=3}{\sum}} d_i\alpha_i\right| \leqslant  d_1 X
\end{array}
 \right\},
\label{auxxx}
\end{eqnarray}
avec $k_2=\overset{n}{\underset{j=3}{\sum}} d_j \alpha_j/d_{1,2}$, où $d_{1,2}$ a été défini en (\ref{dprime}).
Or, on a
$$
\begin{aligned}
&\#\left\{ \frac{k_2-Xd'_2}{d'_1}\leqslant \ell \leqslant \frac{k_2+Xd'_2}{d'_1}  \hspace{1mm} : \hspace{1mm}
\begin{array}{c}
\ell \equiv -\overline{d'_2}k_2 \Mod{d'_2},\\[4mm]
\left|d_{1,2}\left(\ell d'_1-k_2\right)+\overset{n}{\underset{i=3}{\sum}} d_i\alpha_i\right| \leqslant d_1 X
\end{array}
 \right\}\\[2mm]
 & \qquad \qquad  =\frac{1}{d'_2}{\rm{vol}}\left\{t \in \left[\frac{k_2-Xd'_2}{d'_1},\frac{k_2+Xd'_2}{d'_1}\right] \hspace{1mm} :\hspace{1mm} \left|d_{1,2}\left(t d'_1-k_2\right)+\sum_{i=3}^n d_i\alpha_i\right| \leqslant d_1 X\right\}+O(1).\\
 \end{aligned}
$$
Le changement de variable $\alpha_2=\frac{td'_1-k_2}{d'_2}$ fournit alors
$$
\#\mathcal{A}(\mathbf{y};X)=\frac{1}{d'_1}\sum_{\boldsymbol{\alpha} \in \mathcal{A}_2(\mathbf{y},X)} {\rm{vol}}\left\{\alpha_2 \in [-X,X] \quad : \quad  \left|\sum_{i=2}^n d_i\alpha_i\right| \leqslant d_1 X\right\}+O\left(\#\mathcal{A}_2(\mathbf{y},X)\right).
$$
D'après le Lemme \ref{lemme1}, on aboutit à la formule
$$
\#\mathcal{A}(\mathbf{y};X)=\frac{1}{d'_1}\sum_{\boldsymbol{\alpha} \in \mathcal{A}_2(\mathbf{y},X)} {\rm{vol}}\left\{\alpha_2 \in [-X,X] \quad : \quad  \left|\sum_{i=2}^n d_i\alpha_i\right| \leqslant d_1 X\right\}+O\left(\prod_{j=3}^n\frac{X}{ d_1^{(j-1)}}\right).
$$
On pose
$$
S_r(\mathbf{y};X)=\sum_{\boldsymbol{\alpha} \in \mathcal{A}_r(\mathbf{y},X)} {\rm{vol}}\left\{(\alpha_{2},\dots,\alpha_r) \in [-X,X]^{r-1} \quad : \quad  \left|\sum_{i=2}^n d_i\alpha_i\right| \leqslant d_1 X\right\}
$$
Avec les notations (\ref{dprime}) et pour tout $2 \leqslant r \leqslant n-1$, on montre alors l'estimation
\begin{eqnarray}
S_r(\mathbf{y};X)=\frac{1}{d_1^{(r-1)}}S_{r+1}(\mathbf{y};X)+O\left( X^{r-1}\prod_{j=r+2}^n\frac{X}{ d_1^{(j-1)}} \right).
\label{sr}
\end{eqnarray}
En effet, en raisonnant de manière analogue à (\ref{auxxx}), on obtient 
$$
S_r(\mathbf{y};X)=\!\!\!\!\!\!\sum_{\substack{\scriptscriptstyle\boldsymbol{\alpha} \in \mathcal{A}_{r+1}(\mathbf{y},X)\\[1mm]\scriptscriptstyle \frac{k_{r+1}-X d'_{r+1}}{d_1^{(r)}}\leqslant \ell  \leqslant \frac{k_{r+1}+X d'_{r+1}}{d_1^{(r)}} \\[1mm] \scriptscriptstyle\ell \equiv -\overline{d_1^{(r)}}k_{r+1} \Mod{d'_{r+1}}}} 
\!\!\!\!\!\!{\rm{vol}}\left\{\boldsymbol{\alpha} \in [-X,X]^{r-1} \hspace{1mm} : \hspace{1mm}  \left|d_{1,r}(\ell d'_{r+1}-k_{r+1})+\sum_{i=2\atop i \neq r+1}^n d_i\alpha_i\right| \leqslant d_1 X\right\}
$$
avec $k_{r+1}=\overset{n}{\underset{j=r+2}{\sum}} d_j \alpha_j/d_{1,r+1}.$ On écrit alors
$$
\begin{aligned}
& {\rm{vol}}\left\{\boldsymbol{\alpha} \in [-X,X]^{r-1} \hspace{1mm} : \hspace{1mm}  \left|d_{1,r}(\ell d'_{r+1}-k_{r+1})+\sum_{i=2\atop i \neq r+1}^n d_i\alpha_i\right| \leqslant d_1 X\right\}\\
&\qquad \qquad \qquad \qquad \qquad \qquad \qquad \quad=\int_{\boldsymbol{\alpha} \in [-X,X]^{r-1}} \mathds{1}_{\Big|(\ell d'_{r+1}-k_{r+1})+\overset{n}{\underset{i=2\atop i \neq r+1}{\sum}} d_i\alpha_i\Big|  \leqslant d_1 X}(\boldsymbol{\alpha})\mbox{d}\alpha_2\cdots\mbox{d}\alpha_r
\end{aligned}
$$
et on obtient ensuite en intervertissant la sommation sur $\ell$ avec l'intégrale
$$
\begin{aligned}
&\sum_{\frac{k_{r+1}-Xd'_{r+1}}{d_1^{(r)}}\leqslant \ell \leqslant \frac{k_{r+1}+X d'_{r+1}}{d_1^{(r)}}\atop \ell \equiv -\overline{d_1^{(r)}}k_{r+1} \Mod{d'_{r+1}}} {\rm{vol}}\left\{\boldsymbol{\alpha} \in [-X,X]^{r-1} \quad : \quad  \left|d_{1,r}(\ell d'_{r+1}-k_{r+1})+\sum_{i=2\atop i \neq 3}^n d_i\alpha_i\right| \leqslant d_1 X\right\}\\[2mm]
&\qquad \qquad =\int_{\boldsymbol{\alpha} \in [-X,X]^{r-1}} \sum_{\frac{k_{r+1}-Xd'_{r+1}}{d_1^{(r)}}\leqslant \ell \leqslant \frac{k_{r+1}+X d'_{r+1}}{d_1^{(r)}}\atop \ell \equiv -\overline{d_1^{(r)}}k_{r+1} \Mod{d'_{r+1}}}\mathds{1}_{\Big|(\ell d'_{r+1}-k_{r+1})+\overset{n}{\underset{i=2\atop i \neq r+1}{\sum}} d_i\alpha_i\Big|  \leqslant d_1 X}(\boldsymbol{\alpha})\mbox{d}\alpha_2\cdots\mbox{d}\alpha_r.\\
\end{aligned}
$$
Le raisonnement précédent fournit alors
$$
\begin{aligned}
&\sum_{\frac{k_{r+1}-Xd'_{r+1}}{d_1^{(r)}}\leqslant \ell \leqslant \frac{k_{r+1}+X d'_{r+1}}{d_1^{(r)}}\atop \ell \equiv -\overline{d_1^{(r)}}k_{r+1} \Mod{d'_{r+1}}}\mathds{1}_{\Big|(\ell d'_{r+1}-k_{r+1})+\overset{n}{\underset{i=2\atop i \neq r+1}{\sum}} d_i\alpha_i\Big|  \leqslant d_1 X}(\boldsymbol{\alpha})\\[1mm]
&\qquad \qquad \qquad \qquad\qquad \qquad \qquad \qquad \qquad =\frac{1}{d_1^{(r)}}\int_{\alpha_{r+1} \in [-X,X]} \mathds{1}_{\left|\overset{n}{\underset{i=2}{\sum}} d_i\alpha_i\right| \leqslant X|d_1|}\mbox{d}\alpha_{r+1}+O(1).
\end{aligned}
$$
Ainsi en utilisant le Lemme \ref{lemme1}, on obtient bien (\ref{sr}).
Finalement, en itérant (\ref{sr}), il vient
$$
\#\mathcal{A}(\mathbf{y};X)=\frac{1}{d_1} {\rm{vol}}\left\{\boldsymbol{\alpha} \in [-X,X]^n \quad : \quad  \left|\sum_{i=2}^n d_i\alpha_i\right| \leqslant d_1 X\right\}+O\left(\sum_{r=2}^n\prod_{j=2 \atop j \neq r}^n\frac{X}{ d_1^{(j-1)}}\right),
$$
puisqu'on a l'égalité
\begin{eqnarray}
d_1=\prod_{j=2}^n d_1^{(j-1)}.
\label{d1}
\end{eqnarray}
Le résultat en découle après le changement de variables $\boldsymbol{\alpha}'=\frac{\boldsymbol{\alpha}}{X}$ et grâce au fait que, d'après (\ref{d1}), l'on ait
$$
\forall r \in \llbracket 2,n \rrbracket, \quad \frac{d_1}{\overset{n}{\underset{j=2 \atop j \neq r}{\prod}}d_1^{(j-1)}}=d_1^{(r-1)}.
$$
\hfill
$\square$

\subsubsection{Une borne supérieure}

\begin{lemme}
Lorsque $1 \leqslant y_i \leqslant X$ pour tout $i \in \llbracket 1,n \rrbracket$, on a l'estimation
$$
N_{\mathbf{y}}\left(X\right) \ll \frac{X^{n-1}}{||{\mathbf{d}}||_2}
$$
où $N_{\mathbf{y}}\left(X\right)$ a été défini en (\ref{Ny}) et $||.||_2$ désigne la norme euclidienne usuelle sur $\mathbb{R}^n$.
\label{lemme2}
\end{lemme}
\noindent
\textit{Démonstration--}
Soit $X \geqslant 1$. Par symétrie, on peut supposer que $d_1 \geqslant d_2 \geqslant \cdots \geqslant d_n$ si bien que
$$
N_{\mathbf{y}}\left(X\right) \leqslant \#\left\{
(\alpha_2,\dots,\alpha_n) \in \mathbb{Z}^{n-1}  \quad :  \quad
\max_{2\leqslant i \leqslant n}|\alpha_i|\leqslant X,\quad
\sum_{i=2}^n d_i\alpha_i \equiv 0 \Mod{d_1}
\right\}.
$$
En utilisant le Lemme \ref{lemme1} et (\ref{d1}), on obtient bien
\begin{equation}
N_{\mathbf{y}}(X) \ll \frac{X^{n-1}}{d_1} \ll  \frac{X^{n-1}}{||\textbf{d}||_2}. \tag*{$\square$}
\end{equation}
\noindent
\textbf{Remarque.--} Comme mentionné dans \cite{Br}, il est ici important de voir que $||\textbf{d}||_2$ peut être plus grand que $B$ et que l'on fait donc mieux que l'estimation triviale en $O\left(\frac{B^{n-1}}{||\textbf{d}||_2}+B^{n-2}\right)$. Cette estimation découle également de résultats de géométrie des nombres.\\
\newline
\indent
Quitte à appliquer une permutation de $\llbracket 1,n \rrbracket$, il vient que
\begin{eqnarray}
N(B)=n!N_1(B)+R(B),
\label{r}
\end{eqnarray}
avec
\begin{eqnarray}
N_1(B)=\sum_{\mathbf{y} \in \mathbb{N}^n \atop 1 \leqslant y_1\leqslant \cdots\leqslant y_n\leqslant \sqrt[n]{B}} N_{\mathbf{y}}\left(B^{1/n}\right)
\label{n1}
\end{eqnarray}
et
$$
R(B)= \sum_{\substack{\scriptscriptstyle \mathbf{y} \in \mathbb{N}^n \\\scriptscriptstyle 1 \leqslant y_i\leqslant \sqrt[n]{B}\\\scriptscriptstyle \exists i \neq j, \hspace{1mm} y_i=y_j}} N_{\mathbf{y}}\left(B^{1/n}\right)-n! \sum_{\substack{\scriptscriptstyle\mathbf{y} \in \mathbb{N}^n \\\scriptscriptstyle 1 \leqslant y_1\leqslant\cdots\leqslant y_n \leqslant \sqrt[n]{B}\\\scriptscriptstyle \exists i \neq j, \hspace{1mm} y_i=y_j}} N_{\mathbf{y}}\left(B^{1/n}\right).
$$
Le lemme suivant fournit une borne supérieure du bon ordre de grandeur pour la quantité $N_1(B)$.

\begin{lemme}
Lorsque $B \geqslant 2$, on a
$
N_1\left(B\right) \ll B(\log B)^{2^n-n-1}.
$
\label{lemme3}
\end{lemme}
\noindent
\textit{Démonstration--}
Lorsque $B \geqslant 2$, on a 
$$
N_1\left(B\right) \leqslant  \sum_{\mathbf{y} \in \mathbb{N}^n \atop 1 \leqslant y_1\leqslant \cdots\leqslant y_n\leqslant \sqrt[n]{B}} N_{\mathbf{y}}\left(B^{1/n}\right) \ll B^{1-\frac{1}{n}}\sum_{\mathbf{y} \in \mathbb{N}^n \atop 1 \leqslant y_1\leqslant \cdots\leqslant y_n\leqslant \sqrt[n]{B}} \frac{1}{d_1}
$$
d'après le Lemme \ref{lemme2}. On traduit à présent les conditions sur $\mathbf{y}$ en termes de conditions sur $\mathbf{z}$. Lorsque $1 \leqslant y_1\leqslant \cdots\leqslant y_n\leqslant \sqrt[n]{B}$, on a en particulier $d_j \geqslant d_{j+1}$ pour $2 \leqslant j \leqslant n-1$. Pour $j \in \llbracket 2,n-1 \rrbracket$, cette dernière condition se réécrit
\begin{eqnarray}
\prod_{\varepsilon_j(h)=1 \atop\varepsilon_{j+1}(h)=0} z_h \leqslant \prod_{\varepsilon_j(h)=0 \atop\varepsilon_{j+1}(h)=1} z_h.
\label{cond1}
\end{eqnarray}
On constate alors en posant $H_0=\left\{ h_{\ell}=\overset{\ell}{\underset{k=1}{\sum}} 2^{k-1} : 2\leqslant \ell \leqslant n-1\right\}$, que dans l'inégalité (\ref{cond1}), seule la variable $h_j$ de $H_0$ apparaît et que (\ref{cond1}) se réécrit par conséquent
$$
z_{h_j} \leqslant Z_{h_j}:=\prod_{\substack{\scriptscriptstyle\varepsilon_j(h)=0 \\ \scriptscriptstyle\varepsilon_{j+1}(h)=1}} z_h \prod_{\substack{\scriptscriptstyle\varepsilon_j(h)=1 \\ \scriptscriptstyle\varepsilon_{j+1}(h)=0 \\\scriptscriptstyle h \neq h_j}} z_h^{-1}.
$$
La condition $d_1 \geqslant d_2$ est équivalente à
\begin{eqnarray}
\prod_{\varepsilon_2(h)=1 \atop\varepsilon_{1}(h)=0} z_h \geqslant \prod_{\varepsilon_2(h)=0 \atop\varepsilon_1(h)=1} z_h.
\label{cond2}
\end{eqnarray}
En notant $H_1:=\{5\}$ si $n \geqslant 4$ et $H_1:=\{1\}$ si $n=3$, elle peut se réécrire
$$
z_{h_1}:=z_5 \leqslant Z_{h_1}:=\prod_{\substack{\scriptscriptstyle\varepsilon_1(h)=0 \\\scriptscriptstyle \varepsilon_{2}(h)=1}} z_h \prod_{\substack{\scriptscriptstyle \varepsilon_1(h)=1 \\\scriptscriptstyle \varepsilon_{2}(h)=0 \\\scriptscriptstyle h \neq 5}} z_h^{-1}
$$
dans le cas $n \geqslant 4$ et
$$
z_{h_1}:=z_{1} \leqslant Z_{h_1}:=\frac{z_2z_6}{z_5} 
$$
lorsque $n=3$. Enfin, notant $H_2=\left\{h_n=\overset{n}{\underset{k=1}{\sum}} 2^{k-1}\right\}$, on remarque que la seule condition faisant apparaître $z_{h_n}$ est $y_n \leqslant \sqrt[n]{B}$. On peut réécrire cette condition sous la forme
$$
z_{h_n} \leqslant Z_{h_n} \quad \mbox{avec} \quad Z_{h_n}:=\sqrt[n]{B}\prod_{\varepsilon_n(h)=1\atop h \neq h_n} z_h^{-1}.
$$
Par conséquent, en négligeant les conditions de coprimalité provenant du fait que $\mathbf{z}$ est réduit, on obtient
$$
N_1\left(B\right) \ll B^{1-\frac{1}{n}}\sum_{\substack{\scriptscriptstyle z_h \leqslant B^{1/n}\\\scriptscriptstyle d_n \leqslant \cdots \leqslant d_1 \\ \scriptscriptstyle z_{h_n} \leqslant Z_{h_n} }} \frac{1}{d_1}.
$$
La contribution des $z_h$ pour $h \in H_0$ est alors majorée par
$$
\ll \frac{B^{1-\frac{1}{n}}}{d_1}\prod_{2\leqslant j \leqslant n-1} Z_{h_j}=\frac{B^{1-\frac{1}{n}}}{d_1}\prod_{\substack{\scriptscriptstyle{\varepsilon_2(h)=0}\\ \scriptscriptstyle \varepsilon_n(h)=1}}z_h \prod_{\substack{\scriptscriptstyle \varepsilon_2(h)=1\\\scriptscriptstyle \varepsilon_n(h)=0\\\scriptscriptstyle h \not \in H_0}}z_h^{-1}.
$$
En effet, pour $j \in \llbracket 2,n-2 \rrbracket$, on a 
$$
\begin{aligned}
Z_{h_j}Z_{h_{j+1}}&=\frac{\displaystyle\underset{\substack{\scriptscriptstyle\varepsilon_{\scriptscriptstyle j}(h)=0\\ \scriptscriptstyle\varepsilon_{j+1}(h)=1}}{\prod}z_h\underset{\substack{\scriptscriptstyle\varepsilon_{j+1}(h)=0\\ \scriptscriptstyle\varepsilon_{j+2}(h)=1}}{\prod}z_h}{\displaystyle\underset{\substack{\scriptscriptstyle\varepsilon_j(h)=1\\ \scriptscriptstyle\varepsilon_{j+1}(h)=0}}{\prod}z_h\underset{\substack{\scriptscriptstyle\varepsilon_{j+1}(h)=1\\ \scriptscriptstyle\varepsilon_{j+2}(h)=0}}{\prod}z_h}=\frac{\displaystyle\underset{\substack{\scriptscriptstyle\varepsilon_j(h)=0\\\scriptscriptstyle \varepsilon_{j+1}(h)=1\\\scriptscriptstyle \varepsilon_{j+2}(h)=1}}{\prod}z_h\underset{\substack{\scriptscriptstyle\varepsilon_{j+1}(h)=0\\ \scriptscriptstyle\varepsilon_{j+2}(h)=1}}{\prod}z_h}{\displaystyle\underset{\substack{\scriptscriptstyle\varepsilon_j(h)=1\\ \scriptscriptstyle\varepsilon_{j+1}(h)=0}}{\prod}z_h\underset{\substack{\scriptscriptstyle\varepsilon_{j+1}(h)=1\\ \scriptscriptstyle\varepsilon_{j+2}(h)=0\\ \scriptscriptstyle\varepsilon_j(h)=1}}{\prod}z_h}=\frac{\displaystyle\underset{\substack{\scriptscriptstyle\varepsilon_j(h)=0\\ \scriptscriptstyle\varepsilon_{j+1}(h)=1\\ \scriptscriptstyle\varepsilon_{j+2}(h)=1}}{\prod}z_h\underset{\substack{\scriptscriptstyle\varepsilon_{j+1}(h)=0\\ \scriptscriptstyle\varepsilon_{j+2}(h)=1\\ \scriptscriptstyle\varepsilon_j(h)=0}}{\prod}z_h}{\displaystyle\underset{\substack{\scriptscriptstyle\varepsilon_j(h)=1\\ \scriptscriptstyle\varepsilon_{j+1}(h)=0\\ \scriptscriptstyle\varepsilon_{j+2}(h)=0}}{\prod}z_h\underset{\substack{\scriptscriptstyle\varepsilon_{j+1}(h)=1\\ \scriptscriptstyle\varepsilon_{j+2}(h)=0\\ \scriptscriptstyle\varepsilon_j(h)=1}}{\prod}z_h}\\[2mm]
&=\underset{\substack{\scriptscriptstyle\varepsilon_j(h)=0\\ \scriptscriptstyle\varepsilon_{j+2}(h)=1}}{\prod}z_h \underset{\substack{\scriptscriptstyle\varepsilon_j(h)=1\\ \scriptscriptstyle\varepsilon_{j+2}(h)=0\\ \scriptscriptstyle h \not \in \{{\scriptscriptstyle h_j},{\scriptscriptstyle h_{\scriptscriptstyle j+1}}\}}}{\prod}z_h^{-1}\\
\end{aligned}
$$
et il suffit d'itérer ce calcul pour obtenir l'expression
$$
\prod_{2\leqslant j \leqslant n-1} Z_{h_j}=\prod_{\substack{\scriptscriptstyle \varepsilon_2(h)=0\\\scriptscriptstyle \varepsilon_n(h)=1}}z_h \prod_{\substack{\scriptscriptstyle \varepsilon_2(h)=1\\\scriptscriptstyle \varepsilon_n(h)=0\\\scriptscriptstyle h \not \in H_0}}z_h^{-1}.
$$
Remplaçant $d_1$ par son expression,
il vient une contribution des $z_h$ avec $h \in H_0$
\begin{eqnarray}
\ll B^{1-\frac{1}{n}}\prod_{\varepsilon_1(h)=0}z_h^{-1} \prod_{\substack{\scriptscriptstyle\varepsilon_2(h)=0\\\scriptscriptstyle \varepsilon_n(h)=1}}z_h \prod_{\substack{\scriptscriptstyle \varepsilon_2(h)=1\\ \scriptscriptstyle\varepsilon_n(h)=0\\\scriptscriptstyle h \not \in H_0}}z_h^{-1}.
\label{z5}
\end{eqnarray}
On remarque que $z_5$ intervient dans $Z_3$ et $Z_7$ mais disparaît dans le produit des $Z_{h_j}$ et que $z_5$ n'intervient pas dans l'expression de $d_1$ lorsque $n \geqslant 4$. Sommant alors sur $z_{h_1}$, on aboutit à une contribution des $z_h$ pour $h \in H_0 \cup H_1$
$$
\begin{aligned}
\ll& B^{1-\frac{1}{n}}\prod_{\varepsilon_1(h)=0}z_h^{-1} \prod_{\substack{\scriptscriptstyle\varepsilon_2(h)=0\\\scriptscriptstyle \varepsilon_n(h)=1}}z_h \prod_{\substack{\scriptscriptstyle\varepsilon_2(h)=1\\ \scriptscriptstyle\varepsilon_n(h)=0\\\scriptscriptstyle h \not \in H_0}}z_h^{-1}\prod_{\substack{\scriptscriptstyle\varepsilon_1(h)=0 \\ \scriptscriptstyle\varepsilon_{2}(h)=1}} z_h \prod_{\substack{\scriptscriptstyle\varepsilon_1(h)=1 \\\scriptscriptstyle \varepsilon_{2}(h)=0 \\\scriptscriptstyle h \neq 5}} z_h^{-1}\\
\ll & B^{1-\frac{1}{n}}\prod_{\substack{\scriptscriptstyle\varepsilon_1(h)=0\\\scriptscriptstyle \varepsilon_2(h)=0}}z_h^{-1} \prod_{\substack{\scriptscriptstyle\varepsilon_2(h)=0\\\scriptscriptstyle \varepsilon_n(h)=1 \\\scriptscriptstyle \varepsilon_1(h)=0}}z_h \prod_{\substack{\scriptscriptstyle\varepsilon_2(h)=1\\\scriptscriptstyle \varepsilon_n(h)=0\\\scriptscriptstyle h \not \in H_0}}z_h^{-1}
 \prod_{\substack{\scriptscriptstyle\varepsilon_1(h)=1 \\\scriptscriptstyle \varepsilon_{2}(h)=\varepsilon_{n}(h)=0 \\\scriptscriptstyle h \neq 5}} z_h^{-1}\\
 \ll&B^{1-\frac{1}{n}}\prod_{\substack{\scriptscriptstyle \varepsilon_1(h)=0\\ \scriptscriptstyle\varepsilon_2(h)=0\\\scriptscriptstyle \varepsilon_n(h)=0}}z_h^{-1}  \prod_{\substack{\scriptscriptstyle\varepsilon_2(h)=1\\ \scriptscriptstyle\varepsilon_n(h)=0\\\scriptscriptstyle h \not \in H_0}}z_h^{-1}
 \prod_{\substack{\scriptscriptstyle\varepsilon_1(h)=1 \\\scriptscriptstyle \varepsilon_{2}(h)=\varepsilon_{n}(h)=0 \\\scriptscriptstyle h \neq 5}} z_h^{-1}=B^{1-\frac{1}{n}}\prod_{\substack{\scriptscriptstyle \varepsilon_n(h)=0\\\scriptscriptstyle h \not \in H_0 \cup H_1}}z_h^{-1}. \\
 \end{aligned}
 $$ 
 Une dernière sommation sur la variable $z_{h_n}$ fournit de même une contribution des $z_h$ avec $h \in H_0 \cup H_1 \cup H_2$
 \begin{eqnarray}
 \ll B \prod_{\substack{  h \not \in H_0 \cup H_1 \cup H_2}}z_h^{-1}. 
 \label{ineg}
 \end{eqnarray}
 Si $n=3$, c'est la variable $z_{h_1}=z_{1}$ qui n'intervient pas et on somme alors sur $z_{1}$ puis sur $z_{h_3}=z_7$ pour obtenir
 $$
 \ll B\sum_{z_2,z_4,z_5,z_6 \leqslant \sqrt[3]{B}} \frac{1}{z_2z_4z_5z_6}.
 $$
 Dans tous les cas, il reste alors $2^n-n-1$ variables $z_h \leqslant \sqrt[n]{B}$ à sommer. En effet, si $1 \leqslant y_i \leqslant \sqrt[n]{B}$ pour tout $i \in \llbracket 1,n \rrbracket$, alors on a également $1 \leqslant z_j \leqslant \sqrt[n]{B}$ pour tout $j \in \llbracket 1,N \rrbracket$ si bien qu'on obtient finalement
 $$
 N_1\left(B\right) \ll B(\log B)^{2^n-n-1}.
 $$
\hfill
$\square$

On utilise alors ce résultat pour démontrer que la quantité $R(B)$ apparaissant dans (\ref{r}) est bien un terme d'erreur.

\begin{lemme}
Pour tout $B \geqslant 2$, on a
$
R(B) \ll B\left( \log B\right)^{2^{n-1}-n},
$
où $R(B)$ a été défini en (\ref{r}).
\label{erreur}
\end{lemme}
\noindent
\textit{Démonstration--}
Quitte à réordonner, on obtient
$$
R(B) \ll \sum_{\substack{\scriptscriptstyle\mathbf{y} \in \mathbb{N}^n \\\scriptscriptstyle 1 \leqslant y_1 \leqslant \cdots\leqslant y_n\leqslant \sqrt[n]{B}\\\scriptscriptstyle \exists i, \hspace{1mm} y_i=y_{i+1}}} N_{\mathbf{y}}\left(B^{1/n}\right).
$$
De plus, pour $i \in \llbracket 1,n-1 \rrbracket$, la condition $y_i=y_{i+1}$ se réécrit
\begin{eqnarray}
\prod_{\substack{\scriptscriptstyle\varepsilon_i(h)=1 \\\scriptscriptstyle \varepsilon_{i+1}(h)=0 
}} z_h=\prod_{\substack{\scriptscriptstyle\varepsilon_i(h)=0 \\\scriptscriptstyle \varepsilon_{i+1}(h)=1}} z_h.
\label{cond}
\end{eqnarray}
Puisque deux entiers $h$ et $\ell$ de $\llbracket 1,N \rrbracket$ tels que $\varepsilon_i(h)= \varepsilon_{i+1}(\ell)=1$ et $\varepsilon_{i+1}(h)= \varepsilon_{i}(\ell)=0$ ne sont pas comparables pour la relation d'ordre $\preceq$ introduite en section 2, on en déduit que $\mbox{pgcd}(z_h,z_{\ell})=1$. Il s'ensuit que 
\begin{eqnarray}
\forall h \in \llbracket 1,N \rrbracket \quad \mbox{tel que} \quad \varepsilon_i(h)+\varepsilon_{i+1}(h)=1, \quad z_h=1.
\label{un}
\end{eqnarray}
En particulier, on a $Z_{h_i}=z_{h_i}=1$.\\
\indent
Supposons alors que $y_i=y_{i+1}$ pour un certain $i \in \llbracket 1,n-1\rrbracket$. Une sommation sur les $z_h$ pour $h \in H_0$ suivie d'une sommation sur $z_{h_1}$ et $z_{h_n}$ fournit une contribution
$$
\ll B \sum_{\substack{\scriptscriptstyle z_h \leqslant \sqrt[n]{B} \\\scriptscriptstyle h \not \in H_0 \cup H_1 \cup H_2\\\scriptscriptstyle \varepsilon_i(h)+\varepsilon_{i+1}(h)\neq1}} \underset{\substack{\scriptscriptstyle h \not \in H_0 \cup H_1 \cup H_2\\\scriptscriptstyle \varepsilon_i(h)+\varepsilon_{i+1}(h)\neq1}}{\prod}z_h^{-1} \ll B\left( \log B\right)^{2^{n-1}-n},
$$
au vu de (\ref{ineg}) et (\ref{un}).
\hfill
$\square$

\subsubsection{Démonstration de l'asymptotique}

Estimons désormais le cardinal $N_1(B)$ défini en (\ref{n1}). D'après le Lemme \ref{lemme3}, on a
$$
N_1(B) \ll B(\log B)^{2^n-n-1}.
$$
Pour $A >0$ fixé, on peut supposer que
$$
z_h>\log(B)^A \quad \mbox{pour} \quad h \not \in H_0 \cup H_1 \cup H_2
$$
avec les notations de la section précédente. En effet, en reprenant l'inégalité (\ref{ineg}), on obtient une contribution complémentaire (c'est-à-dire pour laquelle au moins un $z_h \leqslant \log(B)^A$ pour $h \not \in H_0 \cup H_1 \cup H_2$) majorée par
$$
\ll B \log (B) ^{2^n-n-2}\log(\log B).
$$
Cela est suffisant pour donner lieu à un terme d'erreur en vue du Théorème \ref{theor1}. Cette réduction du domaine de comptage permet de contrôler le terme d'erreur du Lemme \ref{Lemme} de la façon suivante
$$
\forall j \in \llbracket 2,n-1 \rrbracket, \quad \frac{ d_1^{(j-1)}}{B^{1/n}} \leqslant \frac{d_1^{(j-1)}}{y_j}=\prod_{\varepsilon_1(h)+\cdots+\varepsilon_{j-1}(h)\neq 0\atop \varepsilon_j(h)=1} z_h^{-1} \leqslant \log(B)^{-A}
$$
puisqu'on a $y_j \leqslant \sqrt[n]{B}$ et que tous les indices des variables intervenant dans le produit ci-dessus ne sont pas dans $ H_0 \cup H_1 \cup H_2$.
Par le Lemme \ref{Lemme}, il s'ensuit l'estimation
$$
N_{\mathbf{y}}\left(B^{1/n}\right)=\frac{B^{1-1/n}}{ d_1}b(\mathbf{y})+O\left( \frac{B^{1-1/n}}{d_1}\log(B)^{-A} \right)
$$
puis grâce au Lemme \ref{lemme3}
$$
N_1(B)=B^{1-1/n}\sum_{\mathbf{y} \in \mathbb{N}^n \atop 1 \leqslant y_1\leqslant \cdots\leqslant y_n\leqslant \sqrt[n]{B}}\frac{b(\mathbf{y})}{d_1}+O\left(  B (\log B)^{2^n-n-2}\log(\log B) \right).
$$
On a ici remplacé la somme
$$
\sum_{\substack{\scriptscriptstyle \mathbf{y} \in \mathbb{N}^n \\\scriptscriptstyle 1 \leqslant y_1\leqslant \cdots\leqslant y_n\leqslant \sqrt[n]{B}\\\scriptscriptstyle  z_h>\log(B)^A, \hspace{0.5mm} h \not \in H_0 \cup H_1 \cup H_2}}\frac{b(\mathbf{y})}{d_1}
$$
par la somme 
$$
\sum_{\mathbf{y} \in \mathbb{N}^n \atop 1 \leqslant y_1\leqslant \cdots\leqslant y_n\leqslant \sqrt[n]{B}}\frac{b(\mathbf{y})}{d_1}
$$
au prix d'une contribution négligeable en raisonnant à nouveau comme dans la preuve du Lemme~\ref{lemme3} puisque $b(\mathbf{y}) \ll 1$. On effectue alors la sommation dans le même ordre que lors de la preuve du Lemme \ref{lemme3}. 
Lorsque $n \geqslant 4$, les conditions $Z_{h_j} \geqslant 1$ ne font intervenir la variable $z_5$ que dans $Z_3$ et $Z_7$. Comme on souhaite sommer sur les $z_h$ avec $h \in H_0$ puis sur $z_5$ puis sur $z_{h_n}$, on restreint le domaine de comptage de façon à ne plus avoir cette dépendance en $z_5$ dans $Z_3$ et $Z_7$. Posant 
$$
Z'_5=\prod_{\substack{\scriptscriptstyle\varepsilon_2(h)=1 \\\scriptscriptstyle \varepsilon_3(h)=0 \\ h\neq 3}} z_h \prod_{\substack{\scriptscriptstyle\varepsilon_3(h)=1 \\\scriptscriptstyle \varepsilon_2(h)=0 \\ h\neq 5}} z_h^{-1},
$$
la condition $Z_3=z_5/Z'_5 \geqslant 1$ se réécrit $z_5 \geqslant Z'_5$. Puisque $z_5 \leqslant Z_{h_1}$, on a $Z'_5 \leqslant Z_{h_1}$ d'une part. D'autre part, la contribution des $\mathbf{z}$ tels que $z_5 < Z'_5$ est négligeable. En effet, de la même manière que lors de la preuve du Lemme \ref{lemme2}, on montre que ces $\mathbf{z}$ contribuent pour
\begin{eqnarray}
\ll  B \sum_{z_h \leqslant \sqrt[n]{B} \atop h \not \in H_0 \cup H_1 \cup H_2}\hspace{1mm} \prod_{\substack{ h \not \in H_0 \cup H_1 \cup H_2}}z_h^{-1} \times \frac{Z'_5}{Z_{h_1}}.
\label{contrib}
\end{eqnarray}
Si l'on considère alors la quantité indépendante de la variable $z_9$ définie par $Z_9:=\frac{Z_{h_1}z_9}{Z'_5}$, on écrit $\frac{Z'_5}{Z_{h_1}}=\frac{z_9}{Z_9}$ (avec $9 \not \in H_0 \cup H_1 \cup H_2$) si bien que la contribution de (\ref{contrib}) est
$$
\ll  B \sum_{\substack{\scriptscriptstyle z_h \leqslant \sqrt[n]{B} \\ \scriptscriptstyle h \not \in H_0 \cup H_1 \cup H_2 \\[0.5mm] \scriptscriptstyle z_9 \leqslant Z_9}} \hspace{1mm} \prod_{\substack{  h \not \in H_0 \cup H_1 \cup H_2\cup\{9\}}}z_h^{-1} \times \frac{1}{Z_9}\ll  B\log(B)^{2^n-n-2}.
$$
On peut donc remplacer la condition $Z_3 \geqslant 1$ par $Z'_5 \leqslant Z_{h_1}$.
De la même façon, en posant 
$$
Z''_5=\prod_{\substack{\scriptscriptstyle \varepsilon_4(h)=1 \\\scriptscriptstyle \varepsilon_3(h)=0}} z_h \prod_{\substack{\scriptscriptstyle\varepsilon_3(h)=1 \\ \scriptscriptstyle\varepsilon_4(h)=0 \\\scriptscriptstyle h\neq 5,7}} z_h^{-1},
$$
la condition $Z_7=Z''_5/z_5 \geqslant 1$ se réécrit $z_5 \leqslant Z''_5$. Si on a $Z_{h_1} \leqslant Z''_5$, alors puisque $z_5 \leqslant Z_{h_1}$, on a $z_5 \leqslant Z''_5$. Montrons alors que la condition $Z_7 \geqslant 1$ peut être remplacée par la condition $Z_{h_1} \leqslant Z''_5$. Pour cela, il suffit de voir que la contribution des $\mathbf{z}$ tels que $Z''_5 \leqslant Z_{h_1}$. Si l'on suppose $Z''_5 \leqslant Z_{h_1}$, la condition $z_5 \leqslant Z_{h_1}$ est alors remplacée par $z_5 \leqslant Z''_5$. Le raisonnement ci-dessus fournit alors une contribution
$$
\ll  B \sum_{z_h \leqslant \sqrt[n]{B} \atop h \not \in H_0 \cup H_1 \cup H_2} \hspace{1mm}\prod_{\substack{  h \not \in H_0 \cup H_1 \cup H_2}}z_h^{-1} \times \frac{Z''_5}{Z_{h_1}}.
$$
De même, on écrit alors $\frac{Z''_5}{Z_{h_1}}=\frac{z_9}{Z'_9}$ avec $Z'_9:=\frac{Z_{h_1}z_9}{Z''_5}$ si bien que cette contribution est
$$
\ll  B \sum_{\substack{\scriptscriptstyle z_h \leqslant \sqrt[n]{B} \\ \scriptscriptstyle h \not \in H_0 \cup H_1 \cup H_2 \\[0.5mm] \scriptscriptstyle z_9 \leqslant Z'_9}} \hspace{1mm}\prod_{\substack{  h \not \in H_0 \cup H_1 \cup H_2\cup\{9\}}}z_h^{-1} \times \frac{1}{Z'_9}\ll  B\log(B)^{2^n-n-2},
$$
ce qui est négligeable. Pour finir, on remarque que lorsque $n=3$, la condition $Z_{h_2} \geqslant 1$ ne fait pas intervenir la variable $z_{h_1}$. Il n'est donc pas nécessaire d'imposer de telles restrictions du domaine dans ce cas-là.\\
\newline
\indent
On introduit ensuite la fonction $\tilde{b}: \mathbb{N}^N \rightarrow \mathbb{R}$ définie par 
$$
\tilde{b}(\mathbf{z})={\rm{vol}}\left\{(\alpha_2,\dots,\alpha_n) \in [-1,1]^{n-1} \quad : \quad  \left|\sum_{i=2}^n d_i\alpha_i\right| \leqslant d_1\right\}.
$$
On constate en particulier que si $\mathbf{z}$ est réduit, alors $\tilde{b}(\mathbf{z})=b(\mathbf{y})$ pour $\mathbf{y}$ l'unique $n$-uplet associé à $\mathbf{z}$ à travers la bijection explicitée dans le Lemme \ref{lemmefacto} et où $b$ a été définie en (\ref{by}). On considère également la fonction multiplicative $g:\mathbb{N}^N \rightarrow \mathbb{R}$, indicatrice de l'ensemble des $N$-uplets $\mathbf{z}$ réduits.
Il s'agit ainsi, lorsque $n \geqslant 4$, de sommer $\frac{g(\mathbf{z})\tilde{b}(\mathbf{z})}{d_1}$ sur le domaine $\mathcal{V}$ suivant
\begin{eqnarray}
\left\{
\begin{array}{l}
\forall j \in \llbracket 4,n-1\rrbracket \cup \{1\}, \quad Z_{h_j}\geqslant 1, \quad Z'_5 \leqslant Z_{h_1} \leqslant Z''_5,  \quad  \mbox{et} \quad z_h \leqslant Z_h \quad \mbox{pour} \quad h \in H_0 \cup H_1,\\[2mm]
Z_{h_n} \geqslant 1 \quad \mbox{et} \quad z_{h_n} \leqslant Z_{h_n},\\
 \end{array} 
 \right.
\label{dom2}
\end{eqnarray}
le $N$-uplet $\mathbf{z}$ étant réduit. \`A une contribution négligeable près de l'ordre de 
$$
O\left( B(\log B)^{2^n-n-2}\log(\log B) \right),
$$
on peut, comme dans \cite[section 4]{Br}, se restreindre au domaine~$\mathcal{V}'$ suivant
\begin{eqnarray}
\left\{
\begin{aligned}
&\forall j \in \llbracket 4,n-1\rrbracket \cup \{1\}, \quad Z_{h_j}>\log(B)^5, \quad Z'_5(\log(B))^6 \leqslant Z_{h_1} \leqslant \frac{Z''_5}{\log(B)^6},  \quad  \mbox{et}\\
 & z_h \leqslant Z_h \quad \mbox{pour} \quad h \in H_0  \quad \mbox{et} \quad \frac{Z_{h_1}}{\log(B)}<z_5 \leqslant Z_{h_1},\\
&Z_{h_n} \geqslant \log(B)^5 \quad \mbox{et} \quad z_{h_n} \leqslant Z_{h_n},\\
\end{aligned}
\right.
\label{dom}
\end{eqnarray}
en utilisant la formule (\ref{ineg}) établie lors de la preuve du Lemme \ref{lemme3}. En particulier, on a $Z_3>\log(B)^5$ et $Z_7>\log(B)^5$. Lorsque $n=3$, il s'agit de sommer $\frac{g(\mathbf{z})\tilde{b}(\mathbf{z})}{d_1}$ sur le domaine $\mathcal{V}$ suivant
\begin{eqnarray}
\left\{
\begin{array}{l}
Z_{3}\geqslant 1, \quad Z_{h_1} \geqslant 1,  \quad  \mbox{et} \quad z_h \leqslant Z_h \quad \mbox{pour} \quad h \in H_0 \cup H_1,\\[2mm]
Z_{7} \geqslant 1 \quad \mbox{et} \quad z_{7} \leqslant Z_{7},\\
 \end{array} 
 \right.
\label{dom3}
\end{eqnarray}
le $N$-uplet $\mathbf{z}$ étant réduit. \`A une contribution négligeable près de l'ordre de $O\left(B( \log B)^{3}\log(\log B) \right)$, on peut également se restreindre au domaine~$\mathcal{V}'$ suivant
\begin{eqnarray}
\left\{
\begin{aligned}
 &Z_{3}>\log(B)^5,  \quad z_h \leqslant Z_h \quad \mbox{pour} \quad h \in H_0 \cup H_1,\\
& \sqrt[3]{B}Z_{7} \geqslant \log(B)^5 \quad \mbox{et} \quad z_{7} \leqslant \sqrt[3]{B}Z_{7}.\\
\end{aligned}
\right.
\label{dom4}
\end{eqnarray}
On utilise alors les deux lemmes suivants. Le premier traduit le fait que la fonction $g$ soit très proche, au sens de la convolution, de la fonction constante égale à 1 et le second est tiré de \cite{Br}.

\begin{lemme}
La série de Dirichlet 
$$
\forall \mathbf{s}=(s_1,\dots,s_N), \quad G(\mathbf{s})=\sum_{\mathbf{z} \in \mathbb{N}^N \atop \mathbf{z} \hspace{1mm} {\rm{ r\acute{e}duit}}} \frac{1}{\overset{N}{\underset{h=1}{\prod}} z_h^{s_h}}
$$
est convergente sur le domaine $\mathfrak{Re}(s_i)>1$ pour tout $i \in \llbracket 1,N \rrbracket$. De plus, la fonction 
\begin{eqnarray}
\forall \mathbf{s}=(s_1,\dots,s_N), \quad F(\mathbf{s})=G(\mathbf{s})\prod_{h=1}^N \zeta(s_h)^{-1} 
\label{F}
\end{eqnarray}
admet un prolongement holomorphe à la région $\mathfrak{Re}(s_h)>\frac{1}{2}$ pour tout $h \in \llbracket 1,N \rrbracket$.
\label{holo}
\end{lemme}
\noindent
\textit{Démonstration--}
On constate que la variable $z_{N}$ n'est soumise à aucune condition de coprimalité et que toutes les autres variables $z_h$ sont soumises à au moins une condition de coprimalité. On note alors $E_n$ l'ensemble des couples $(k,\ell)$ de $\llbracket 1,N-1\rrbracket^2$ tels que, si $\mathbf{z}$ est réduit, $\mbox{pgcd}(z_k,z_{\ell})=1$. Pour tout $ \mathbf{s}$ tel que $\mathfrak{Re}(s_i)>1$ pour tout $i \in \llbracket 1,N \rrbracket$, on a ainsi
\begin{eqnarray}
\begin{aligned}
G(\mathbf{s})&=\zeta(s_N)\prod_p \left( 1 +\sum_{\substack{\scriptscriptstyle (\nu_1,\dots,\nu_{N-1}) \in \mathbb{Z}_{\geqslant 0}^{N-1} \smallsetminus\{\mathbf{0}\} \\\scriptscriptstyle \nu_i \nu_j=0 \hspace{0.5mm} {\rm{avec}} \hspace{0.5mm} (i,j) \in E_n}} \frac{1}{p^{\nu_1s_1+\cdots+\nu_{N-1}s_{N-1}}}\right)\\
&=\zeta(s_N)\prod_p \left( 1 +\sum_{i=1}^{N-1} \frac{1}{p^{s_i}}+\sum_{\substack{\scriptscriptstyle(\nu_1,\dots,\nu_{N-1}) \in \mathbb{Z}_{\geqslant 0}^{N-1} \\\scriptscriptstyle \nu_i \nu_j=0 \hspace{0.5mm} {\rm{avec}} \hspace{0.5mm} (i,j) \in E_n\\\scriptscriptstyle \nu_1+\cdots+\nu_{N-1} \geqslant 2}} \frac{1}{p^{\nu_1s_1+\cdots+\nu_{N-1}s_{N-1}}}\right).\\
\end{aligned}
\label{pe}
\end{eqnarray}
Ainsi,
$$
G(\mathbf{s})\prod_{h=1}^N \zeta(s_h)^{-1} =\prod_p\Bigg(\prod_{i=1}^{N-1} \left(1-\frac{1}{p^{s_i}}\right)\Bigg) \left( 1 +\sum_{i=1}^{N-1} \frac{1}{p^{s_i}}+\sum_{\substack{\scriptscriptstyle(\nu_1,\dots,\nu_{N-1}) \in \mathbb{Z}_{\geqslant 0}^{N-1}\\\scriptscriptstyle \nu_i \nu_j=0 \hspace{0.5mm} {\rm{avec}} \hspace{0.5mm} (i,j) \in E_n\\\scriptscriptstyle \nu_1+\cdots+\nu_{N-1} \geqslant 2}} \frac{1}{p^{\nu_1s_1+\cdots+\nu_{N-1}s_{N-1}}}\right).
$$
Le produit de droite étant convergent lorsque $\mathfrak{Re}(s_i)>\frac{1}{2}$ pour tout $i \in \llbracket 1,N \rrbracket$, on obtient bien le résultat annoncé.
\hfill
$\square$

\begin{lemme}
Soient $f$ une fonction multiplicative en une variable dont la série de Dirichlet est absolument convergente pour $\mathfrak{Re}(s)\geqslant \frac{2}{3}$ et $v$ une fonction bornée et dérivable sur $[0,1]$. On a alors pour tout $Z \geqslant 1$
$$
\sum_{z \leqslant Z}(\mathds{1} \ast f)(z)v\left(\frac{z}{Z}\right)=Z\sum_{k\geqslant 1} \frac{f(k)}{k}\int_0^1 v(u)\mbox{d}u+O\left(  Z^{2/3}\sum_{k\geqslant 1} \frac{|f(k)|}{k^{2/3}}\int_0^1\left|v'(u)\right|u^{2/3}\mbox{d}u\right).
$$
\label{lemmeaux}
\end{lemme}
\noindent
\textit{Démonstration--}
La démonstration de ce lemme s'obtient aisément à l'aide d'une sommation d'Abel (voir \cite[section 4]{Br}).
\hfill
$\square$\\
\newline
On notera dans la suite $f:\mathbb{N}^N \rightarrow \mathbb{R}$ la fonction arithmétique associée à la série de Dirichlet~$F$. Autrement dit, pour tout $\mathbf{s}$ tel que $\mathfrak{Re}(s_h)>1$ pour tout $h \in \llbracket 1,N \rrbracket$, on a
$$
F(\mathbf{s})=\sum_{\mathbf{z} \in \mathbb{N}^N} \frac{f(\mathbf{z})}{\overset{N}{\underset{h=1}{\prod}} z_h^{s_h}}.
$$
et $g=\mathds{1} \ast f$. Une application des Lemmes \ref{holo} et \ref{lemmeaux} en sommant d'abord sur les $z_h$ pour $h \in H_0$ et le fait que $d_1$ ne fasse intervenir aucune variable $z_{h_j}$ pour $j \in \llbracket 2,n-1 \rrbracket$ fournit alors
$$
\sum_{z_{h_j} \leqslant Z_{h_j}\atop h_j \in H_0} g(\mathbf{z})\tilde{b}(\mathbf{z})=Z^{(0)}\beta\left(\frac{z_{h_1}}{Z_{h_1}}\right)\left(\sum_{k_h \mid z_h \atop h \not \in H_0}\sum_{k_h \geqslant 1 \atop h \in H_0} \frac{f(\mathbf{k})}{\displaystyle\underset{\scriptscriptstyle h \in H_0}{\prod}k_h}+O\left(\frac{1}{\log(B)}\sum_{k_h \mid z_h \atop h \not \in H_0}\sum_{k_h \geqslant 1 \atop h \in H_0} \frac{|f(\mathbf{k})|}{\displaystyle\underset{\scriptscriptstyle h \in H_0}{\prod}k_h^{2/3}}\right)\right)
$$
avec $Z^{(0)}=\overset{n-1}{\underset{j=2}{\prod}} Z_{h_j}$ et
$$
\beta(u_1)=\int_{[0,1]^{n-2}}{\rm{vol}}\left\{(\alpha_1,\dots,\alpha_n) \in [-1,1]^n \quad : \quad  \left|\sum_{i=2}^n \left(\prod_{\ell=1}^{i-1}u_{\ell}\right)\alpha_i\right| \leqslant 1\right\}\mbox{d}u_2\cdots\mbox{d}u_{n-1}.
$$
En effet, on a les égalités
$$
\forall 2 \leqslant \ell \leqslant n-1, \quad \frac{d_{\ell+1}}{d_{\ell}}=\frac{z_{h_{\ell}}}{Z_{h_{\ell}}} \quad \mbox{et} \quad  \frac{d_{2}}{d_{1}}=\frac{z_{h_1}}{Z_{1}}
$$
et, pour tout $\mathbf{z}$,
$$
\begin{aligned}
\tilde{b}(\mathbf{z})&={\rm{vol}}\left\{(\alpha_1,\dots,\alpha_n) \in [-1,1]^n \quad : \quad  \left|\sum_{i=2}^n \frac{d_i}{d_1}\alpha_i\right| \leqslant 1\right\}\\
&={\rm{vol}}\left\{(\alpha_1,\dots,\alpha_n) \in [-1,1]^n \quad :\quad  \left|\sum_{i=2}^n \left(\prod_{\ell=1}^{i-1}\frac{d_{\ell+1}}{d_{\ell}}\right)\alpha_i\right| \leqslant 1\right\}.\\
\end{aligned}
$$
est une fonction différentiable à dérivées partielles bornées sur $[0,1]$ si bien que chacun des termes
$$
\int_0^1\left|v'(u)\right|u^{2/3}\mbox{d}u \ll 1.
$$
De plus, grâce aux conditions (\ref{dom}), dans le terme d'erreur du Lemme \ref{lemmeaux}, on a bien
$$
\left(Z^{(0)}\right)^{2/3}=Z^{(0)}\left(Z^{(0)}\right)^{-1/3} \ll \frac{Z^{(0)}}{\log(B)^{5(n-2)/3}} \ll \frac{Z^{(0)}}{\log(B)}. 
$$
On pourra noter que $Z^{(0)}$ ne dépend pas de $Z_{h_1}$ d'après (\ref{z5}). On effectue alors la sommation par rapport à $z_{h_1}$ pour obtenir
$$
\!\!\!\sum_{z_{h} \leqslant Z_{h}\atop h \in H_0 \cup H_1} \!\!\!g(\mathbf{z})\tilde{b}(\mathbf{z})=Z^{(1)}\tilde{\beta}\!\left(\sum_{k_h \mid z_h \atop h \not \in H_0 \cup H_1}\sum_{k_h \geqslant 1 \atop h \in H_0 \cup H_1} \!\frac{f(\mathbf{k})}{\displaystyle\underset{\scriptscriptstyle h \in H_0 \cup H_1}{\prod}k_h}+O\left(\!\!\frac{1}{\log(B)}\!\!\!\sum_{k_h \mid z_h \atop h \not \in H_0 \cup H_1}\sum_{k_h \geqslant 1 \atop h \in H_0 \cup H_1} \frac{|f(\mathbf{k})|}{\displaystyle\underset{\scriptscriptstyle h \in H_0 \cup H_1}{\prod}k_h^{2/3}}\right)\right),
$$
où $Z^{(1)}=Z^{(0)}Z_{h_1}$ et 
\begin{eqnarray}
\tilde{\beta}=\int_0^1\beta(u_1)\mbox{d}u_1.
\label{betatilde}
\end{eqnarray}
Les quantités $Z^{(1)}$ et $d_1$ étant indépendantes de $z_{h_n}$, en sommant sur $z_{h_n}$, il vient
$$
\begin{aligned}
\sum_{z_{h_j} \leqslant Z_{h_j}\atop h_j \in H_0 \cup H_1 \cup H_2} g(\mathbf{z})\tilde{b}(\mathbf{z})&=Z^{(2)}\tilde{\beta}\left(\sum_{k_h \mid z_h \atop h \not \in H_0 \cup H_1 \cup H_2}\sum_{k_h \geqslant 1 \atop h \in H_0 \cup H_1 \cup H_2} \frac{f(\mathbf{k})}{\displaystyle\underset{\scriptscriptstyle h \in H_0 \cup H_1 \cup H_2}{\prod}k_h}\right.\\
& \quad \left.+O\left(\frac{1}{\log(B)}\sum_{k_h \mid z_h \atop h \not \in H_0 \cup H_1\cup H_2}\sum_{k_h \geqslant 1 \atop h \in H_0 \cup H_1\cup H_2} \frac{|f(\mathbf{k})|}{\displaystyle\underset{\scriptscriptstyle h \in H_0 \cup H_1 \cup H_2}{\prod}k_h^{2/3}}\right)\right),
\end{aligned}
$$
avec $Z^{(2)}=Z^{(1)}Z_{h_n}$. On a ainsi
$$
\begin{aligned}
\!\!\!\!\sum_{\mathbf{y} \in \mathbb{N}^n \atop 1 \leqslant y_1\leqslant \cdots\leqslant y_n\leqslant \sqrt[n]{B}/k}\!\!\!\!\frac{b(\mathbf{y})}{d_1}&=\sqrt[n]{B}\tilde{\beta}\!\!\!\sum_{z_h  \in \mathcal{V}', \hspace{1mm} z_h \leqslant \sqrt[n]{B}/k \atop h \not \in H_0 \cup H_1 \cup H_2} \prod_{h \not \in H_0 \cup H_1 \cup H_2}z_h^{-1}\left(\sum_{k_h \mid z_h \atop h \not \in H_0 \cup H_1 \cup H_2}\sum_{k_h \geqslant 1 \atop h \in H_0 \cup H_1 \cup H_2} \frac{f(\mathbf{k})}{\displaystyle\underset{\scriptscriptstyle h \in H_0 \cup H_1 \cup H_2}{\prod}k_h}\right.\\
& \quad \left.+O\left(\frac{1}{\log(B)}\sum_{k_h \mid z_h \atop h \not \in H_0 \cup H_1\cup H_2}\sum_{k_h \geqslant 1 \atop h \in H_0 \cup H_1\cup H_2} \frac{|f(\mathbf{k})|}{\displaystyle\underset{\scriptscriptstyle h \in H_0 \cup H_1 \cup H_2}{\prod}k_h^{2/3}}\right)\right).
\end{aligned}
$$
Dans un premier temps, on remarque que l'on peut sommer sur $\mathcal{V}$ défini en (\ref{dom2}) lorsque $n \geqslant 4$ et en (\ref{dom3}) lorsque $n=3$ quitte à rajouter une contribution négligeable. En effet, considérons $(z_h)_{h \not\in H_0 \cup H_1 \cup H_2} \in \mathcal{V} \smallsetminus\mathcal{V}'$. S'il existe $4 \leqslant \ell \leqslant n-1$ ou $\ell=1$ tel que $Z_{h_{\ell}} \leqslant \log(B)^5$, alors
$$
 \prod_{\substack{\scriptscriptstyle \varepsilon_j(h)=1 \\\scriptscriptstyle \varepsilon_{j+1}(h)=0 \\\scriptscriptstyle h \neq h_j}} z_h^{-1} \leqslant  \log(B)^5 \prod_{\substack{\scriptscriptstyle\varepsilon_j(h)=0 \\\scriptscriptstyle \varepsilon_{j+1}(h)=1}} z_h^{-1}
$$
si bien que
$$
\begin{aligned}
\sum_{z_h  \in \mathcal{V} \smallsetminus \mathcal{V}', \hspace{1mm} z_h \leqslant \sqrt[n]{B} \atop h \not \in H_0 \cup H_1 \cup H_2} \prod_{h \not \in H_0 \cup H_1 \cup H_2}z_h^{-1} &\ll \log(B)^5\sum_{z_h \leqslant \sqrt[n]{B} \atop h \not \in H_0 \cup H_1 \cup H_2} \prod_{\substack{\scriptscriptstyle h \not \in H_0 \cup H_1 \cup H_2\\\scriptscriptstyle \varepsilon_j(h)=\varepsilon_{j+1}(h) } }z_h^{-1}\prod_{\substack{\scriptscriptstyle h \not \in H_0 \cup H_1 \cup H_2\\\scriptscriptstyle \varepsilon_j(h)=0 \\\scriptscriptstyle \varepsilon_{j+1}(h)=1 } }z_h^{-2}\\
&\ll \log(B)^5 B^{2^{n-2}-1}\frac{1}{B^{2^{n-2}}}\log(B)^{2^n-n-1-2^{n-1}-1}\\
& \ll \frac{\log(B)^{2^{n-1}-n+2}}{B} \ll 1. \\
\end{aligned}
$$
De même, si $Z_5 \leqslant Z'_5(\log(B))^6$, on a
$$
\frac{Z_5}{Z'_5}=\prod_{\substack{\scriptscriptstyle\varepsilon_1(h)=0 \\\scriptscriptstyle \varepsilon_{3}(h)=1 \\\scriptscriptstyle h \neq 5 } }z_h\prod_{\substack{\scriptscriptstyle\varepsilon_1(h)=1 \\\scriptscriptstyle \varepsilon_{3}(h)=0 \\\scriptscriptstyle h \neq 3,5 } }z_h^{-1}
$$
et lorsque $Z_5>\frac{Z''_5}{\log(B)^6}$, il vient
$$
\frac{Z''_5}{Z''_5}=\prod_{\substack{\scriptscriptstyle\varepsilon_1(h)=0 \\ \scriptscriptstyle\varepsilon_{2}(h)=1 \\\scriptscriptstyle \varepsilon_4(h)\neq 1 \hspace{1mm} {\rm{ou}} \hspace{1mm}\varepsilon_3(h)\neq 0 } }z_h\prod_{\substack{\scriptscriptstyle\varepsilon_1(h)=1 \\\scriptscriptstyle \varepsilon_{2}(h)=0 \\\scriptscriptstyle \varepsilon_4(h)\neq 0 \hspace{1mm} {\rm{ou}} \hspace{1mm}\varepsilon_3(h)\neq 1} }z_h\prod_{\substack{\scriptscriptstyle\varepsilon_4(h)=1 \\\scriptscriptstyle \varepsilon_{3}(h)=0 \\\scriptscriptstyle \varepsilon_1(h)\neq 0 \hspace{1mm} {\rm{ou}} \hspace{1mm}\varepsilon_2(h)\neq 1} }z_h^{-1}\prod_{\substack{\scriptscriptstyle\varepsilon_4(h)=0 \\\scriptscriptstyle \varepsilon_{3}(h)=1 \\\scriptscriptstyle \varepsilon_1(h)\neq 1 \hspace{1mm} {\rm{ou}} \hspace{1mm}\varepsilon_2(h)\neq 0 \\\scriptscriptstyle h \neq 5} }z_h^{-1}
$$
de sorte que le même raisonnement
permet de conclure à une contribution $\ll 1$. On a donc 
$$
\begin{aligned}
\!\!\!\!\sum_{\mathbf{y} \in \mathbb{N}^n \atop 1 \leqslant y_1\leqslant \cdots\leqslant y_n\leqslant \sqrt[n]{B}}\!\!\!\!\frac{b(\mathbf{y})}{d_1}&=\sqrt[n]{B}\tilde{\beta}\!\!\sum_{z_h  \in \mathcal{V}, \hspace{1mm} z_h \leqslant \sqrt[n]{B} \atop h \not \in H_0 \cup H_1 \cup H_2} \prod_{h \not \in H_0 \cup H_1 \cup H_2}z_h^{-1}\left(\sum_{k_h \mid z_h \atop h \not \in H_0 \cup H_1 \cup H_2}\sum_{k_h \geqslant 1 \atop h \in H_0 \cup H_1 \cup H_2} \frac{f(\mathbf{k})}{\displaystyle\underset{\scriptscriptstyle h \in H_0 \cup H_1 \cup H_2}{\prod}k_h}\right.\\
& \quad \left.+O\left(\frac{1}{\log(B)}\sum_{k_h \mid z_h \atop h \not \in H_0 \cup H_1\cup H_2}\sum_{k_h \geqslant 1 \atop h \in H_0 \cup H_1\cup H_2} \frac{|f(\mathbf{k})|}{\displaystyle\underset{\scriptscriptstyle h \in H_0 \cup H_1 \cup H_2}{\prod}k_h^{2/3}}\right)\right).
\end{aligned}
$$
Le domaine $\mathcal{V}$ se réécrit sous la forme
$$
\left\{
\begin{aligned}
&\prod_{\substack{\scriptscriptstyle\varepsilon_j(h)=0 \\\scriptscriptstyle \varepsilon_{j+1}(h)=1}} z_h^{-1} \prod_{\substack{\scriptscriptstyle\varepsilon_j(h)=1 \\\scriptscriptstyle \varepsilon_{j+1}(h)=0 \\\scriptscriptstyle h \neq h_j}} z_h \leqslant 1 \quad \mbox{pour} \quad 4 \leqslant j \leqslant n-1,\\
&\prod_{\substack{\scriptscriptstyle\varepsilon_n(h)=1 \\\scriptscriptstyle h \neq h_n}} z_h \leqslant \sqrt[n]{B}, \quad \prod_{\substack{\scriptscriptstyle\varepsilon_1(h)=0 \\\scriptscriptstyle \varepsilon_{3}(h)=1 } }z_h^{-1}\prod_{\substack{\scriptscriptstyle\varepsilon_1(h)=1 \\\scriptscriptstyle \varepsilon_{3}(h)=0 \\\scriptscriptstyle h \neq 3 } }z_h\leqslant 1, \quad \prod_{\substack{\scriptscriptstyle\varepsilon_1(h)=0 \\ \scriptscriptstyle\varepsilon_{2}(h)=1 } }z_h^{-1}\prod_{\substack{\scriptscriptstyle\varepsilon_1(h)=1 \\\scriptscriptstyle \varepsilon_{2}(h)=0 \\\scriptscriptstyle h \neq 5 } }z_h\leqslant 1\\
&\prod_{\substack{\scriptscriptstyle\varepsilon_1(h)=0 \\\scriptscriptstyle \varepsilon_{2}(h)=1  } }z_h\prod_{\substack{\scriptscriptstyle\varepsilon_1(h)=1 \\\scriptscriptstyle \varepsilon_{2}(h)=0\\\scriptscriptstyle h \neq 5} }z_h^{-1}\prod_{\substack{\scriptscriptstyle\varepsilon_4(h)=1 \\\scriptscriptstyle \varepsilon_{3}(h)=0 } }z_h^{-1}\prod_{\substack{\scriptscriptstyle\varepsilon_4(h)=0 \\\scriptscriptstyle \varepsilon_{3}(h)=1 \\\scriptscriptstyle h \neq 5,7} }z_h \leqslant 1,\\
&z_h \leqslant \sqrt[n]{B} \quad \mbox{pour} \quad h \not \in H_0 \cup H_1 \cup H_2.
\end{aligned}
\right.
$$
lorsque $n \geqslant 4$ et
$$
\left\{
\begin{aligned}
&\frac{z_2}{z_4z_5} \leqslant 1,\\
&z_4z_5z_6 \leqslant \sqrt[3]{B}, \quad \frac{z_5}{z_2z_6} \leqslant 1,\\
&z_h \leqslant \sqrt[3]{B} \quad \mbox{pour} \quad h \not \in H_0 \cup H_1 \cup H_2.
\end{aligned}
\right.
$$
lorsque $n=3$. En écrivant $z_h=B^{\frac{t_h}{n}}$ avec $t_h \geqslant 0$ pour tout $h \not \in H_0 \cup H_1 \cup H_2$, on obtient que l'ensemble des $(t_h)_{h \not \in H_0 \cup H_1 \cup H_2}$ vérifie
$$
\left\{
\begin{aligned}
&\sum_{\substack{\scriptscriptstyle\varepsilon_j(h)=1 \\\scriptscriptstyle \varepsilon_{j+1}(h)=0\\\scriptscriptstyle h \neq h_j}} t_h \leqslant \sum_{\substack{\scriptscriptstyle\varepsilon_j(h)=0 \\\scriptscriptstyle \varepsilon_{j+1}(h)=1}} t_h \quad \mbox{pour} \quad 4 \leqslant j \leqslant n-1,\\
&\sum_{\substack{\scriptscriptstyle\varepsilon_n(h)=1 \\\scriptscriptstyle h \neq h_n}} t_h\leqslant 1, \quad \sum_{\substack{\scriptscriptstyle\varepsilon_1(h)=1 \\\scriptscriptstyle \varepsilon_{3}(h)=0 \\\scriptscriptstyle h \neq 3 } }t_h \leqslant \sum_{\substack{\scriptscriptstyle\varepsilon_1(h)=0 \\\scriptscriptstyle \varepsilon_{3}(h)=1 } }t_h, \quad \sum_{\substack{\scriptscriptstyle\varepsilon_1(h)=1 \\\scriptscriptstyle \varepsilon_{2}(h)=0 \\\scriptscriptstyle h \neq 5 } }t_h \leqslant \sum_{\substack{\scriptscriptstyle\varepsilon_1(h)=0 \\\scriptscriptstyle \varepsilon_{2}(h)=1 } }t_h\\
&\sum_{\substack{\scriptscriptstyle\varepsilon_1(h)=0 \\\scriptscriptstyle \varepsilon_{2}(h)=1  } }t_h+\sum_{\substack{\scriptscriptstyle\varepsilon_4(h)=0 \\\scriptscriptstyle \varepsilon_{3}(h)=1 \\\scriptscriptstyle h \neq 5,7} }t_h \leqslant \sum_{\substack{\scriptscriptstyle\varepsilon_1(h)=1 \\\scriptscriptstyle \varepsilon_{2}(h)=0\\\scriptscriptstyle h \neq 5} }t_h+\displaystyle\sum_{\substack{\scriptscriptstyle\varepsilon_4(h)=1 \\\scriptscriptstyle \varepsilon_{3}(h)=0 } }t_h,\\
&t_h \leqslant 1 \quad \mbox{pour} \quad h \not \in H_0 \cup H_1 \cup H_2,
\end{aligned}
\right.
$$ 
lorsque $n \geqslant 4$ et
$$
\left\{
\begin{aligned}
&t_2\leqslant t_4+t_5\\
&t_4+t_5+t_6 \leqslant 1, \quad t_5 \leqslant t_2+t_6,\\
&t_2,t_4,t_5,t_6 \leqslant 1\\
\end{aligned}
\right.
$$
si $n=3$ si bien que dans tous les cas, il est bien inclus dans $[0,1]^{2^n-n-1}$. Une application directe de~\cite[lemma 8]{LB} fournit alors l'estimation
$$
N\left(\frac{B}{k^n}\right)=\frac{2^{n-1} n!}{n^{2^n-n-1}k^n}\tilde{\beta}VF(\mathbf{1})B\log(B)^{2^n-n-1}+O\left( \frac{B}{k^{n}} (\log B)^{2^n-n-2}\log(\log B) \right),
$$
avec 
\begin{eqnarray}
\!V={\rm{vol}}
\left\{ (t_h) \in [0,1]^{2^n-n-1} \hspace{1mm} :
\begin{array}{c}
\displaystyle\underset{\substack{\scriptscriptstyle\varepsilon_j(h)=1 \\\scriptscriptstyle \varepsilon_{j+1}(h)=0\\\scriptscriptstyle h \neq h_j}}{\sum} t_h \leqslant \underset{\substack{\scriptscriptstyle\varepsilon_j(h)=0 \\\scriptscriptstyle \varepsilon_{j+1}(h)=1}}{\sum} t_h \quad \mbox{pour} \quad 4 \leqslant j \leqslant n-1,\\[3mm]
\displaystyle\underset{\substack{\scriptscriptstyle\varepsilon_n(h)=1 \\\scriptscriptstyle h \neq h_n}}{\sum} t_h\leqslant 1, \quad \underset{\substack{\scriptscriptstyle\varepsilon_1(h)=1 \\\scriptscriptstyle \varepsilon_{3}(h)=0 \\\scriptscriptstyle h \neq 3 } }{\sum}t_h \leqslant \underset{\substack{\scriptscriptstyle\varepsilon_1(h)=0 \\\scriptscriptstyle \varepsilon_{3}(h)=1 } }{\sum}t_h, \hspace{1mm} \\[5mm]
\displaystyle\underset{\substack{\scriptscriptstyle\varepsilon_1(h)=1 \\ \scriptscriptstyle\varepsilon_{2}(h)=0 \\\scriptscriptstyle h \neq 5 } }{\sum}t_h \leqslant \underset{\substack{\scriptscriptstyle\varepsilon_1(h)=0 \\ \scriptscriptstyle\varepsilon_{2}(h)=1 } }{\sum}t_h, \\[5mm]
\displaystyle\underset{\substack{\scriptscriptstyle\varepsilon_1(h)=0 \\ \scriptscriptstyle\varepsilon_{2}(h)=1  } }{\sum}t_h+\underset{\substack{\scriptscriptstyle\varepsilon_4(h)=0 \\ \scriptscriptstyle\varepsilon_{3}(h)=1 \\\scriptscriptstyle h \neq 5,7} }{\sum}t_h \leqslant \underset{\substack{\scriptscriptstyle\varepsilon_1(h)=1 \\\scriptscriptstyle \varepsilon_{2}(h)=0\\\scriptscriptstyle h \neq 5} }{\sum}t_h+\underset{\substack{\scriptscriptstyle\varepsilon_4(h)=1 \\\scriptscriptstyle \varepsilon_{3}(h)=0 } }{\sum}t_h\\[3mm]
\end{array}
\right\}
\label{V}
\end{eqnarray}
lorsque $n \geqslant 4$ et
$$
V={\rm{vol}}
\left\{ (t_2,t_4,t_5,t_6) \in [0,1]^{4} \quad :
\begin{array}{l}
t_2\leqslant t_4+t_5\\[2mm]
t_4+t_5+t_6 \leqslant 1, \quad t_5 \leqslant t_2+t_6
\end{array}
\right\}
$$
si $n=3$. Finalement, il vient
\begin{eqnarray}
N(B;U_n)=\frac{2^{n-1} n!\tilde{\beta}V}{n^{2^n-n-1}\zeta(n)}F(\mathbf{1}) B\log(B)^{2^n-n-1}+O\left( B(\log B)^{2^n-n-2}\log(\log B) \right).
\label{cn}
\end{eqnarray}
\textbf{Remarque.--} 
Géométriquement, on a transformé le problème de comptage sur la variété $W_n$ en un problème de comptage sur la sous-variété (\ref{torsor}) 
\begin{eqnarray}
\sum_{j=1}^n d_jx_j=0
\quad
\mbox{avec}
\quad 
\forall i \in \llbracket 1,n \rrbracket, \quad d_i=\prod_{1 \leqslant h \leqslant N} z_h^{1-\varepsilon_i(h)}.
\label{torsor}
\end{eqnarray}
On montrera plus tard que cette variété est liée au torseur versel d'une résolution crépante de $W_n$. Cela est cohérent avec les résultats et le torseur versel obtenus dans \cite{Blomer2014}.


\section{Vérification de la conjecture de Peyre pour $W_n$}

\subsection{Résolution crépante des singularités de $W_n$ et forme conjecturale de la constante de Peyre}
On s'appuie ici sur le travail de Per Salberger présenté en Annexe de cet article. L'objet de cette partie est de construire à partir de cet Annexe, une résolution crépante des singularités de $W_n$ puis de détailler tous les éléments de la géométrie de cette résolution crépante nécessaires à la vérification du fait que $c_n=c_{{\rm Peyre}}$ pour tout $n \geqslant 3$ où $c_n$ est la constante obtenue dans le Théorème \ref{theor1} et $c_{{\rm Peyre}}$ est définie en (\ref{cp}) \textit{infra} ou par \cite[formule 5.1]{P95}.

\subsubsection{Une résolution crépante de $W_n$} 
On reformule dans cette section le résultat principal de l'Annexe de Salberger en utilisant les notations de cet article sans donner aucune preuve. On renvoie le lecteur intéressé par ces dernières à cet Annexe en fin d'article.\\
\par
On rappelle ici que $N=2^n-1$ et on introduit pour tout $h \in \llbracket 1,N \rrbracket$, $\mathbb{P}^{(h)}\times \mathbb{P}^{(h)}$ comme étant l'espace biprojectif $\mathbb{P}^{s(h)-1}\times \mathbb{P}^{s(h)-1}$ de coordonnées bihomogènes $$\left(\mathbf{Y}^{(h)};\mathbf{Z}^{(h)}\right)=\left( Y^{(h)}_{i_1},\dots,Y^{(h)}_{i_{s(h)}};Z^{(h)}_{i_1},\dots,Z^{(h)}_{i_{s(h)}} \right)$$ pour $\varepsilon_{i_1}(h)=\cdots=\varepsilon_{i_{s(h)}}(h)=1$ et $B^{(h)} \subseteq \mathbb{P}^{(h)}\times \mathbb{P}^{(h)}$ la sous-variété fermée définie par les équations
\begin{eqnarray}
Y^{(h)}_{i_1}Z^{(h)}_{i_1}=\cdots=Y^{(h)}_{i_{s(h)}}Z^{(h)}_{i_{s(h)}} \quad \mbox{pour} \quad \varepsilon_{i_1}(h)=\cdots=\varepsilon_{i_{s(h)}}(h)=1.\label{4}
\end{eqnarray}
On pose également $B_{0,n}$ comme étant la sous-variété fermée de $\prod_{1 \leqslant h \leqslant N} B^{(h)}$ définie par les équations 
\begin{numcases}{}
Y^{(h)}_i Y^{(\ell)}_{j}=Y^{(h)}_{j} Y^{(\ell)}_{i} \quad \mbox{ pour } \quad  \ell \preceq h \in \llbracket 1,N \rrbracket \quad \mbox{ et } \quad \varepsilon_i(\ell)=\varepsilon_j(\ell)=1, \label{5}\\
Z^{(h)}_i Z^{(\ell)}_{j}=Z^{(h)}_{j} Z^{(\ell)}_{i} \quad \mbox{ \hspace{0.75mm}pour } \quad \ell \preceq h \in \llbracket 1,N \rrbracket \quad \mbox{ et } \quad \varepsilon_i(\ell)=\varepsilon_j(\ell)=1.  \label{6}
\end{numcases} 
On a un morphisme évident donné par la projection sur le dernier facteur $p_0: B_{0,n} \rightarrow B^{(N)}$ défini par $\prod_{1 \leqslant h \leqslant N} \left( \mathbf{Y}^{(h)};\mathbf{Z}^{(h)} \right) \mapsto\left( \mathbf{Y}^{(N)};\mathbf{Z}^{(N)} \right)$.\\
\par
On introduit également $C_n \subseteq \prod_{1 \leqslant h \leqslant N} \mathbb{P}^{(h)}$, la variété torique de Coxeter de $\mathfrak{S}_n$ de coordonnées multihomogènes $\prod_{1 \leqslant h \leqslant N} \left( \mathbf{Y}^{(h)} \right)$ définie par les équations $(\ref{4})$. Enfin, on définit la sous-variété fermée $X_{0,n} \subseteq \mathbb{P}^{2n-1} \times \prod_{1 \leqslant h \leqslant N} \mathbb{P}^{(h)}\times \mathbb{P}^{(h)}$ de coordonnées multihomogènes $$\left(x_1,\dots,x_n,y_1,\dots,y_n;\prod_{1 \leqslant h \leqslant N} \left( \mathbf{Y}^{(h)};\mathbf{Z}^{(h)}\right)\right)$$ définie par les équations suivantes
\begin{numcases}{}
Y^{(h)}_{i_1}Z^{(h)}_{i_1}=\cdots=Y^{(h)}_{i_{s(h)}}Z^{(h)}_{i_{s(h)}} \quad \mbox{pour} \quad h \in \llbracket 1,N \rrbracket \quad \mbox{et} \quad \varepsilon_{i_1}(h)=\cdots=\varepsilon_{i_{s(h)}}(h)=1, \label{4}\\
Y^{(h)}_i Y^{(\ell)}_{j}=Y^{(h)}_{j} Y^{(\ell)}_{i} \quad \mbox{ pour } \quad  \ell \preceq h \in \llbracket 1,N \rrbracket \quad \mbox{ et } \quad \varepsilon_i(\ell)=\varepsilon_j(\ell)=1, \label{5}\\
Z^{(h)}_i Z^{(\ell)}_{j}=Z^{(h)}_{j} Z^{(\ell)}_{i} \quad \mbox{\hspace{0.75mm} pour } \quad \ell \preceq h \in \llbracket 1,N \rrbracket \quad \mbox{ et } \quad \varepsilon_i(\ell)=\varepsilon_j(\ell)=1, \label{6}\\
x_1 Z^{(N)}_1+\cdots+x_n Z^{(N)}_n=0, \label{7}\\
 y_iY^{(N)}_j-y_jY^{(N)}_i=0 \quad \mbox{ pour } \quad  1 \leqslant i<j\leqslant n. \label{8}
\end{numcases}
Le résultat principal de l'Annexe, dû à Per Salberger, est alors le suivant et permet d'obtenir une résolution crépante de $W_n$.
\begin{theor}[Salberger, {[Annexe]}]
Soit $n \geqslant 1$. La restriction de la projection sur le premier facteur $$\mbox{pr}_1: \mathbb{P}^{2n-1} \times \prod_{1 \leqslant h \leqslant N} \mathbb{P}^{(h)}\times \mathbb{P}^{(h)} \rightarrow \mathbb{P}^{2n-1}$$ à la variété $X_{0,n}$ définie par les équations (\ref{4}), (\ref{5}), (\ref{6}), (\ref{7}) et (\ref{8}) fournit alors une résolution crépante $f_{0,n}:X_{0,n} \rightarrow W_n$ des singularités de $W_n$. De plus, $X_{0,n}$ est un $\mathbb{P}^{n-1}$-fibré sur une variété $B_{0,n}$ isomorphe à la variété torique de Coxeter $C_n$ de $\mathfrak{S}_n$.
\label{crep}
\end{theor}

On rappelle également que dans le cas $n=3$ une résolution crépante est construite dans \cite{Blomer2014}. Cette dernière est isomorphe à celle fournie par le Théorème \ref{crep} comme le remarque Salberger à la suite du théorème 1 de son Annexe. Cette résolution crépante est un $\mathbb{P}^2$-fibré sur la variété torique $B^{(3)}$. Plus généralement, Blomer, Brüdern et Salberger considèrent dans \cite{Blomer2014} la variété triprojective $X_n \subseteq \mathbb{P}^{2n-1} \times \mathbb{P}^{n-1} \times \mathbb{P}^{n-1}$ de coordonnées homogènes $(x_1,\dots,x_n,y_1,\dots,y_n;Y_1,\dots,Y_n;Z_1,\dots,Z_n)$ définie par les équations suivantes
%
\begin{numcases}{}
 x_1Z_1+\cdots+x_nZ_n=0, \label{1}\\
 y_iY_j-y_jY_i=0 \quad \mbox{ pour } \quad  1 \leqslant i<j\leqslant n, \label{2}\\
 Y_1Z_1=\cdots=Y_nZ_n. \label{3}
\end{numcases}
Il est alors établi dans \cite{Blomer2014} que la projection $\mbox{pr}_1: \mathbb{P}^{2n-1} \times \mathbb{P}^{n-1} \times \mathbb{P}^{n-1} \rightarrow \mathbb{P}^{2n-1}$ donne lieu par restriction à un morphisme propre, $G_n-$équivariant et crépant $f_n:X_n \rightarrow W_n$. Le morphisme $f_3$ est alors la résolution crépante obtenue pour $n=3$ dans \cite{Blomer2014}. En revanche, $X_n$ n'est pas lisse dès que $n \geqslant 4$.

\subsubsection{Les hypersurfaces $W_n$ sont "presque de Fano"}
La proposition suivante permet de justifier que la formule empirique \cite[formule 5.1]{P03} de Peyre s'applique bien dans le cas des variétés $W_n$ pour $n \geqslant 3$.
\begin{prop}
La variété $X_{0,n} \subseteq \mathbb{P}^{2n-1} \times \prod_{1 \leqslant h \leqslant N} \mathbb{P}^{(h)}\times \mathbb{P}^{(h)}$ définie par les équations (\ref{4}), (\ref{5}), (\ref{6}), (\ref{7}) et (\ref{8}) est "presque de Fano" au sens de \cite[Définition 3.1]{P03}.
\label{propfano}
\end{prop}
\noindent
\textit{Démonstration.--}
Le Lemme \ref{lemmegeo} entraîne aisément, de la même façon que dans \cite[lemma 6]{Blomer2014}, que le groupe de Picard géométrique de $X_{0,n}$ est sans torsion et que la classe anticanonique est dans l'intérieur de ${\rm C}_{{\rm eff}}(X_{0,n})$. Il ne reste donc qu'à justifier le fait que
$$
H^1(X_{0,n},O_{X_{0,n}})=H^2(X_{0,n},O_{X_{0,n}})=\{0\}.
$$
On utilise alors le fait que 
$$
H^1(B_{0,n},O_{B_{0,n}})=H^2(B_{0,n},O_{B_{0,n}})=\{0\}.
$$
d'après \cite[section 3.3]{Bour} et le fait que $X_{0,n}$ soit un $\mathbb{P}^{n-1}$-fibré sur $B_{0,n}$ permet de conclure.
\hfill
$\square$

\subsubsection{Forme conjecturale de la constante de Peyre}
Les conjectures originales de Manin \cite{FMT} et de Peyre \cite{P95} sur les variétés de Fano non singulières ne s'appliquent pas directement au problème de comptage associé à $W_n$ puisque cette dernière est une hypersurface singulière pour $n \geqslant 3$. Néanmoins, puisque $f_{0,n}:X_{0,n} \rightarrow W_n$ est une résolution crépante de $W_n$ dont la restriction à $X^{\circ}_0\rightarrow U_n$ est un isomorphisme où $X^{\circ}_{0,n}$ est l'ouvert de $X_{0,n}$ défini par les conditions $y_1\cdots y_n \neq 0$, alors on a
$$
N(B;U_n)=\#\left\{ x \in X^{\circ}_{0,n}(\mathbb{Q}) \hspace{1mm} : \hspace{1mm} \left(H \circ f_{0,n}\right)(x) \leqslant B \right\}
$$
où la hauteur $H \circ f_{0,n}$ est une hauteur anticanonique puisque $f_{0,n}$ est crépante. Comme d'après la Proposition \ref{propfano}, la variété $X_{0,n}$ est une variété "presque de Fano" au sens de \cite[Définition 3.1]{P03}, alors la conjecture de Manin prend la forme suivante où la constante de Peyre est donnée par la formule empirique \cite[formule 5.1]{P03} 
$$
N(B;U_n)=c_{{\rm Peyre}} B\left( \log(B)\right)^{{\rm rk}({\rm Pic}(X_{0,n}))-1}(1+o(1))
$$
avec
\begin{eqnarray}
c_{{\rm Peyre}}=\alpha(X_{0,n})\beta(X_{0,n})\omega_H\left( X_{0,n}(\mathbb{A}_{\mathbb{Q}})^{{\rm Br}(X_{0,n})} \right)
\label{cp}
\end{eqnarray}
et
\begin{eqnarray}
\beta(X_{0,n})=\#H^1\left(\mbox{Gal}(\overline{\mathbb{Q}},\mathbb{Q}),\mbox{Pic}(\overline{X_{0,n}})\right)=\mbox{Coker}\big( \mbox{Br}(\mathbb{Q}) \rightarrow \mbox{Br}(X_{0,n}) \big),
\label{beta}
\end{eqnarray}
$\alpha(X_{0,n})$ est le volume d'un certain polytope dans le dual du cône effectif  et $\omega_H\left(X_{0,n}(\mathbb{A}_{\mathbb{Q}})^{{\rm{Br}}(X_{0,n})}\right)$ est un nombre de Tamagawa que l'on détaillera en section 4.3.

\subsubsection{Le facteur $\beta(X_{0,n})$}
En raisonnant comme dans \cite{Blomer2014}, on obtient que le groupe de Brauer cohomologique de $X_{0,n}$, à savoir $\mbox{Br}(X_{0,n})=H^2_{{\footnotesize\mbox{\'et}}}(\overline{X}_0,\mathbb{G}_m)$, est trivial. En effet, c'est un invariant birationnel \cite{Gr} et on combine alors le fait que l'hypersurface $W_n$ considérée soit rationnelle avec le fait que $\mbox{Br}\big(\mathbb{P}^r_{\mathbb{Q}}\big)=\{0\}$ pour tout $r \geqslant 1$ pour obtenir le résultat. Il s'ensuit alors de (\ref{beta}) que $\beta(X_{0,n})=1$.

\subsubsection{Le facteur $\alpha(X_{0,n})$} Dans le but de calculer le facteur $\alpha(X_{0,n})$ apparaissant dans la constante de Peyre, il est nécessaire de décrire un peu plus précisément le géométrie de $X_{0,n}$ et notamment son groupe de Picard et son cône pseudo-effectif. On s'appuie pour ce faire sur \cite{Blomer2014} où ce travail est effectué dans le cas $n=3$ et sur les résultats de Salberger présentés en Annexe. On prouve pour commencer la proposition suivante. Les arguments reposent essentiellement sur le fait que la variété $W_n$ soit la compactification équivariante d'un groupe algébrique.

\begin{prop}
Le morphisme $f_{0,n}: X_{0,n} \rightarrow W_n$ défini dans l'Annexe et en section 4.1 est un morphisme propre, $G_n$-équivariant d'une variété normale $X_{0,n}$ vers $W_n$. De plus, il s'agit d'une résolution crépante de $W_n$.
\end{prop}
\noindent
\textit{Démonstration.--} Le fait qu'il s'agisse d'une résolution crépante résulte du théorème 4 de l'Annexe de Salberger et le fait que $X_{0,n}$ soit normale résulte du fait que $X_{0,n}$ soit lisse. On tire aussi de la construction par Salberger de $f_{0,n}$ en Annexe le fait que $f_{0,n}=p_{X_n} \circ f_n$ avec 
$$
p_{X_n}:X_{0,n} \rightarrow X_n \quad  \mbox{et} \quad f_n:X_n \rightarrow W_n
$$
définie en \cite[section 3]{Blomer2014}. D'après le théorème 6 de \cite{Blomer2014}, il vient que $f_n$ est un morphisme propre et \-$G_n$-équivariant. Il suffit donc de démontrer que $p_{X_n}$ est un morphisme propre et $G_n$-équivariant afin de conclure. Pour ce faire on s'inspire de la preuve de ce théorème 6 de \cite{Blomer2014}. On considère alors 
$$
U_n=j_n(G_n)=\{(x_1,\dots,x_n,y_1,\dots,y_n) \in W_n \hspace{1mm} :\hspace{1mm} y_1\cdots y_n \neq 0\}
$$
pour $j_n:G_n \hookrightarrow W_n$. On a alors $f_{0,n}^{-1}(U_n) \subseteq B^{\ast}_0$ et $f_{0,n}$ est un isomorphisme de $f_{0,n}^{-1}(U_n)$ sur $U_n$. En effet, l'application inverse est donnée par
$$
(x_1,\dots,x_n,y_1,\dots,y_n) \in U_n \longmapsto \left( x_1,\dots,x_n,y_1,\dots,y_n; \prod_{1 \leqslant h \leqslant N}\left(\mathbf{y}^{(h)};\frac{1}{\mathbf{y}^{(h)}}\right)\right) \in f_{0,n}^{-1}(U_n),
$$
où
$$
\mathbf{y}^{(h)}=\left(y_{i_1},\dots,y_{i_{s(h)}}\right) \quad \mbox{pour} \quad \varepsilon_{i_1}(h)=\cdots=\varepsilon_{i_{s(h)}}(h)=1
$$
et
$$
\frac{1}{\mathbf{y}^{(h)}}=\left(\frac{1}{y_{i_1}},\dots,\frac{1}{y_{i_{s(h)}}}\right) \quad \mbox{pour} \quad \varepsilon_{i_1}(h)=\cdots=\varepsilon_{i_{s(h)}}(h)=1.
$$
On peut donc considérer $G_n$ comme un ouvert de $X_{0,n}$ et que $X_{0,n}$ est muni d'une $G_n$-action naturelle $\beta_0: G_n \times X_{0,n} \rightarrow X_{0,n}$ pour laquelle
$$
\left( \begin{pmatrix}
 b_1&a_1\\
 0&b_1\\
 \end{pmatrix}, \dots,\begin{pmatrix}
 b_n&a_n\\
 0&b_n\\
 \end{pmatrix}
 \right). \left(x_1,\dots,x_n,y_1,\dots,y_n; \prod_{1 \leqslant h \leqslant N} \left(\mathbf{Y}^{(h)};\mathbf{Z}^{(h)}\right)  \right)
 $$
est donné par
$$
\left(b_1x_1+a_1y_1,\dots,b_nx_n+a_ny_n,b_1y_1,\dots,b_ny_n;  \prod_{1 \leqslant h \leqslant N}\left(\mathbf{b}\mathbf{Y}^{(h)};\frac{\mathbf{Z}^{(h)}}{\mathbf{b}} \right)\right),
$$
où l'on note
$$
\mathbf{b}\mathbf{Y}^{(h)}=\left(b_{i_1}Y^{(h)}_{i_1},\dots,b_{i_{s(h)}}Y^{(h)}_{i_{s(h)}}\right) \quad \mbox{pour} \quad \varepsilon_{i_1}(h)=\cdots=\varepsilon_{i_{s(h)}}(h)=1
$$
et
$$
\frac{\mathbf{Z}^{(h)}}{\mathbf{b}}=\left(\frac{Z^{(h)}_{i_1}}{b_{i_1}},\dots,\frac{Z^{(h)}_{i_{s(h)}}}{b_{i_{s(h)}}}\right) \quad \mbox{pour} \quad \varepsilon_{i_1}(h)=\cdots=\varepsilon_{i_{s(h)}}(h)=1.
$$
La restriction de $\beta_0$ à $G_n \times  f_{0,n}^{-1}(U_n)$ donne simplement la loi de groupe de $G_n$ et le diagramme suivant 
\begin{center}
\hspace{0.1mm}\xymatrix{
    C_n \times X_{0,n} \ar[r]^{{\hspace{1mm} \beta_0}} \ar[d]_{{({\rm Id},f_{0,n})}} & X_{0,n} \ar[d]^{{f_{0,n}}} \\
    G_n \times W_n \ar[r]_{{\hspace{1mm} \alpha}} & W_n
  }
\end{center}
est commutatif, ce qui permet de compléter la preuve de la proposition.
\hfill
$\square$\\
\newline
\indent
Passons maintenant à une description du groupe de Picard de $X_{0,n}$ et de son cône pseudo-effectif. De la même façon que dans \cite{Blomer2014}, le fait que $X_{0,n}$ soit une compactification équivariante de $G_n$ nous permet d'exploiter le résultat de \cite{TT} suivant.

\begin{prop}
Soient $Y$ une compactification équivariante lisse et propre d'un groupe alg\-ébrique linéaire connexe et résoluble $G$, ${\rm{Div}}_{Y \smallsetminus G}(Y)$ le groupe abélien libre des diviseurs à support dans $Y \smallsetminus G$ et ${\rm C}_{{\rm eff}}(Y) \subseteq {\rm{Pic}}(Y) \otimes_{\mathbb{Z}}\mathbb{R}$ le cône pseudo-effectif engendré par les classes de diviseurs effectifs. Alors la frontière
$$
D=Y \smallsetminus G=\bigcup_{\iota \in I} D_{\iota},
$$
où $I$ indexe les composantes irréductibles de $D$ et $D_{\iota}$ est une composante irréductible de $D$ pour $\iota \in I$, est un diviseur de Weil à croisements normaux. De plus,
\begin{enumerate}[(i)]
\item
On a une suite exacte
$
0 \rightarrow {\rm{Hom}}(G,\mathbb{G}_m) \rightarrow {\rm{Div}}_{Y \smallsetminus G}(Y) \rightarrow {\rm{Pic}}(Y) \rightarrow 0.
$
\item
On a ${\rm{Pic}}(Y)=\bigoplus_{\iota \in I} \mathbb{Z}D_{\iota}$.
\item
On a ${\rm C}_{{\rm eff}}(Y)=\sum_{\iota \in I}\mathbb{R}_{\geqslant 0}D_{\iota}.$
\end{enumerate}
\label{proptt}
\end{prop}
\noindent
\textit{Démonstration.--} La preuve se trouve dans \cite[proposition 1.1]{TT}.\hfill
$\square$\\
\newline
On écrira dans toute la suite $[D]$ pour la classe dans $\mbox{Pic}(X_{0,n})$ d'un diviseur $D$ de $X_{0,n}$, $\left[-K_{X_{0,n}}\right]$ pour la classe anticanonique, $\omega_{X_{0,n}}^{-1}$ pour le faisceau anticanonique et ${\rm C}_{{\rm eff}}(X_{0,n}) \subseteq \mbox{Pic}(X_{0,n}) \otimes_{\mathbb{Z}}\mathbb{R}$ le cône pseudo-effectif engendré par les classes de diviseurs effectifs. On introduit alors $D_N$ la sous-variété de $X_{0,n}$ définie par les équations $y_1=\cdots=y_n=0$ et, pour tout $h \in \llbracket 1,N-1\rrbracket$, on pose également $D_h$ comme étant la sous-variété de $X_{0,n}$ définie par les équations
\begin{eqnarray}
\forall \ell \in \llbracket 1,N\rrbracket \quad \mbox{tel que} \quad \ell \not\preceq h \quad \mbox{et} \quad \ell \not \preceq N-h: \quad 
\left\{
\begin{array}{l}
Y_i^{(\ell)}=0 \quad \mbox{si} \quad \varepsilon_i(\ell)=\varepsilon_i(h)=1\\
Z_j^{(\ell)}=0 \quad \mbox{si} \quad \varepsilon_j(\ell)=1 \quad \mbox{et} \quad \varepsilon_j(h)=0.\\
\end{array}
\right.
\label{div}
\end{eqnarray}
On remarque en particulier que la condition $\ell \not \preceq h$ implique l'existence d'un $i_0 \in \llbracket 1,n \rrbracket$ tel que $\varepsilon_{i_0}(\ell)=1$ et $\varepsilon_{i_0}(h)=0$ et la condition $\ell \not \preceq N-h$ implique l'existence d'un $j_0 \in \llbracket 1,n \rrbracket$ tel que $\varepsilon_{j_0}(\ell)=\varepsilon_{j_0}(h)=1$ si bien que les conditions (\ref{div}) sont bien définies. On démontre alors le lemme essentiel suivant, inspiré du lemme 4 de \cite{Blomer2014}.

\begin{lemme}
\begin{enumerate}[(i)]
\item
On a $X_{0,n} \smallsetminus G_n=\bigcup_{1 \leqslant h \leqslant N} D_h$ et ${\rm{Div}}_{X_{0,n} \smallsetminus G_n}(X_{0,n})=\bigoplus_{1 \leqslant h \leqslant N} \mathbb{Z}D_h$.
\item
Le morphisme canonique de ${\rm{Div}}_{X_{0,n} \smallsetminus G_n}(X_{0,n})$ vers ${\rm{Pic}}(X_{0,n})$ est surjectif de noyau engendré par
$$
\sum_{1 \leqslant h \leqslant N \atop \varepsilon_i(h)=1}D_h-\sum_{1 \leqslant h \leqslant N \atop \varepsilon_1(h)=1}D_h \quad \mbox{pour} \quad 2\leqslant i \leqslant n.
$$
En particulier ${\rm{rang}}\left({\rm{Pic}}(X_{0,n})\right)=2^n-n.$
\item
On a $${\rm C}_{{\rm eff}}(X_{0,n})=\sum_{1 \leqslant h \leqslant N}\mathbb{R}_{\geqslant 0}D_{h}.$$
\item
Pour tout $1 \leqslant i \leqslant n$, le diviseur
$$n\sum_{1 \leqslant h \leqslant N \atop \varepsilon_i(h)=1}D_h$$ est un diviseur anticanonique.
\end{enumerate}
\label{lemmegeo}
\end{lemme}
\noindent
\textbf{Remarque.--} On déduit alors du point $(\textit{ii})$ du Lemme \ref{lemmegeo} (et plus particulièrement du fait que $\mbox{rang}\left(\mbox{Pic}(X_{0,n})\right)=2^n-n$) que la puissance de log obtenue dans le Théorème \ref{theor1} est bien conforme à la conjecture de Manin.\\
\newline
\noindent
\textit{Démonstration.--}
Commençons par le point (\textit{i}). Par définition, on a $$X_{0,n} \smallsetminus G_n=\left\{ \left(x_1,\dots,x_n,y_1,\dots,y_n; \prod_{1 \leqslant h \leqslant N} \left(Y^h;Z^h\right)  \right)
 \in X_{0,n} \hspace{1mm} : \hspace{1mm} y_1\cdots y_n=0\right\}.$$ L'inclusion 
$$
 \bigcup_{1 \leqslant h \leqslant N} D_h \subseteq X_{0,n} \smallsetminus G_n
$$
est claire. Supposons que 
$$
P=\left(x_1,\dots,x_n,y_1,\dots,y_n; \prod_{1 \leqslant h \leqslant N} (\mathbf{Y}^{(h)};\mathbf{Z}^{(h)})  \right) \in X_{0,n} \smallsetminus G_n.
$$
On a alors deux cas. Soit $y_1=\cdots=y_n=0$, auquel cas on a immédiatement $P \in D_N$. Soit il existe $k \in \left\llbracket 1,n-1\right\rrbracket$ tel que $y_{i_1}=\cdots=y_{i_k}=0$ pour $\{i_1,\dots,i_k\}\subseteq \llbracket 1,n\rrbracket$ et $y_j\neq 0$ pour $j \not \in \{i_1,\dots,i_k\}$. Alors, d'après (\ref{8}) et (\ref{7}), on a $Y^{(N)}_{i_1}=\cdots=Y^{(N)}_{i_k}=0$ et $Z^{(N)}_{i}=0$ pour tout $i \not \in \{i_1,\dots,i_k\}$. Grâce aux équations (\ref{4}), (\ref{5}), (\ref{6}), (\ref{7}) et (\ref{8}), on en déduit que $P \in D_h$ pour $h=\sum_{j=1}^k 2^{i_j-1}$. Il suffit alors de constater que chacune des variétés $D_h$ pour $h \in \llbracket 1,N\rrbracket$ est irréductible et on peut alors conclure grâce à la Proposition~\ref{proptt}.\\
\par
Passons au point (\textit{ii}). Le point $(i)$ de la Proposition \ref{proptt} permet d'obtenir la suite exacte
$$
0 \rightarrow \mbox{Hom}(G_n,\mathbb{G}_m) \rightarrow \mbox{Div}_{X_{0,n} \smallsetminus G_n}(X_{0,n}) \rightarrow \mbox{Pic}(X_{0,n}) \rightarrow 0
$$
où l'application $\mbox{Hom}(G_n,\mathbb{G}_m) \rightarrow \mbox{Div}_{X_{0,n} \smallsetminus G_n}(X_{0,n})$ est l'application diviseur des fonctions rationnelles. Comme dans \cite[Lemma 4.(\textit{i})]{Blomer2014}, on obtient que $\mbox{Hom}(G_n,\mathbb{G}_m) $ est le groupe abélien libre engendré par les $\frac{y_i}{y_1}=\frac{Y^{(N)}_i}{Y^{(N)}_1}$ pour $i \in \llbracket 2,n\rrbracket$. Il vient alors facilement que
$$
\forall i \in \llbracket 2,n\rrbracket, \quad \mbox{div}\left( \frac{Y^{(N)}_i}{Y^{(N)}_1} \right)=
\sum_{1 \leqslant h \leqslant N \atop \varepsilon_i(h)=1}D_h-\sum_{1 \leqslant h \leqslant N \atop \varepsilon_1(h)=1}D_h \quad \mbox{pour} \quad 2\leqslant i \leqslant n,
$$
ce qui permet de conclure la démonstration de $(ii)$.\\
\par
Le point (\textit{iii}) résulte immédiatement du point \textit{(iii)} de la Proposition \ref{proptt}.\\
\par
Enfin, démontrons le point (\textit{iv}). Comme $f_{0,n}:X_{0,n} \rightarrow W_n$ est crépante et que $\omega_{W_n} \cong O_{W_n}(-n)$, d'après \cite[II.6.17.1]{Ha}, il suffit d'établir que le sous-schéma fermé défini par $y_i=0$ pour $i \in \llbracket 1,n \rrbracket$ donne lieu au diviseur $$\sum_{1 \leqslant h \leqslant N \atop \varepsilon_i(h)=1}D_h.$$ Comme remarqué dans \cite[lemma 4]{Blomer2014}, il est facile de voir que $y_i$ a multiplicité 1 le long de $D_0$ et on est par conséquent ramené à établir que le sous-schéma fermé défini par $Y_i=0$ donne lieu au diviseur $$\sum_{1 \leqslant h \leqslant N-1 \atop \varepsilon_i(h)=1}D_h,$$ ce qui est clair au vu de (\ref{div}).
\hfill
$\square$\\
\newline
\indent
On est désormais en mesure de calculer, de façon analogue à \cite{Blomer2014}, le facteur $\alpha(X_{0,n})$ défini dans \cite{P95}. On introduit ${\rm C}_{{\rm eff}}(X_{0,n})^{\vee} \subseteq \mbox{Hom}\left( \mbox{Pic}(X_{0,n})\otimes_{\mathbb{Z}}\mathbb{R},\mathbb{R} \right)$, le cône dual de ${\rm C}_{{\rm eff}}(X_{0,n})$ constitué de toutes les applications linéaires $\Lambda: \mbox{Pic}(X_{0,n})^{\vee}\otimes_{\mathbb{Z}} \mathbb{R} \rightarrow \mathbb{R}$ telles que $\lambda\left([D]\right) \geqslant 0$ pour tout diviseur effectif $D$ de $X_{0,n}$. Considérons de plus $\ell:\mbox{Hom}\left( \mbox{Pic}(X_{0,n})\otimes_{\mathbb{Z}}\mathbb{R},\mathbb{R} \right)\rightarrow \mathbb{R}$ l'application linéaire qui à tout $\Lambda \in \mbox{Hom}\left( \mbox{Pic}(X_{0,n})\otimes_{\mathbb{Z}}\mathbb{R},\mathbb{R} \right)$ associe $\Lambda\left(\left[ -K_{X_{0,n}}\right]\right)$. On munit alors $\mbox{Hom}\left( \mbox{Pic}(X_{0,n})\otimes_{\mathbb{Z}}\mathbb{R},\mathbb{R} \right)$ de la mesure de Lebesgue $\mbox{d}s$ normalisée telle que $L=\mbox{Hom}\left( \mbox{Pic}(X_{0,n}),\mathbb{Z} \right)$ soit de covolume 1. On munit ainsi l'hyperplan
$
\mathcal{H}=\ell^{-1}\left( \{1\}\right)
$
de la mesure $\frac{{\rm{d}}s}{{\rm{d}}(\ell-1)}$. En particulier, si $z_1,\dots, z_r$ sont des coordonnées de $\mbox{Hom}\left( \mbox{Pic}(X_{0,n})\otimes_{\mathbb{Z}}\mathbb{R},\mathbb{R} \right)$ correspondant à une $\mathbb{Z}-$base de $L$ et si $\ell(z_1,\dots,z_r)=\alpha_1z_1+\cdots+\alpha_rz_r$ pour $(\alpha_1,\dots,\alpha_r) \in \mathbb{R}^r$, alors 
\begin{eqnarray}
\frac{\mbox{d}s}{\mbox{d}(\ell-1)}=\frac{\mbox{d}z_1\cdots\mbox{d}z_{i-1}\widehat{\mbox{d}z_i}\mbox{d}z_{i+1}\cdots\mbox{d}z_r}{|\alpha_i|}=\frac{\mbox{d}z_1\cdots\mbox{d}z_{i-1}\mbox{d}z_{i+1}\cdots\mbox{d}z_r}{|\alpha_i|} \quad \mbox{dès que } \quad \alpha_i \neq 0.
\label{mesure}
\end{eqnarray}
On a alors, en accord avec \cite{P95},
$$
\alpha(X_{0,n})=\mbox{Vol}\left\{ \Lambda \in {\rm C}_{{\rm eff}}(X_{0,n})^{\vee} \quad | \quad \Lambda\left(\left[ -K_{X_{0,n}} \right]\right)=1 \right\}=\int_{{\rm C}_{{\rm eff}}(X_{0,n})^{\vee} \cap \mathcal{H}} \frac{\mbox{d}s}{\mbox{d}(\ell-1)}.
$$

\begin{lemme}
On a
$
\alpha(X_{0,n})=\frac{1}{n^{2^n-n}}V
$
où $V$ est défini en (\ref{V}).
\label{alpha}
\end{lemme}
\noindent
\textit{Démonstration.--}
Supposons que $n \geqslant 4$, le cas $n=3$ étant couvert par \cite[lemma 5]{Blomer2014}. Grâce au Lemme \ref{lemmegeo}, on sait que $\mbox{Pic}(X_{0,n})$ est engendré par les $D_h$ pour $h \in \llbracket 1,N \rrbracket$ avec les relations
$$
\forall i \neq j \in \llbracket 1,n \rrbracket, \quad \sum_{1 \leqslant h \leqslant N \atop \varepsilon_{i}(h)=1}D_h-\sum_{1 \leqslant h \leqslant N \atop \varepsilon_j(h)=1}D_h=0.
$$
Ces dernières fournissent les relations
\begin{eqnarray}
\forall 4 \leqslant j \leqslant n-1, \quad D_{h_j}=\sum_{\varepsilon_{j+1}(h)=1 \atop \varepsilon_{j}(h)=0 } D_h-\sum_{\substack{\scriptscriptstyle\varepsilon_{j+1}(h)=0 \\\scriptscriptstyle \varepsilon_{j}(h)=1\\\scriptscriptstyle h \neq h_j }} D_h
\label{rel1}
\end{eqnarray}
et
\begin{eqnarray}
D_{5}=\sum_{\varepsilon_{2}(h)=1 \atop \varepsilon_{1}(h)=0 } D_h-\sum_{\substack{\scriptscriptstyle\varepsilon_{2}(h)=0 \\\scriptscriptstyle \varepsilon_{1}(h)=1\\\scriptscriptstyle h \neq 5 }} D_h \label{rel2}
\end{eqnarray}
ainsi que
\begin{eqnarray}
D_3=\sum_{\varepsilon_{3}(h)=1 \atop \varepsilon_{1}(h)=0 } D_h-\sum_{\substack{\scriptscriptstyle \varepsilon_{3}(h)=0 \\\scriptscriptstyle \varepsilon_{1}(h)=1\\\scriptscriptstyle h \neq 3 }} D_h 
\label{rel3}
\end{eqnarray}
et enfin, en utilisant également la relation (\ref{rel2}),
\begin{eqnarray}
\begin{aligned}
D_7&=\sum_{\varepsilon_{4}(h)=1 \atop \varepsilon_{3}(h)=0 } D_h-\sum_{\substack{\scriptscriptstyle\varepsilon_{4}(h)=0 \\\scriptscriptstyle \varepsilon_{3}(h)=1\\\scriptscriptstyle h \neq 7 }} D_h=\sum_{\varepsilon_{4}(h)=1 \atop \varepsilon_{3}(h)=0 } D_h-D_5-\sum_{\substack{\scriptscriptstyle \varepsilon_{4}(h)=0 \\\scriptscriptstyle \varepsilon_{3}(h)=1\\\scriptscriptstyle h \neq 5,7 }} D_h\\
&=\sum_{\varepsilon_{4}(h)=1 \atop \varepsilon_{3}(h)=0 } D_h+\sum_{\substack{\scriptscriptstyle\varepsilon_{2}(h)=0 \\ \scriptscriptstyle\varepsilon_{1}(h)=1\\\scriptscriptstyle h \neq 5 }} D_h-\sum_{\varepsilon_{2}(h)=1 \atop \varepsilon_{1}(h)=0 } D_h-\sum_{\substack{\scriptscriptstyle\varepsilon_{4}(h)=0 \\\scriptscriptstyle \varepsilon_{3}(h)=1\\\scriptscriptstyle h \neq 5,7 }} D_h.
\end{aligned}
\label{rel4}
\end{eqnarray}
On note bien que chacun des membres de droite des relations (\ref{rel1}), (\ref{rel2}), (\ref{rel3}) et (\ref{rel4}) ne fait intervenir aucune des variables $D_{h_j}$ pour $h_j \in H_0 \cup H_1$. On en déduit que $(D_h)_{h \not \in H_0 \cup H_1}$ forme une $\mathbb{Z}$-base de $\mbox{Pic}(X_{0,n})$. On considère alors $(e_h)_{h \not \in H_0 \cup H_1}$ la $\mathbb{Z}$-base duale de $L$ vérifiant $e_h([D_{\ell}])=\delta_{h,\ell}$ pour tous $(h, \ell) \in \llbracket 1,N \rrbracket^2$ et on notera $(z_h)_{h \not \in H_0 \cup H_1}$ les coordonnées de $\mbox{Hom}\left( \mbox{Pic}(X_{0,n})\otimes_{\mathbb{Z}}\mathbb{R},\mathbb{R} \right)$ relatives à cette base. Le Lemme \ref{lemmegeo} entraîne que 
$$
{\rm C}_{{\rm eff}}(X_{0,n})^{\vee}=\left\{ (z_h)_{h \not \in H_0 \cup H_1} \in \left(\mathbb{R}_{\geqslant 0}\right)^{2^n-n} \hspace{1mm} : \hspace{1mm}
\begin{array}{c}
\displaystyle\underset{\substack{\scriptscriptstyle\varepsilon_{j+1}(h)=1 \\ \scriptscriptstyle\varepsilon_{j}(h)=0 }}{\sum} z_h-\underset{\substack{\scriptscriptstyle\varepsilon_{j+1}(h)=0 \\\scriptscriptstyle \varepsilon_{j}(h)=1\\\scriptscriptstyle h \neq h_j }}{\sum} z_h \geqslant 0, \quad 4 \leqslant j \leqslant n-1,\\[5mm]
\displaystyle\underset{\varepsilon_{2}(h)=1 \atop \varepsilon_{1}(h)=0 }{\sum} z_h-\underset{\substack{\scriptscriptstyle\varepsilon_{2}(h)=0 \\\scriptscriptstyle \varepsilon_{1}(h)=1\\\scriptscriptstyle h \neq 5 }}{\sum} z_h \geqslant 0,\hspace{1mm} \underset{\varepsilon_{3}(h)=1 \atop \varepsilon_{1}(h)=0 }{\sum} z_h-\underset{\substack{\scriptscriptstyle\varepsilon_{3}(h)=0 \\ \scriptscriptstyle\varepsilon_{1}(h)=1\\\scriptscriptstyle h \neq 3 }}{\sum} z_h \geqslant 0,\\[5mm]
\displaystyle\underset{\varepsilon_{4}(h)=1 \atop \varepsilon_{3}(h)=0 }{\sum} z_h+\underset{\substack{\scriptscriptstyle\varepsilon_{2}(h)=0 \\\scriptscriptstyle \varepsilon_{1}(h)=1\\\scriptscriptstyle h \neq 5 }} {\sum}z_h-\underset{\varepsilon_{2}(h)=1 \atop \varepsilon_{1}(h)=0 }{\sum} z_h-\underset{\substack{\scriptscriptstyle\varepsilon_{4}(h)=0 \\\scriptscriptstyle \varepsilon_{3}(h)=1\\\scriptscriptstyle h \neq 5,7 }}{\sum} z_h \geqslant 0.
\end{array}
 \right\}
$$
et le point (\textit{iv}) du Lemme \ref{lemmegeo} entraîne que $\mathcal{H}$ a pour équation $$n\sum_{\varepsilon_n(h)=1} z_h=0,$$ équation qui ne fait pas intervenir aucune variable $z_h$ pour $h \in H_0 \cup H_1$. D'après (\ref{mesure}) et en éliminant $z_N$, on obtient
$$
\frac{\mbox{d}s}{\mbox{d}(\ell-1)}=\frac{\displaystyle\underset{\scriptscriptstyle 1 \leqslant h \leqslant N-1}{\prod} \mbox{d}z_h}{n}
$$
si bien que
$$
\alpha(X_{0,n})=\frac{1}{n}\int_{\Delta} \prod_{1 \leqslant h \leqslant N-1}\mbox{d}z_h
$$
où
$$
\Delta=\left\{ (z_h)_{h \not \in H_0 \cup H_1 \cup H_2} \in \left(\mathbb{R}_{\geqslant 0}\right)^{2^n-n-1} \hspace{1mm} : \hspace{1mm}
\begin{array}{l}
\displaystyle n\underset{\substack{\scriptscriptstyle\varepsilon_n(h)=1 \\\scriptscriptstyle h \neq h_n}}{\sum} z_h \leqslant 1,\\[5mm]
\displaystyle\underset{\substack{\scriptscriptstyle\varepsilon_{j+1}(h)=1 \\\scriptscriptstyle \varepsilon_{j}(h)=0 }}{\sum} z_h-\underset{\substack{\scriptscriptstyle\varepsilon_{j+1}(h)=0 \\\scriptscriptstyle \varepsilon_{j}(h)=1\\\scriptscriptstyle h \neq h_j }}{\sum} z_h \geqslant 0 \qquad 4 \leqslant j \leqslant n-1,\\[5mm]
\displaystyle\underset{\scriptscriptstyle\varepsilon_{2}(h)=1 \atop \scriptscriptstyle\varepsilon_{1}(h)=0 }{\sum} z_h-\underset{\substack{\scriptscriptstyle\varepsilon_{2}(h)=0 \\\scriptscriptstyle \varepsilon_{1}(h)=1\\\scriptscriptstyle h \neq 5 }}{\sum} z_h \geqslant 0,\quad \underset{\scriptscriptstyle\varepsilon_{3}(h)=1 \atop\scriptscriptstyle \varepsilon_{1}(h)=0 }{\sum} z_h-\underset{\substack{\scriptscriptstyle\varepsilon_{3}(h)=0 \\\scriptscriptstyle \varepsilon_{1}(h)=1\\\scriptscriptstyle h \neq 3 }}{\sum} z_h \geqslant 0,\\[5mm]
\displaystyle\underset{\scriptscriptstyle\varepsilon_{4}(h)=1 \atop \scriptscriptstyle\varepsilon_{3}(h)=0 }{\sum} z_h+\underset{\substack{\scriptscriptstyle\varepsilon_{2}(h)=0 \\\scriptscriptstyle \varepsilon_{1}(h)=1\\\scriptscriptstyle h \neq 5 }} {\sum}z_h-\underset{\scriptscriptstyle\varepsilon_{2}(h)=1 \atop\scriptscriptstyle \varepsilon_{1}(h)=0 }{\sum} z_h-\underset{\substack{\scriptscriptstyle\varepsilon_{4}(h)=0 \\\scriptscriptstyle \varepsilon_{3}(h)=1\\\scriptscriptstyle h \neq 5,7 }}{\sum} z_h \geqslant 0.
\end{array}
 \right\}.
$$
Le changement de variables $t_h=n z_h$ pour tout $h \not\in H_0 \cup H_1 \cup H_2$ fournit alors immédiatement le résultat escompté.
\hfill
$\square$

\subsubsection{Construction du torseur versel associé à $W_n$}

Soit $n \geqslant 3$. Si l'on note 
\begin{eqnarray}
\mathcal{W}_n=\#\left\{ (x_1,\dots,x_n,y_1,\dots,y_n) \in \mathbb{Z}^{2n} \hspace{1mm} : \hspace{1mm} \begin{array}{c}
 y_1y_2 \cdots y_n \neq 0,\\
(\mathbf{x},\mathbf{y}) \hspace{2mm} {\rm{ v\acute{e}rifie }} \hspace{2mm} (\ref{eq})
 \end{array}
 \right\}
\label{wn}
\end{eqnarray}
et
$$
\mathcal{A}_n=\#\left\{ \left(\mathbf{x}';\left(z_h\right)_{1 \leqslant h \leqslant N}\right) \in \mathbb{Z}^n \times \left(\mathbb{Z}\smallsetminus\{0\}\right)^{2^n-1} \hspace{1mm} : \hspace{1mm}   \begin{array}{c}
 \left(z_h\right)_{1 \leqslant h \leqslant N} \hspace{2mm} N\mbox{-uplet réduit},\\[1mm]
 z_h > 0 \hspace{2mm} {\rm si} \hspace{2mm} s(h) \geqslant 2,\\[1mm]
\displaystyle\sum_{j=1}^n \Bigg(\prod_{1 \leqslant h \leqslant N} z_h^{1-\varepsilon_j(h)}\Bigg)z_{2^{j-1}}x'_j=0
 \end{array}
\right\},
$$
la section 3.2.2 et en particulier la relation (\ref{torsor}) permet d'obtenir le lemme suivant, qui coïncide avec \cite[Lemma 7]{Blomer2014} dans le cas $n=3$. 
\begin{lemme}
L'application $\mathcal{A}_n \rightarrow \mathbb{Z}^{2n}$ définie par
\begin{eqnarray}
\left(\mathbf{x}';\left(z_h\right)_{1 \leqslant h \leqslant N}\right) \in \mathcal{A}_n \longmapsto \left(z_1x'_1,\dots,z_{2^{n-1}}x'_n,\prod_{1 \leqslant h \leqslant N} z_h^{\varepsilon_1(h)},\dots,\prod_{1 \leqslant h \leqslant N} z_h^{\varepsilon_n(h)}  \right)
\label{bij2}
\end{eqnarray}
réalise une bijection de $\mathcal{A}_n$ sur $\mathcal{W}_n$.
\label{bij}
\end{lemme}

C'est ce paramètrage de $W_n$ qui nous a permis d'obtenir le Théorème \ref{theor1}. On explicite maintenant le torseur versel $\mathcal{T}$ associé à la résolution crépante $X_{0,n}$ explicitée en section 4.2 et on établit de façon analogue à \cite{Blomer2014} que ce paramètrage du problème de comptage est en réalité une descente sur ce torseur versel $\mathcal{T}$. La constante de Peyre s'interprétera alors à la manière de \cite{Sal} et \cite{Blomer2014} comme un volume adélique de $\mathcal{T}$. On renvoie à la section 4.2 de \cite{Blomer2014} et au Lemme \ref{lemmegeo} pour un rappel concernant la notion de torseur versel et la justification qu'il n'existe, à isomorphisme près, qu'un seul torseur versel pour la variété $X_{0,n}$.\\
\par
On suit ici le schéma de la démonstration dans le cas $n=3$ de \cite{Blomer2014}. La variété $X_{0,n}$ est une hypersurface de la variété $\Xi_0 \subseteq \mathbb{P}^{2n-1} \times \prod_{1 \leqslant h \leqslant N} \mathbb{P}^{(h)}\times \mathbb{P}^{(h)}$ définie par les équations (\ref{4}), (\ref{5}), (\ref{6}) et (\ref{8}). De la même façon que pour $X_{0,n}$, il est facile de voir que la restriction à $\Xi_0$ de la projection $\mathbb{P}^{2n-1} \times \prod_{1 \leqslant h \leqslant N} \mathbb{P}^{(h)}\times \mathbb{P}^{(h)} \rightarrow \prod_{1 \leqslant h \leqslant N} \mathbb{P}^{(h)}\times \mathbb{P}^{(h)}$ donne lieu à un morphisme $\gamma:\Xi_0 \rightarrow B_{0,n}$ qui munit la variété $\Xi_0$ d'une structure de $\mathbb{P}^n$-fibré sur $B_{0,n}$. On commence, comme dans \cite{Blomer2014}, à décrire le torseur versel de $\Xi_0$. Pour ce faire, on tire parti du lemme fondamental suivant.

\begin{lemme}
La variété $\Xi_0$ est une variété torique projective lisse et déployée. Le tore associé est l'ouvert $U$ de $\Xi_0$ défini par
$$
x_1 \cdots x_ny_1 \cdots y_n \prod_{1 \leqslant h \leqslant N} \prod_{1 \leqslant i \leqslant n \atop \varepsilon_i(h)=1} Y^{(h)}_i Z^{(h)}_i\neq 0.
$$
\label{torique}
\end{lemme}
\noindent
\textit{Démonstration.--} On considère $T_1$ le tore déployé de dimension $2n-1$ défini comme le quotient du tore $\mathbb{G}^{2n}_m$ par le plongement diagonal de $\mathbb{G}_m$ dans $\mathbb{G}_m^{2n}.$ Lorsque la multiplication sur $U$ est donnée par la multiplication coordonnée par coordonnée, l'application 
$$
\left(x_1,\dots,x_n,y_1,\dots,y_n; \prod_{1 \leqslant h \leqslant N} \left(\mathbf{Y}^{(h)};\mathbf{Z}^{(h)}\right)  \right) \in U \longmapsto \left[(x_1,\dots,x_n,y_1,\dots,y_n)\right] \in T_1
$$
est un isomorphisme de réciproque
$$
\left[(t_1,\dots,t_n,u_1,\dots,u_n)\right] \in T_1 \longmapsto \left(t_1,\dots,t_n,u_1,\dots,u_n; \prod_{1 \leqslant h \leqslant N} \left(\mathbf{U}^{(h)};\frac{1}{\mathbf{U}^{(h)}}\right)  \right) \in U,
$$
avec
$$
\mathbf{U}^{(h)}=\left(u_{i_1},\dots,u_{i_{s(h)}}\right) \quad \mbox{pour} \quad \varepsilon_{i_1}(h)=\cdots=\varepsilon_{i_{s(h)}}(h)=1
$$
et
$$
\frac{1}{\mathbf{U}^{(h)}}=\left(\frac{1}{u_{i_1}},\dots,\frac{1}{u_{i_{s(h)}}}\right) \quad \mbox{pour} \quad \varepsilon_{i_1}(h)=\cdots=\varepsilon_{i_{s(h)}}(h)=1.
$$
Le tore $\mathbb{G}_m^{2n}$ agit sur $\Xi_0$ par multiplication coordonnée par coordonnée. Le produit
$$
(t_1,\dots,t_n,u_1,\dots,u_n).\left(x_1,\dots,x_n,y_1,\dots,y_n; \prod_{1 \leqslant h \leqslant N} \left(\mathbf{Y}^{(h)};\mathbf{Z}^{(h)}\right)  \right)
$$
est donné par
$$
\left(t_1x_1,\dots,t_nx_n,u_1y_1,\dots,u_ny_n; \prod_{1 \leqslant h \leqslant N} \left(\mathbf{u}\mathbf{Y}^{(h)};\frac{\mathbf{Z}^{(h)}}{\mathbf{u}}\right)  \right),
$$
avec
$$
\mathbf{u}\mathbf{Y}^{(h)}=\left(u_{i_1}Y^{(h)}_{i_1},\dots,u_{i_{s(h)}}Y^{(h)}_{i_{s(h)}}\right) \quad \mbox{pour} \quad \varepsilon_{i_1}(h)=\cdots=\varepsilon_{i_{s(h)}}(h)=1
$$
et
$$
\frac{\mathbf{u}}{\mathbf{Z}^{(h)}}=\left(\frac{Z^{(h)}_{i_1}}{u_{i_1}},\dots,\frac{Z^{(h)}_{i_{s(h)}}}{u_{i_{s(h)}}}\right) \quad \mbox{pour} \quad \varepsilon_{i_1}(h)=\cdots=\varepsilon_{i_{s(h)}}(h)=1.
$$
Cela fournit après passage au quotient une action $\rho: T \times \Xi_0 \rightarrow \Xi_0$ dont la restriction à $U$ donne la loi de groupe de $T$ et cela permet de conclure que $\Xi_0$ est un variété torique projective de dimension $2n-1$. Elle est lisse puisqu'il s'agit d'un $\mathbb{P}^n$-fibré sur $B_{0,n}$.\hfill
$\square$\\
\newline
\indent
Comme dans \cite{Blomer2014}, on utilise alors, le résultat de \cite{Sal} suivant qui permet de décrire le torseur versel d'une variété torique déployée projective lisse. Le torseur versel $\mathcal{T}$ de $\Xi_0$ est ainsi donné par le morphisme de Cox de la sous-variété torique $\mathbb{A}^n \smallsetminus F$ de $\mathbb{A}^n$ dans $\Xi_0$ décrit dans \cite{Co} où $n$ est le nombre de cônes de dimension un (appelés arêtes) de l'éventail $\Delta$ de $\Xi_0$ (voir \cite{Fu} pour plus de détails à ce sujet) et $F \subseteq \mathbb{A}^n$ est le sous-ensemble fermé défini par les monômes $t^{\sigma}=\prod_{\rho \not\in \sigma(1)} t_{\rho}$ où $t_{\rho}$ est un système de coordonnées de $\mathbb{A}^n$ lorsque $\rho$ décrit l'ensemble des arêtes de $\Delta$, $\sigma \in \Delta$ un cône maximal et $\sigma(1)$ est l'ensemble des arêtes de $\sigma$. On applique alors ce résultat pour en déduire la proposition suivante.

\begin{prop}
Soit $\Omega \subset \mathbb{A}^{2^n+n-1}$ la sous-variété ouverte  de coordonnées $\left(\mathbf{x};\left(z_h\right)_{1 \leqslant h \leqslant N}\right)$ définie par les conditions
\begin{eqnarray}
z_{h_n}\prod_{1 \leqslant h \leqslant N \atop h \not \in \mathcal{H}_{\mathbf{j}}}z_h\neq 0 \quad \mbox{ou} \quad x_i\prod_{1 \leqslant h \leqslant N \atop h \not \in \mathcal{H}_{\mathbf{j}}}z_h\neq 0
\label{condd1}
\end{eqnarray}
pour $i \in \llbracket 1,n \rrbracket$ et $\mathbf{j}=(j_1,\dots,j_n)$ tel que $\{j_1,\dots,j_n\}=\{1,\dots,n\}$ et 
$$
\mathcal{H}_{\mathbf{j}}=\left\{ h_{1,{\mathbf{j}}}=2^{j_1-1}\preceq \cdots \preceq h_{k,{\mathbf{j}}}=\sum_{s=1}^k 2^{j_s-1} \preceq \cdots \preceq h_{n,{\mathbf{j}}}=N \right\}.
$$
Le morphisme $\varphi: \Omega \rightarrow \Xi_0$ défini par
\begin{eqnarray}
\left(\mathbf{x};\left(z_h\right)_{1 \leqslant h \leqslant N}\right) \longmapsto \left( x_1,\dots,x_n,\prod_{1 \leqslant h \leqslant N \atop \varepsilon_1(h)=1}z_h,\dots, \prod_{1 \leqslant h \leqslant N \atop \varepsilon_n(h)=1}z_h;\prod_{1 \leqslant h \leqslant N} \left(\mathbf{Y}^{(h)} ; \mathbf{Z}^{(h)}\right)\right)
\label{app}
\end{eqnarray}
où
$$
\mathbf{Y}^{(h)}=\left( \frac{\displaystyle\underset{\scriptscriptstyle 1 \leqslant \ell \leqslant N \atop\scriptscriptstyle \varepsilon_{i_1}(\ell)=1}{\prod}z_{\ell}}{\displaystyle\underset{\scriptscriptstyle 1 \leqslant \ell \leqslant N \atop\scriptscriptstyle h \preceq \ell}{\prod}z_{\ell}},\dots,\frac{\displaystyle\underset{\scriptscriptstyle 1 \leqslant \ell \leqslant N \atop\scriptscriptstyle \varepsilon_{i_{s(h)}}(\ell)=1}{\prod}z_{\ell}}{\displaystyle\underset{\scriptscriptstyle 1 \leqslant \ell \leqslant N \atop\scriptscriptstyle h \preceq \ell}{\prod}z_{\ell}} \right) \quad \mbox{pour} \quad \varepsilon_{i_1}(h)=\cdots=\varepsilon_{i_{s(h)}}(h)=1
$$
et
$$
\mathbf{Z}^{(h)}= \left(\frac{\displaystyle\underset{\scriptscriptstyle 1 \leqslant \ell \leqslant N \atop\scriptscriptstyle \varepsilon_{i_1}(\ell)=0}{\prod}z_{\ell}}{\displaystyle\underset{\scriptscriptstyle 1 \leqslant \ell \leqslant N \atop\scriptscriptstyle \varepsilon_{i}(\ell)=0 \hspace{1mm} {\rm si} \hspace{1mm} \varepsilon_i(h)=1}{\prod}z_{\ell}},\dots,\frac{\displaystyle\underset{\scriptscriptstyle 1 \leqslant \ell \leqslant N \atop\scriptscriptstyle \varepsilon_{i_{s(h)}}(\ell)=0}{\prod}z_{\ell}}{\displaystyle\underset{\scriptscriptstyle 1 \leqslant \ell \leqslant N \atop\scriptscriptstyle \varepsilon_{i}(\ell)=0 \hspace{1mm} {\rm si} \hspace{1mm} \varepsilon_i(h)=1}{\prod}z_{\ell}} \right) \quad \mbox{pour} \quad \varepsilon_{i_1}(h)=\cdots=\varepsilon_{i_{s(h)}}(h)=1
$$
est alors le morphisme sous-jacent d'un torseur versel pour $\Xi_0$.
\label{torseurversel}
\end{prop}
\noindent
\textit{Démonstration.--} On renvoie à \cite{Blomer2014} et \cite{Fu} pour les trois résultats rappelés ci-dessous. Soient $\Delta$ l'éventail associé à la variété torique $\Xi_0$ et $\sigma$ un cône de $\Delta$. On notera alors $O_{\sigma}$ l'orbite associée à $\sigma$ et $V(\sigma)=\overline{O_{\sigma}}$. On a alors les trois résultats suivants qui vont s'avérer essentiels.\\
\begin{itemize}
\item[\textit{(i)}] Tout d'abord, il existe une bijection entre les arêtes $\rho$ de $\Delta$ et les composantes irréductibles de $\Xi_0 \smallsetminus U$.
\item[\textit{(ii)}] Ensuite, il existe également une bijection entre les cônes maximaux $\sigma$ de $\Delta$ et les points fixes de $\Xi_0 \smallsetminus U$ sous l'action de $U$. 
\item[\textit{(iii)}] Enfin, pour $\sigma \in \Delta$ et $\rho \in \Delta$ une arête, on a l'équivalence suivante:
$$
\rho \in \sigma(1) \quad \Longleftrightarrow \quad V(\sigma) \subseteq D_{\rho} \quad \mbox{où } D_{\rho} \mbox{ est la composante irréductible de } \Xi_0 \smallsetminus U \mbox{ associée à } \rho.
$$
\end{itemize}
\vspace{0.3cm}
\indent
Dans le cas de $\Xi_0$, les composantes irréductibles de $\Xi_0 \smallsetminus U$ s'obtiennent de la même façon que dans la démonstration du Lemme \ref{lemmegeo}. On a $2^n+n-1$ composantes irréductibles. Les $n$ premières que l'on notera $X_i$ dans la suite et à qui on associera la coordonnée $\xi_i$ sont définies par l'équation $x_i=0$ pour $i \in \llbracket 1,n \rrbracket$. On a également $D_N$, associé à la coordonnée $z_N$, définie par l'équation $y_1=\cdots=y_n=0$. Enfin, pour tout $h\in \llbracket 1,N-1 \rrbracket$, on a $D_h$, associé à la coordonnée $z_h$, définie par les équations
$$
\forall \ell \in \llbracket 1,N\rrbracket \quad \mbox{tel que} \quad \ell \not\preceq h \quad \mbox{et} \quad \ell \not \preceq N-h: \quad 
\left\{
\begin{array}{l}
Y_i^{(\ell)}=0 \quad \mbox{si} \quad \varepsilon_i(\ell)=\varepsilon_i(h)=1\\
Z_j^{(\ell)}=0 \quad \mbox{si} \quad \varepsilon_j(\ell)=1 \quad \mbox{et} \quad \varepsilon_j(h)=0.\\
\end{array}
\right.
$$
On déduit alors de \cite{Sal} et de la discussion \textit{supra} \textit{(i)} qu'il existe un plongement naturel du torseur versel $\mathcal{T}$ de $\Xi_0$ dans l'espace affine $\mathbb{A}^{2^n+n-1}$ de coordonnées $\left(\boldsymbol{\xi};\left(z_h\right)_{1 \leqslant h \leqslant N}\right)$.\\
\indent
Au vue de l'action de $U$ sur $\Xi_0$, un point $P$ de $\Xi_0\smallsetminus U$ est fixé par cette action si, et seulement si, son image a exactement une coordonnée non nulle par chacune des projections suivantes
$$
\begin{array}{ccc}
\displaystyle\mathbb{P}^{2n-1} \times \prod_{1 \leqslant h \leqslant N} \displaystyle\mathbb{P}^{(h)}\times \mathbb{P}^{(h)} &\longrightarrow & \mathbb{P}^{2n-1}\\
\displaystyle\left(x_1,\dots,x_n,y_1,\dots,y_n; \prod_{1 \leqslant h \leqslant N} (\mathbf{Y}^{(h)};\mathbf{Z}^{(h)})  \right)& \longmapsto & (x_1,\dots,x_n,y_1,\dots,y_n)\\ 
\end{array}
$$
et pour tout $h \in \llbracket 1,N \rrbracket$
$$
\begin{array}{ccc}
\displaystyle\mathbb{P}^{2n-1} \times \prod_{1 \leqslant h \leqslant N} \displaystyle\mathbb{P}^{(h)}\times \mathbb{P}^{(h)} &\longrightarrow & \mathbb{P}^{(h)}\\
\displaystyle\left(x_1,\dots,x_n,y_1,\dots,y_n; \prod_{1 \leqslant h \leqslant N} (\mathbf{Y}^{(h)};\mathbf{Z}^{(h)})  \right)& \longmapsto & \left(\mathbf{Y}^{(h)}\right)\\ 
\end{array}
$$
et
$$
\begin{array}{ccc}
\displaystyle\mathbb{P}^{2n-1} \times \prod_{1 \leqslant h \leqslant N} \displaystyle\mathbb{P}^{(h)}\times \mathbb{P}^{(h)} &\longrightarrow & \mathbb{P}^{(h)}\\
\displaystyle\left(x_1,\dots,x_n,y_1,\dots,y_n; \prod_{1 \leqslant h \leqslant N} (\mathbf{Y}^{(h)};\mathbf{Z}^{(h)})  \right)& \longmapsto & \left(\mathbf{Z}^{(h)}\right).\\ 
\end{array}
$$
On en déduit que pour point $\Xi_0\smallsetminus U$ fixé par l'action de $U$, soit il existe $j_n \in \llbracket 1,n \rrbracket$ tel que $y_{j_n} \neq 0$ soit il existe $\ell \in \llbracket 1,n \rrbracket$ tel que $x_{\ell} \neq 0$. On admet à présent dans la suite de ce paragraphe que dès qu'un indice $i \in \llbracket 1,n \rrbracket$ a été fixé tel que $Y^{(h)}_i\neq 0$ pour un certain $1 \leqslant h \leqslant N$ donné, alors pour tous les autres indices $j$ tels que $\varepsilon_j(h)=1$, on a $Y^{(h)}_j= 0$.\\
\indent
Plaçons-nous pour commencer dans le cas où $y_{j_n} \neq 0$. On a alors nécessairement $Y^{(N)}_{j_n} \neq 0$ d'après (\ref{8}) et plus généralement, pour tout $h \in \llbracket 1,N \rrbracket$ tel que $\varepsilon_{j_n}(h)=1$, on a $Y^{(h)}_{j_n} \neq 0$ d'après (\ref{5}). De plus, il existe $\ell_n \in \llbracket 1,n \rrbracket \smallsetminus \{j_n\}$ tel que $Z^{(N)}_{\ell_n} \neq 0$ et plus généralement, pour tout $h \in \llbracket 1,N \rrbracket$ tel que $\varepsilon_{\ell_n}(h)=1$, on a $Z^{(h)}_{\ell_n} \neq 0$ au vu de (\ref{4}) et de (\ref{6}). On a donc que pour tout $h \in \llbracket 1,N \rrbracket$ tel que $\varepsilon_{\ell_n}(h)=\varepsilon_{j_n}(h)=1$, $Y^{(h)}_{j_n} \neq 0$ et $Z^{(h)}_{\ell_n} \neq 0$. Intéressons-nous à présent aux $h \in \llbracket 1,N \rrbracket$ tel que $\varepsilon_{\ell_n}(h)=1$ et $\varepsilon_{j_n}(h)=0$. Pour ces $h$, on a nécessairement $Z^{(h)}_{\ell_n} \neq 0$. On pose alors $h_{n-1,{\mathbf{j}}}=N-2^{j_n-1}$. Puisque $P$ est fixé sous l'action de $U$, il existe $j_{n-1} \in \llbracket 1,n \rrbracket\smallsetminus\{ j_n,\ell_n\}$ tel que $Y^{(h_{n-1,{\mathbf{j}}})}_{j_{n-1}} \neq 0$ et plus généralement pour tout $h \in \llbracket 1,N \rrbracket$ tel que $\varepsilon_{\ell_n}(h)=1$ et $h \preceq h_{n-1,\mathbf{j}}$, on a $Y^{(h)}_{j_{n-1}} \neq 0$ d'après (\ref{4}) et (\ref{5}). On construit alors en itérant ce procédé $h_{2,\mathbf{j}}, \dots,h_{n-1,\mathbf{j}}$ ainsi que $j_1,\dots,j_{n-1}$ tel que pour tout $r \in \llbracket 2,n-1\rrbracket$,
$$
h_{r,\mathbf{j}}=N-\sum_{s=r+1}^{n}2^{j_s-1}=\sum_{s=1}^r 2^{j_s-1}
$$
et
$$
Y^{(h)}_{j_r} \neq 0 \quad \forall h \in \llbracket 1,N \rrbracket \quad \mbox{tel que} \quad h \preceq h_{r,\mathbf{j}} \quad \mbox{et} \quad \varepsilon_{j_r}(h)=1
$$
et $j_r \in \llbracket 1,n \rrbracket \smallsetminus \{ j_n, \dots, j_{r+1}, \ell_n \}$ (ce qui est bien toujours possible). Passons à présent aux cas des $h \in \llbracket 1,N \rrbracket$ tels que $\varepsilon_{\ell_n}(h)=0$ et $\varepsilon_{j_n}(h)=1$. Pour ces $h$, on a nécessairement $Y^{(h)}_{j_n} \neq 0$ et on note que tout $h \in \llbracket 1,N \rrbracket$ tel que $\varepsilon_{\ell_n}(h)=0$ et $\varepsilon_{j_n}(h)=1$ vérifie $h \not \preceq \ell$ et $\ell \not \preceq h$ pour tout $\ell \in \llbracket 1,N \rrbracket$ tel que $\varepsilon_{j_n}(\ell)=0$ et $\varepsilon_{j_{n-1}}(\ell)=1$. On pose alors ${h'_{n-1,\boldsymbol{\ell}}}=N-2^{\ell_n-1}$. Il existe ainsi $\ell_{n-1} \in \llbracket 1,n \rrbracket\smallsetminus\{ j_n,\ell_n\}$ tel que $Z^{({h'_{n-1,\boldsymbol{\ell}}})}_{\ell_{n-1}} \neq 0$ et plus généralement pour tout $h \in \llbracket 1,N \rrbracket$ tel que $\varepsilon_{\ell_n}(h)=0$ et $h \preceq {h'_{n-1,\boldsymbol{\ell}}}$, on a $Z^{(h)}_{\ell_{n-1}} \neq 0$. On construit alors ${h'_{3,\boldsymbol{\ell}}}, \dots,{h'_{n-2,\boldsymbol{\ell}}}$ ainsi que $\ell_3,\dots,\ell_{n-1}$ tel que pour tout $r \in \llbracket 3,n-1\rrbracket$,
$$
h'_{r,\boldsymbol{\ell}}=N-\sum_{s=r+1}^{n}2^{\ell_s-1}
$$
et
$$
Z^{(h)}_{\ell_r} \neq 0 \quad \forall h \in \llbracket 1,N \rrbracket \quad \mbox{tel que} \quad h \preceq h'_{r,\boldsymbol{\ell}} \quad \mbox{et} \quad \varepsilon_{\ell_r}(h)=1
$$
et $j_r \in \llbracket 1,n \rrbracket \smallsetminus \{ \ell_n, \dots, \ell_{r+1}, j_n, j_{n-1} \}$. Ensuite si on pose $h'_{2,\boldsymbol{\ell}}=h'_{3,\boldsymbol{\ell}}-2^{\ell_{n-2}-1}$, alors on a $Y^{(h'_{2,\boldsymbol{\ell}})}_{j_n} \neq 0$ et $Z^{(h'_{2,\boldsymbol{\ell}})}_{j_{n-1}} \neq 0$ car les chiffres de $h'_{2,\boldsymbol{\ell}}$ sont exactement $j_n$ et $j_{n-1}$. Reste alors à traiter le cas des $h \in \llbracket 1,N \rrbracket$ tel que $\varepsilon_{\ell_n}(h)=0$ et $\varepsilon_{j_n}(h)=0$. Puisqu'on a $h_{2,\mathbf{j}} \preceq h_{3,\mathbf{j}} \preceq \cdots \preceq h_{n-1,\mathbf{j}}$ et $h'_{2,\boldsymbol{\ell}} \preceq h'_{3,\boldsymbol{\ell}} \preceq \cdots \preceq h'_{n-1,\boldsymbol{\ell}}$, il existe $r$ et $s$ tels que $h \preceq h_{r,\mathbf{j}}$ et $h \preceq h'_{s,\boldsymbol{\ell}}$ et $h \not \preceq h_{i,\mathbf{j}}$ pour $i >r$ et $h \not \preceq h'_{j,\boldsymbol{\ell}}$ pour $j>s$. Il s'ensuit alors que $Y^{(h)}_{j_r} \neq 0$ et $Z^{(h)}_{\ell_s} \neq 0$. Par construction, $h_{i,\mathbf{j}}$ et $h'_{j,\boldsymbol{\ell}}$ ne sont pas comparables pour la relation $\preceq$ pour $i \geqslant 3$ et $j \geqslant 3$. 
On en déduit finalement que $P$ n'appartient pas aux composantes irréductibles $D_h$ avec $h \in \llbracket 1,N-1\rrbracket$ telles que
\begin{eqnarray}
\left\{
\begin{array}{l}
\varepsilon_{j_n}(h)=1 \\
\mbox{ou} \quad h \prec h_{n-1,\mathbf{j}} \quad \mbox{et} \quad \varepsilon_{j_{n-1}}(h)=1\\
\mbox{ou} \quad h \prec h_{n-2,\mathbf{j}} \quad \mbox{et} \quad \varepsilon_{j_{n-2}}(h)=1\\
\quad \quad \quad \quad \vdots \\
\mbox{ou} \quad h \prec h_{2,\mathbf{j}} \quad \mbox{et} \quad \varepsilon_{j_{2}}(h)=1\\
\end{array}
\right. 
\label{cc}
\end{eqnarray}
ou
\begin{eqnarray}
\left\{
\begin{array}{l}
\varepsilon_{\ell_n}(h)=0 \\
\mbox{ou} \quad h \prec h'_{n-1,\boldsymbol{\ell}} \quad \mbox{et} \quad \varepsilon_{\ell_{n-1}}(h)=0\\
\mbox{ou} \quad h \prec h'_{n-2,\boldsymbol{\ell}} \quad \mbox{et} \quad \varepsilon_{\ell_{n-2}}(h)=0\\
\vdots \\
\mbox{ou} \quad h \prec h'_{3,\boldsymbol{\ell}} \quad \mbox{et} \quad \varepsilon_{\ell_{3}}(h)=0.\\
\end{array}
\right.
\label{ccc}
\end{eqnarray}
Les seuls $h \in \llbracket 1,n \rrbracket \smallsetminus\{ h_{2,\mathbf{j}},\dots,h_{n-1,\mathbf{j}}\}$ qui ne vérifient pas les conditions (\ref{cc}) sont les $h$ tels que $\varepsilon_{j_2}(h)=\varepsilon_{j_2}(h)=\cdots=\varepsilon_{j_{n}}(h)=0$. En effet, les conditions $\varepsilon_{j_n}(h)=\cdots=\varepsilon_{j_{n-k+1}}(h)=0$ pour $k \in \llbracket 1,n-1\rrbracket$ impliquent que $h \prec h_{k,\mathbf{j}}$. Il s'ensuit que $h=2^{\ell_n-1}$ et les conditions (\ref{cc}) entraînent que $P$ n'appartient pas aux $D_h$ tels que $h \not \in \{ 2^{\ell_n-1},h_{2,\mathbf{j}},\dots,h_{n-1,\mathbf{j}} \}$. Les conditions (\ref{ccc}) ne sont pas vérifiées par un élément $h \in  \{ 2^{\ell_n-1},h_2,\dots,h_{n-1} \}$. En effet, un tel élément vérifie $\varepsilon_{\ell_n}(h)=1$ et ne vérifie jamais $h \preceq h'_{k,\boldsymbol{\ell}}$ pour tout $k \in \llbracket 3,n-1 \rrbracket$. Ainsi les conditions (\ref{ccc}) impliquent que le point $P$ n'appartient pas à $D_h$ pour $h \in \llbracket 1,N-1\rrbracket \smallsetminus \{ 2^{\ell_n-1},h_2,\dots,h_{n-1} \}$.\\
\indent
Le cas où $x_{\ell} \neq 0$ se traite de façon parfaitement analogue. Le sous-ensemble exceptionnel défini dans \cite{Co} et dans la discussion précédant la proposition est donc donné par les monômes 
$$
z_{h_n}\prod_{1 \leqslant h \leqslant N \atop h \not \in \mathcal{H}_{\mathbf{j}}}z_h\neq 0 \quad \mbox{ou} \quad x_i\prod_{1 \leqslant h \leqslant N \atop h \not \in \mathcal{H}_{\mathbf{j}}}z_h\neq 0
\label{condd1}
$$
pour $i \in \llbracket 1,n \rrbracket$ et $\mathbf{j}=(j_1,\dots,j_n)$ tel que $\{j_1,\dots,j_n\}=\{1,\dots,n\}$ et 
$$
\mathcal{H}_{\mathbf{j}}=\left\{ h_{1,{\mathbf{j}}}=2^{j_1-1}\preceq \cdots \preceq h_{k,{\mathbf{j}}}=\sum_{s=1}^k 2^{j_s-1} \preceq \cdots \preceq h_{n,{\mathbf{j}}}=N \right\}.
$$
Par conséquent, $\mathcal{T}=\mathbb{A}^{2^n+n-1} \smallsetminus F$ est bien l'ouvert $\Omega$ défini dans l'énoncé de la Proposition \ref{torseurversel}.\\
Grâce à la description du morphisme $\varphi: \mathcal{T} \rightarrow \Xi_0$ rappelée dans \cite{Blomer2014} dans leur preuve du lemme 10, on sait que la restriction de $\varphi$ à $\mathbb{G}_m^{2^n+n-1}$ est le morphisme de tores $\mathbb{G}_m^{2^n+n-1} \rightarrow U$ dual du morphisme "diviseur" $\mathbb{Q}[U]^{\ast}/\mathbb{Q}^{\ast}\rightarrow \mbox{Div}_{\Xi_0 \smallsetminus U}
(\Xi_0)$, où $\mbox{Div}_{\Xi_0 \smallsetminus U}(\Xi_0)$ désigne le groupe abélien libre des diviseurs de $\Xi_0$ à support dans $\Xi_0 \smallsetminus U$. Le groupe $\mathbb{Q}[U]^{\ast}/\mathbb{Q}^{\ast}$ est engendré par $x_1/y_n, \dots,x_n/y_n,y_1/y_n,\dots,$ $y_{n-1}/y_n$. Puisqu'on a 
$$
\mbox{div}(x_i)=X_i, \quad \mbox{div}(y_i)=\sum_{1 \leqslant h \leqslant N \atop \varepsilon_i(h)=1} D_h \quad \mbox{pour} \quad i \in \llbracket 1,n \rrbracket
$$
ainsi que pour tout $h \in \llbracket 1,N-1 \rrbracket$,
$$
\mbox{div}\left(Y^{(h)}_i\right)=\sum_{1 \leqslant\ell \leqslant N \atop \varepsilon_i(\ell)=1} D_{\ell}-\sum_{1 \leqslant\ell \leqslant N \atop h \preceq \ell} D_{\ell} \quad \mbox{pour} \quad \varepsilon_i(h)=1
$$
et
$$
\mbox{div}\left(Z^{(h)}_i\right)=\sum_{1 \leqslant\ell \leqslant N \atop \varepsilon_i(\ell)=0} D_{\ell}-\sum_{1 \leqslant\ell \leqslant N \atop \varepsilon_j(\ell)=0 \hspace{1mm} {\rm si} \hspace{1mm} \varepsilon_j(h)=1} D_{\ell} \quad \mbox{pour} \quad \varepsilon_i(h)=1,
$$
on obtient aisément que la restriction de $\varphi$ à l'ouvert dense de $\Omega$ sur lequel chacune des coordonnées est non nulle est bien donnée par (\ref{app}) et on conclut alors par densité la démonstration de cette proposition.
\hfill
$\square$\\
\newline
\indent
On déduit de la Proposition \ref{torseurversel} le théorème suivant de manière analogue à \cite{Blomer2014}.

\begin{theor}
Soit $O \subset \mathbb{A}^{2^n+n-1}$ la sous-variété ouverte  de coordonnées $\left(\mathbf{x}';\left(z_h\right)_{1 \leqslant h \leqslant N}\right)$ définie par les conditions
\begin{eqnarray}
z_{h_n}\prod_{1 \leqslant h \leqslant N \atop h \not \in \mathcal{H}_{\mathbf{j}}}z_h\neq 0 \quad \mbox{ou} \quad x'_i z_{2^{i-1}}\prod_{1 \leqslant h \leqslant N \atop h \not \in \mathcal{H}_{\mathbf{j}}}z_h\neq 0
\label{condd2}
\end{eqnarray}
pour $i \in \llbracket 1,n \rrbracket$ et $\mathbf{j}=(j_1,\dots,j_n) $ tel que $\{j_1,\dots,j_n\}=\{1,\dots,n\}$ et 
$$
\mathcal{H}_{\mathbf{j}}=\left\{ h_{1,{\mathbf{j}}}=2^{j_1-1}\preceq \cdots \preceq h_{k,{\mathbf{j}}}=\sum_{s=1}^k 2^{j_s-1} \preceq \cdots \preceq h_{n,{\mathbf{j}}}=N \right\}.
$$
Le morphisme $\varphi_O: O \rightarrow \Xi_0$ défini par
\begin{eqnarray}
\left(\mathbf{x}';\left(z_h\right)_{1 \leqslant h \leqslant N}\right) \longmapsto \left( z_1x'_1,\dots,z_{2^{n-1}}x'_n,\prod_{1 \leqslant h \leqslant N \atop \varepsilon_1(h)=1}z_h,\dots, \prod_{1 \leqslant h \leqslant N \atop \varepsilon_n(h)=1}z_h;\prod_{1 \leqslant h \leqslant N} \left(\mathbf{Y}^{(h)} ; \mathbf{Z}^{(h)}\right)\right)
\label{app2}
\end{eqnarray}
où
$$
\mathbf{Y}^{(h)}=\left( \frac{\displaystyle\underset{\scriptscriptstyle 1 \leqslant \ell \leqslant N \atop\scriptscriptstyle \varepsilon_{i_1}(\ell)=1}{\prod}z_{\ell}}{\displaystyle\underset{\scriptscriptstyle 1 \leqslant \ell \leqslant N \atop\scriptscriptstyle h \preceq \ell}{\prod}z_{\ell}},\dots,\frac{\displaystyle\underset{\scriptscriptstyle 1 \leqslant \ell \leqslant N \atop\scriptscriptstyle \varepsilon_{i_{s(h)}}(\ell)=1}{\prod}z_{\ell}}{\displaystyle\underset{\scriptscriptstyle 1 \leqslant \ell \leqslant N \atop\scriptscriptstyle h \preceq \ell}{\prod}z_{\ell}} \right) \quad \mbox{pour} \quad \varepsilon_{i_1}(h)=\cdots=\varepsilon_{i_{s(h)}}(h)=1
$$
et
$$
\mathbf{Z}^{(h)}= \left(\frac{\displaystyle\underset{\scriptscriptstyle 1 \leqslant \ell \leqslant N \atop\scriptscriptstyle \varepsilon_{i_1}(\ell)=0}{\prod}z_{\ell}}{\displaystyle\underset{\scriptscriptstyle 1 \leqslant \ell \leqslant N \atop\scriptscriptstyle \varepsilon_{i}(\ell)=0 \hspace{1mm} {\rm si} \hspace{1mm} \varepsilon_i(h)=1}{\prod}z_{\ell}},\dots,\frac{\displaystyle\underset{\scriptscriptstyle 1 \leqslant \ell \leqslant N \atop\scriptscriptstyle \varepsilon_{i_{s(h)}}(\ell)=0}{\prod}z_{\ell}}{\displaystyle\underset{\scriptscriptstyle 1 \leqslant \ell \leqslant N \atop\scriptscriptstyle \varepsilon_{i}(\ell)=0 \hspace{1mm} {\rm si} \hspace{1mm} \varepsilon_i(h)=1}{\prod}z_{\ell}} \right) \quad \mbox{pour} \quad \varepsilon_{i_1}(h)=\cdots=\varepsilon_{i_{s(h)}}(h)=1
$$
est alors le morphisme sous-jacent d'un torseur versel pour $X_{0,n}$.
\label{torseurversel2}
\end{theor}
\noindent
\textit{Démonstration.--}
On obtient grâce à la suite spectrale de Leray (voir \cite{Blomer2014}) le diagramme commutatif de $\mathfrak{g}$-modules triviaux suivants
$$
\xymatrix{
0 \ar[r] &\mbox{Pic}(\overline{B}_0)\ar[r]^{{\overline{\gamma}^{\ast}}} \ar[d]^{{{\rm Id}}}&\mbox{Pic}(\overline{\Xi}_0) \ar[r] \ar[d]&  \mathbb{Z} \ar[r] \ar[d]^{{\rm Id}}&0\\
0 \ar[r] & \mbox{Pic}(\overline{B}_0) \ar[r]^{{\overline{\lambda}^{\ast}}}  & \mbox{Pic}(\overline{X}_0) \ar[r]&  \mathbb{Z} \ar[r] & 0 \\
}
$$
où $\mathfrak{g}=\mbox{Gal}\left( \overline{\mathbb{Q}}/\mathbb{Q}\right)$, $\overline{\lambda}^{\ast}: \overline{X}_0 \rightarrow \overline{B}_0$ et $\overline{\gamma}^{\ast}: \overline{\Xi}_0 \rightarrow \overline{B}_0$ sont les morphismes sur $\overline{\mathbb{Q}}$ provenant respectivement des fibrations $\lambda: X_{0,n} \rightarrow B_{0,n}$ et $\gamma: \Xi_0 \rightarrow B_{0,n}$. La flèche $\mbox{Pic}(\overline{\Xi}_0) \rightarrow \mbox{Pic}(\overline{X}_0)$ est un isomorphisme et on a une suite exacte duale de $\mathbb{Q}$-tores algébriques
$$
\xymatrix{
1 \ar[r] &\mathbb{G}_m \ar[r]&T \ar[r]& S \ar[r] &1\\
}
$$
avec $T$ le tore dont le groupe des caractères est $\hat{T}=\mbox{Pic}(\overline{\Xi}_0)=\mbox{Pic}(\overline{X}_0)$ et $S$ le tore dont le groupe des caractères est $\hat{S}=\mbox{Pic}(\overline{B}_0)$. De la fonctorialité de la suite exacte (voir \cite{CS})
$$
\xymatrix{
0 \ar[r] &H^1_{{\footnotesize\mbox{\'et}}}(\mathbb{Q},T) \ar[r]&H^1_{{\footnotesize\mbox{\'et}}}(X_{0,n},T) \ar[r]^{\hspace{-7mm}\chi}& \mbox{Hom}_{\mathfrak{g}}\left(\hat{T}, \mbox{Pic}\left( \overline{X}_0\right)\right)   \\
}
$$
par rapport à $X_{0,n} \rightarrow \Xi_0$, il s'ensuit que le $T$-torseur versel pour $\Xi_0$ se restreint au sous-ensemble $\varphi^{-1}(X_{0,n}) \subseteq \Omega$ à un $T$-torseur versel pour $X_{0,n}$. Le sous-ensemble $\varphi^{-1}(X_{0,n})$ est défini par l'équation (\ref{4}). Après avoir remarqué que (\ref{cccond}) et (\ref{4}) sont équivalentes lorsque (\ref{condd1}) est vérifiée 
et en utilisant la Proposition \ref{torseurversel}, cette dernière prend la forme
\begin{eqnarray}
\sum_{i=1}^n \xi_i \left(\prod_{1 \leqslant h \leqslant N} z_h^{1-\varepsilon_i(h)} \right)=0.
\label{cccond}
\end{eqnarray}
\`A la manière de \cite{Blomer2014}, on définit les $n$ fonctions régulières $x'_1,\dots,x'_n$ sur $\varphi^{-1}(X_{0,n})$. Sur l'ouvert
$$
\prod_{1 \leqslant h \leqslant N} z_h^{1-\varepsilon_i(h)} \neq 0
$$
pour $i \in \llbracket 1,n \rrbracket$ fixé, on pose
\begin{eqnarray}
x'_i=-\frac{\displaystyle\overset{n}{\underset{\scriptscriptstyle j=1 \atop\scriptscriptstyle j \neq i}{\sum}} \xi_j \left(\underset{\scriptscriptstyle 1 \leqslant h \leqslant N\atop\scriptscriptstyle h \neq z_{2^{i-1}}}{\prod} z_h^{1-\varepsilon_j(h)} \right)}{\displaystyle\underset{\scriptscriptstyle 1 \leqslant h \leqslant N}{\prod} z_h^{1-\varepsilon_i(h)}} \quad \mbox{et} \quad x'_j=\frac{\xi_j}{z_{2^{j-1}}} \quad \mbox{pour} \quad j \in \llbracket 1,n\rrbracket \smallsetminus\{i\}.
\label{condd3}
\end{eqnarray}
Si l'on impose de plus $z_{2^{i-1}}\neq 0$, alors on vérifie aisément que $x'_i=\frac{\xi_i}{z_{2^{i-1}}}$, ce qui assure que les définitions sur les différents ouverts lorsque $i$ varie dans $\llbracket 1,n \rrbracket$ sont compatibles et se recollent pour donner des fonctions régulières bien définies sur $\varphi^{-1}(X_{0,n})$. Cela permet d'obtenir la proposition.
\hfill
$\square$\\
\newline
\indent
On adapte alors la fin de la section 4 de \cite{Blomer2014} pour obtenir des informations sur les points entiers qui seront nécessaires afin de calculer le nombre de Tamagawa apparaissant dans la constante de Peyre. Soit $\underline{X_{0,n}} \subseteq \mathbb{P}_{\mathbb{Z}}^{2n-1} \times \prod_{1 \leqslant h \leqslant N} \mathbb{P}_{\mathbb{Z}}^{(h)}\times \mathbb{P}_{\mathbb{Z}}^{(h)}$ défini par les équations (\ref{4}), (\ref{5}), (\ref{6}), (\ref{7}) et (\ref{8}) et $\underline{\Xi_0} \subseteq \mathbb{P}_{\mathbb{Z}}^{2n-1} \times \prod_{1 \leqslant h \leqslant N} \mathbb{P}_{\mathbb{Z}}^{(h)}\times \mathbb{P}_{\mathbb{Z}}^{(h)}$ défini par les équations (\ref{4}), (\ref{5}), (\ref{6}) et (\ref{8}). D'après \cite[pages 22-23]{Fu}, on peut étendre le morphisme $\varphi: \Omega \rightarrow \Xi_0$ de la Proposition \ref{torseurversel} en un morphisme $\underline{\varphi}:\underline{\Omega} \rightarrow \underline{\Xi_0}$ entre deux schémas toriques puisque le morphisme de Cox provient d'un morphisme d'éventails. Les arguments utilisés lors de la démonstration de la Proposition \ref{torseurversel} mais sur $\mathbb{Z}$ fournissent que $\underline{\Omega}$ est le sous-schéma de $\mathbb{A}_{\mathbb{Z}}^{2^n+n-1}$ de coordonnées $\left(\mathbf{x},\left( z_h\right)_{1 \leqslant h \leqslant N}\right)$ défini par les conditions (\ref{condd1}) et le morphisme $\underline{\varphi}:\underline{\Omega} \rightarrow \underline{\Xi_0}$ est défini par (\ref{app}) est le morphisme sous-jacent d'un torseur $\underline{\varphi}_{\underline{\mathcal{T}}}:\underline{\mathcal{T}} \rightarrow \underline{\Xi_0}$ sous un tore déployé sur $\mathbb{Z}$ que l'on notera $\underline{T} \cong \mathbb{G}_{m,\mathbb{Z}}^{2^n-n}$ avec $H^1_{{\footnotesize\mbox{\'et}}}(\mathbb{Z},\underline{T})$. Ainsi il y a une bijection entre les $\underline{T}(\mathbb{Z})$-orbites de points entiers de $\underline{\Omega}$ et les points entiers de $\underline{\Xi_0}$ (voir \cite{Blomer2014} et \cite[III.4.9]{Mi} pour plus de détails). De même, en restreignant $\underline{\varphi}$ au fermé $\underline{O}=\underline{\varphi}^{-1}(\underline{X_{0,n}})$ de $\underline{\Omega}$ défini par l'équation (\ref{7}), on peut également introduire des coordonnées $\left(\mathbf{x}';\left(z_h\right)_{1 \leqslant h \leqslant N}\right)$ telles que $\underline{O}$ soit le sous-schéma localement fermé   de $\mathbb{A}_{\mathbb{Z}}^{2^n+n-1}$ défini par (\ref{cccond}) et (\ref{condd2}). On obtient également un morphisme $\underline{\varphi}_{\underline{O}}:\underline{O} \rightarrow \underline{X_{0,n}}$ donné par (\ref{app2}). On dispose alors du lemme suivant qui permet d'obtenir la Proposition \ref{propeq} \textit{infra}, équivalente au Lemme \ref{bij}. Mais avant de pouvoir établir cette Proposition \ref{propeq}, on a besoin du lemme auxiliaire suivant.

\begin{lemme}
La condition de coprimalité
\begin{eqnarray} 
\underset{\mathbf{j}}{\rm{pgcd}}\left(\prod_{1 \leqslant h \leqslant N \atop h \not \in \mathcal{H}_{\mathbf{j}}}z_h \right)=1
\label{condcoprim}
\end{eqnarray}
où $\mathbf{j}$ décrit tous les $n$-uplets $\mathbf{j}=(j_1,\dots,j_n)$ tel que $\{j_1,\dots,j_n\}=\{1,\dots,n\}$ et 
$$
\mathcal{H}_{\mathbf{j}}=\left\{ h_{1,{\mathbf{j}}}=2^{j_1-1}\preceq \cdots \preceq h_{k,{\mathbf{j}}}=\sum_{s=1}^k 2^{j_s-1} \preceq \cdots \preceq h_{n,{\mathbf{j}}}=N \right\}
$$
est équivalente au fait que le $N$-uplet $(z_h)_{1 \leqslant N \leqslant N}$ soit réduit.
\label{lemmecoprime}
\end{lemme}
\noindent
\textit{Démonstration.--}
Pour simplifier les notations, on notera dans toute cette démonstration
$$
P_{\mathbf{j}}=\prod_{1 \leqslant h \leqslant N \atop h \not \in \mathcal{H}_{\mathbf{j}}}z_h
$$
pour $\mathbf{j}=(j_1,\dots,j_n)$ tel que $\{j_1,\dots,j_n\}=\{1,\dots,n\}$ et 
\begin{eqnarray}
\mathcal{H}_{\mathbf{j}}=\left\{ h_{1,{\mathbf{j}}}=2^{j_1-1}\preceq \cdots \preceq h_{k,{\mathbf{j}}}=\sum_{s=1}^k 2^{j_s-1} \preceq \cdots \preceq h_{n,{\mathbf{j}}}=N \right\}.
\label{hj}
\end{eqnarray}
Supposons dans un premier temps que la condition de coprimalité suivante
$$
{\rm{pgcd}}\left(P_{\mathbf{j}} \hspace{1mm} : \hspace{1mm} \mathbf{j}\right)=1 
$$
pour $\mathbf{j}=(j_1,\dots,j_n)$ tel que $\{j_1,\dots,j_n\}=\{1,\dots,n\}$ et 
$
\mathcal{H}_{\mathbf{j}}$ comme en (\ref{hj})) est vérifiée et considérons $h$ et $\ell$ deux éléments de $\llbracket 1,N-1 \rrbracket$ non comparables, c'est-à-dire tels que $h \not \preceq \ell$ et $\ell \not \preceq h$. Soit $\mathbf{j}=(j_1,\dots,j_n)$ tel que $\{j_1,\dots,j_n\}=\{1,\dots,n\}$. On constate alors que soit
$
z_h \mid P_{\mathbf{j}}
$
soit $z_{\ell} \mid P_{\mathbf{j}}$. En effet, supposons que $z_h \nmid P_{\mathbf{j}}$. Alors 
$
h \in \mathcal{H}_{\mathbf{j}}$. Clairement tous les éléments de $ \mathcal{H}_{\mathbf{j}}$ sont alors comparables à $h$ si bien que 
$
\ell \not \in \mathcal{H}_{\mathbf{j}}$ et ainsi $z_{\ell} \mid P_{\mathbf{j}}$ d'après la définition de $P_{\mathbf{j}}$. On en déduit par conséquent que la condition de coprimalité (\ref{condcoprim}) implique que le $N$-uplet $(z_h)_{1 \leqslant N \leqslant N}$ soit réduit.\\
\par
Réciproquement, supposons que le $N$-uplet $(z_h)_{1 \leqslant N \leqslant N}$ soit réduit. Raisonnons alors par l'absurde en supposant qu'il existe un nombre premier $p$ tel que 
$
p \mid P_{\mathbf{j}}
$
pour tout $\mathbf{j}=(j_1,\dots,j_n)$ tel que $\{j_1,\dots,j_n\}=\{1,\dots,n\}$ et 
$
\mathcal{H}_{\mathbf{j}}=\left\{ h_{1,{\mathbf{j}}}=2^{j_1-1}\preceq \cdots \preceq h_{k,{\mathbf{j}}}=\sum_{s=1}^k 2^{j_s-1} \preceq \cdots \preceq h_{n,{\mathbf{j}}}=N \right\}$ comme en (\ref{hj}). En particulier, $p \mid  P_{\{ 1,3,\dots,n-1,2 \}}$ et par conséquent il existe $H_1 \in \llbracket 1,N-1 \rrbracket$ tel que
$$
\exists i \in \llbracket 1,n \rrbracket, \quad z_{H_1} \prec h_{i,\{ 1,3,\dots,n-1,2 \}}.
$$ 
Posons $k=s(H_1)$ et $c_1<c_2<\cdots<c_k$ les chiffres de $H_1$. Il existe alors $(n-k)!k!$ facteurs $P_{\mathbf{j}}$ qui ne sont pas divisibles par $z_{H_1}$. En effet, il s'agit des $\mathbf{j}$ tels que $\{j_1,\dots,j_k\}=\{c_1,\dots,c_k\}$. Soit alors un tel $\mathbf{j}$. Puisque $p \mid P_{\mathbf{j}}$, il existe $H_2 \in \llbracket 1,N-1 \rrbracket\smallsetminus\{H_1\}$ tel que $p \mid z_{H_2}$ et
$$
\exists i \in \llbracket 1,n \rrbracket, \quad z_{H_2} \prec h_{i,\mathbf{j}}.
$$
On a alors deux cas. Puisque $p \mid (z_{H_1},z_{H_2})$ et que le $N$-uplet $(z_h)_{1 \leqslant N \leqslant N}$ est réduit, on en déduit que $H_1$ et $H_2$ sont comparables. Supposons alors par exemple que $H_2\prec H_1$ et posons $\ell=s(H_2)<k=s(H_1)$. Le nombre de termes $P_{\mathbf{j}}$ qui ne sont ni divisibles par $z_{H_1}$ ni divisibles par $z_{H_2}$ est de $(n-k)!(k-\ell)!\ell!$. En effet, il s'agit des $\mathbf{j}$ tels que $j_1,\dots,j_{\ell}$ soient les chiffres de $H_2$ et $j_{\ell+1},\dots,j_k$ soient les chiffres de $H_1$ qui ne sont pas des chiffres de $H_2$. Le cas $H_1\prec H_2$ se traite de façon analogue. En itérant ce procédé, on construit $H_1,H_2,\dots,H_{n-1} \in \llbracket 1,N \rrbracket^{n-1}$ tels que
$$
\exists \sigma \in \mathfrak{S}_{n-1}, \quad H_{\sigma(1)}\prec H_{\sigma(2)} \prec \cdots \prec H_{\sigma(n-1)}
$$ 
et tel que le seul facteur $P_{{\mathbf{j}_0}}$ qui ne soit divisible par aucun des $z_{H_{i}}$ pour $i \in \llbracket 1,n-1\rrbracket$ soit le $n$-uplet $\mathbf{j}_0$ défini par $j_1$ l'unique chiffre de $H_{\sigma(1)}$ puis pour tout $i \in \llbracket 2,n-1\rrbracket$, $j_i$ l'unique chiffre de $H_{\sigma(i)}$ qui n'est pas un chiffre de $H_{\sigma(i-1)}$. On pose alors $j_n$ tel que $\{j_1,\dots,j_n\}=\{1,\dots,n\}$. Puisque $p \mid P_{{\mathbf{j}_0}}$, il existe un facteur $z_h$ de $P_{\mathbf{j}_0}$ tel que $p \mid z_h$. Montrons alors que pour tout facteur $z_h$ de $P_{{\mathbf{j}_0}}$, il existe un $i \in \llbracket 1,n-1\rrbracket$ tel que $h$ et $H_i$ ne soient pas comparables. Soit $h$ tel que $z_h \mid P_{{\mathbf{j}_0}}$. On a alors que
$$
z_h \not \in \mathcal{H}_{{\mathbf{j}_0}}=\left\{ H_{\sigma(1)}\prec \cdots \prec H_{\sigma(k)} \prec \cdots \prec H_{\sigma(n-1)} \prec N \right\}.
$$
Si $h \not \prec H_{\sigma(n-1)}$, alors $z_h$ et $z_{H_{\sigma(n-1)}}$ ne sont pas comparables. On peut donc supposer que $h \prec H_{\sigma(n-1)}$. On n'a alors que $h \not \succ H_{\sigma(n-2)}$ car il n'existe aucun $r \in \llbracket 1,N-1\rrbracket$ tel que $z_r \succ H_{\sigma(n-2)}$ et $z_r \prec H_{\sigma(n-1)}$. Si maintenant $h$ et $H_{\sigma(n-2)}$ ne sont pas comparable, alors on a le résultat. On peut donc désormais supposer que $h \prec H_{\sigma(n-2)}$. Alors de même soit $h$ et $H_{\sigma(n-3)}$ ne sont pas comparables soit $h \prec H_{\sigma(n-3)}$ et de proche en proche, ce raisonnement permet de conclure que dans tous les cas, il existe un $i \in \llbracket 1,n-1\rrbracket$ tel que $h$ et $H_i$ ne soient pas comparables. On obtient donc finalement un 
$$
h \not \in \mathcal{H}_{{\mathbf{j}_0}}=\left\{ H_{\sigma(1)}\prec \cdots \prec H_{\sigma(k)} \prec \cdots \prec H_{\sigma(n-1)} \prec N \right\}
$$
et un $i \in \llbracket 1,n-1\rrbracket$ tel que $p \mid (z_h,z_{H_i})$. Mais puisque $h$ et $H_i$ ne sont pas comparables, on a une contradiction et finalement, on a bien la condition de coprimalité (\ref{condcoprim}). Cela conclut la preuve de ce lemme.
\hfill
$\square$

\begin{prop}
Soit $\mathcal{A}_0$ l'ensemble des $2^n+n-1$-uplets $\left(\mathbf{x}';\left(z_h\right)_{1 \leqslant h \leqslant N}\right)$ tels que $x'_i \in \mathbb{Z}$, $z_h \in \mathbb{N}$ si $s(h)>1$ et $z_h \in \mathbb{Z}\smallsetminus\{0\}$ si $s(h)=1$ vérifiant la relation 
\begin{eqnarray}
\sum_{i=1}^n x'_i \left(\prod_{\substack{\scriptscriptstyle 1 \leqslant h \leqslant N\\\scriptscriptstyle h \neq h_n \\\scriptscriptstyle s(h)\geqslant 2 }} z_h^{1-\varepsilon_i(h)} \right)=0.
\label{cccond2}
\end{eqnarray}
ainsi que les conditions de coprimalité
\begin{eqnarray}
{\rm{pgcd}}\left( z_{h_n},x'_1z_1,\dots,x'_nz_{2^{n-1}}\right)=1 \quad 
\label{condcoprim2}
\end{eqnarray}
et que le $N$-uplet $\left(z_h\right)_{1 \leqslant h \leqslant N}$ est réduit. Alors, l'application $\mathcal{A}_0 \rightarrow \mathbb{Z}^{2n}$ définie par (\ref{bij2}) est une bijection entre $\mathcal{A}_0$ et $\mathcal{W}_n$ défini en (\ref{wn}), l'ensemble des solutions entières primitives de (\ref{eq}). De plus, dans ce cas, la condition (\ref{condcoprim2}) est équivalente à la condition $\mbox{pgcd}(x_1,\dots,x_n,y_1,\dots,y_n)=1$.
\label{propeq}
\end{prop}
\noindent
\textit{Démonstration.--}
La démonstration est très fortement inspirée de celle de \cite[lemma 11]{Blomer2014}. Soit $\underline{W_n}\subseteq \mathbb{P}^{2n-1}_{\mathbb{Z}}$ le sous-schéma défini par l'équation (\ref{eq}) et $\underline{f}: \underline{X_{0,n}} \rightarrow \underline{W_n}$ le morphisme induit par la projection
$$
\mathbb{P}_{\mathbb{Z}}^{2n-1} \times \prod_{1 \leqslant h \leqslant N} \mathbb{P}_{\mathbb{Z}}^{(h)}\times \mathbb{P}_{\mathbb{Z}}^{(h)} \longrightarrow \mathbb{P}_{\mathbb{Z}}^{2n-1},
$$
provenant de $f:X_{0,n} \rightarrow W_n$ après changement de base. On a alors des bijections naturelles $ \underline{X_{0,n}}(\mathbb{Z})=X_{0,n}(\mathbb{Q})$, $\underline{W_n}(\mathbb{Z})=W_n(\mathbb{Q})$ et $X_{0,n}^{\circ}(\mathbb{Q})=W_n^{\circ}(\mathbb{Q})$ pour $X_{0,n}^{\circ} \subseteq X_{0,n}$ et $W_n^{\circ} \subseteq W_n$ les deux ouverts définis par les conditions $y_1y_2\cdots y_n \neq 0$. Si l'on considère à présent les ouverts $\underline{X_{0,n}}^{\circ} \subseteq X_{0,n}$ et $\underline{W_n}^{\circ} \subseteq W_n$ provenant respectivement de $X_{0,n}^{\circ} \subseteq X_{0,n}$ et $W_n^{\circ} \subseteq W_n$, il vient la bijection $\underline{X_{0,n}}^{\circ}(\mathbb{Z})=\underline{W_n}^{\circ}(\mathbb{Z})$. Enfin, on pose $O^{\circ} \subseteq O$ l'ouvert défini par la condition
$$
\prod_{1 \leqslant h \leqslant N} z_h \neq 0.
$$
On introduit alors $\underline{O}^{\circ}(\mathbb{Z}) \subseteq \underline{O}(\mathbb{Z})$ correspondant à $O^{\circ}(\mathbb{Q}) \subseteq O(\mathbb{Q})$ et tel qu'on ait une bijection $\underline{O}^{\circ}(\mathbb{Z})=O^{\circ}(\mathbb{Q})$. Comme $\varphi^{-1}\left(W_n^{\circ}\right)=O^{\circ}$, il s'ensuit une bijection entre $\underline{X_{0,n}}^{\circ}$ et les $\underline{T}(\mathbb{Z})$-orbites de $\underline{O}^{\circ}(\mathbb{Z})$. L'application de $\underline{O}^{\circ}(\mathbb{Z})$ vers $\underline{W_n}^{\circ}(\mathbb{Z})$ s'obtient grâce à $\underline{f} \circ \underline{\varphi}$ donnée par (\ref{app2})\\
\par
Un point
$
\left(\mathbf{x}';\left(z_h\right)_{1 \leqslant h \leqslant N}\right) \in \underline{O}^{\circ}(\mathbb{Z})
$ si, et seulement si, il vérifie (\ref{cccond2}) et (\ref{cccond}) modulo $p$ pour tout nombre premier $p$, ce qui est bien équivalent à (\ref{cccond2}), (\ref{condcoprim}) et (\ref{condcoprim2}). On obtient alors les conditions de coprimalité de l'énoncé de la proposition grâce au Lemme \ref{lemmecoprime}. On identifie alors à présent l'ensemble $\underline{X_{0,n}}^{\circ}(\mathbb{Z})=\underline{W_n}^{\circ}(\mathbb{Z})$ avec l'ensemble des $(2n-1)$-uplets de la forme $\pm (x_1,\dots,x_n,y_1,\dots,y_n)$ à coordonnées premières entre elles et vérifiant (\ref{eq}) et $y_1\cdots y_n \neq 0$. Chaque $\underline{T}(\mathbb{Z})$-orbites de $\underline{O}^{\circ}(\mathbb{Z})$ comporte $2^{\dim(T)}=2^{2^n-n}$ éléments qui ne différent que par le signe de certaines composantes. Soit $\left(\mathbf{x}';\left(z_h\right)_{1 \leqslant h \leqslant N}\right) \in \underline{O}^{\circ}(\mathbb{Z})$ vérifiant (\ref{condcoprim}) et (\ref{condcoprim2}). D'après le Lemme \ref{lemmecoprime}, cela équivaut au fait que $\left(\mathbf{x}';\left(z_h\right)_{1 \leqslant h \leqslant N}\right)$ vérifie les conditions (\ref{condcoprim}) et que $\left(z_h\right)_{1 \leqslant h \leqslant N}$ est réduit. On a alors 
$$
\underline{f} \circ \underline{\varphi}=(x_1,\dots,x_n,y_1,\dots,y_n) \in \underline{W_n}^{\circ}(\mathbb{Z})
$$
avec $\mbox{pgcd}(x_1,\dots,x_n,y_1,\dots,y_n)=1$. D'après les propriétés générales des torseurs (voir \cite[Chapter 2]{Sk}), on sait que la $\underline{T}(\mathbb{Z})$-orbites de $\left(\mathbf{x}';\left(z_h\right)_{1 \leqslant h \leqslant N}\right)$ est la fibre du point $(x_1,\dots,x_n,y_1,\dots,y_n)$. D'après le Lemme \ref{lemmefacto}, les points de cette fibre ne diffèrent que par leurs signes. \'Etudions alors combien de points de $\underline{O}^{\circ}(\mathbb{Z})$ appartiennent à cette fibre. On commence par fixer $x'_1, \dots, x'_n$, ce qui fixe $z_1,z_2,\dots,z_{2^{n-1}}$ et donne $2^n$ choix. Il s'agit alors de déterminer les signes possibles de $\left(z_h\right)_{1 \leqslant h \leqslant N \atop s(h) \geqslant 2}$ tels que
$$
\forall i \in \llbracket 1,n \rrbracket, \quad \prod_{1 \leqslant h \leqslant N \atop \varepsilon_i(h)=1, \hspace{1mm} h \neq 2^{i-1}}z_h
$$
soit de signe prescrit. Autrement dit, on cherche le cardinal des $(i_h)_{1 \leqslant h \leqslant N \atop s(h) \geqslant 2} \in  \left(\mathbb{Z}/2\mathbb{Z}\right)^{2^n-n-1}$ tels que
$$
\forall i \in \llbracket 1,n \rrbracket, \quad \sum_{1 \leqslant h \leqslant N \atop \varepsilon_i(h)=1, \hspace{1mm} h \neq 2^{i-1}}i_h \equiv k_i\Mod{2}
$$
pour $(k_1,\dots,k_n) \in \left(\mathbb{Z}/2\mathbb{Z}\right)^n$ fixé. On peut choisir librement dans la première équation tous les $i_h$ tels que $s(h) \geqslant 2$ et $\varepsilon_1(h)=1$ excepté $i_{N-1}$ qui est alors fixé par l'équation, ce qui laisse
$$
\#\{h \in \llbracket 1,N\rrbracket \hspace{1mm}: \hspace{1mm} \varepsilon_1(h)=1\}-2=2^{n-1}-2
$$
choix. De même, dans la deuxième équation, on peut choisir librement tous les $i_h$ tels que $s(h) \geqslant 2$, $\varepsilon_1(h)=0$ et $\varepsilon_2(h)=1$ excepté $i_{N-2}$ qui est alors fixé par l'équation, ce qui laisse
$$
\#\{h \in \llbracket 1,N\rrbracket \hspace{1mm}: \hspace{1mm} \varepsilon_2(h)=1, \hspace{1mm} \varepsilon_1(h)=0\}-2=2^{n-2}-2
$$
choix. De même, on peut choisir comme cela librement dans l'équation $k \in \llbracket 1,n-2\rrbracket$, $2^{n-k}-2$ variables $i_h$ librement tels que $s(h) \geqslant 2$ et $\varepsilon_k(h)=1$, $\varepsilon_1(h)=\cdots=\varepsilon_{k-1}(h)=0$. On constate alors que toutes les variables des deux dernières équations sont fixée excepté $i_{2^{n-2}+2^{n-1}}$. Il y a alors deux cas de figure, soit cette variable ne peut pas être fixée de sorte que les équations $n-1$ et $n$ soient vérifiées, soit il existe une unique valeur de $i_{2^{n-2}+2^{n-1}}$ telle que les deux dernières équations $n-1$ et $n$ soient vérifiées. De plus, cette variable peut être fixée si, et seulement si,
$$
\sum_{\substack{\varepsilon_{n-1}(h)=1, \hspace{1mm} \varepsilon_n(h)=0\\ h \neq 2^{n-2}}} i_h \equiv \sum_{\substack{\varepsilon_{n-1}(h)=0, \hspace{1mm} \varepsilon_n(h)=1\\ h \neq 2^{n-2}}} i_h\Mod{2}
$$
où les variables en question dans les deux sommes ci-dessus ne font jamais intervenir les variables $i_h$ pour $h=N-2^{i-1}$ et $i \in \llbracket 1,n-2 \rrbracket$. On obtient ainsi
$$
\left(\sum_{k=1}^{n-2}\left(2^{n-k}-2\right)\right)-1=2^n-2n-1
$$
choix. Finalement, on aboutit à $2^{2^n-n-1}$ éléments au-dessus de $(x_1,\dots,x_n,y_1,\dots,y_n)$. On obtient donc $2^{2^n-n}$ éléments au-dessus du point $\pm (x_1,\dots,x_n,y_1,\dots,y_n)$ de $\underline{W_n}^{\circ}(\mathbb{Z})$.\\
\par
Si on considère à présent, les points d'une $\underline{T}(\mathbb{Z})$-orbites de $\underline{O}^{\circ}(\mathbb{Z})$ tels que $z_h \in \mathbb{N}^{\ast}$ dès que $s(h)>1$, on obtient exactement deux points dont toutes les coordonnées sont égales mis à part les $\left(z_{2^{i-1}}\right)_{1 \leqslant i \leqslant n}$ qui sont opposés. De plus, ces deux points ont pour images par $\underline{f} \circ \underline{\varphi}$ exactement $\pm (x_1,\dots,x_n,y_1,\dots,y_n)$. Cela permet bien d'obtenir le résultat de la proposition.
\hfill
$\square$\\
\subsubsection{Le nombre de Tamagawa $\omega_H\left( X_{0,n}(\mathbb{A}_{\mathbb{Q}})^{{\rm Br}(X_{0,n})} \right)$}
Cette construction, de manière complètement analogue à la section 5 de  \cite{Blomer2014}, permet alors de construire explicitement le nombre de Tamagawa 
$$
\omega_H\left( X_{0,n}(\mathbb{A}_{\mathbb{Q}})^{{\rm Br}(X_{0,n})} \right)=\mu_{\infty}(X_{0,n}(\mathbb{R})) \prod_p \left(1-\frac{1}{p}\right)^{2^n-n}\mu_p(X_{0,n}(\mathbb{Q}_p))
$$
et ainsi d'obtenir l'expression conjecturale complète de la constante de Peyre. On obtient notamment les deux lemmes suivants dont on ne détaillera pas les démonstrations ici puisqu'elles découlent presqu'à la lettre de celles de la section 5 de \cite{Blomer2014}.

\begin{lemme}
On a
$$
\mu_{\infty}(X_{0,n}(\mathbb{R}))\!=\!\!\!\int_{N_v} \!\!\frac{\mbox{d}t_1\cdots \mbox{d}t_{n-1}\mbox{d}t_{n+1} \cdots \mbox{d}t_{2n-1}}{\left| t_{n+1}\cdots t_{2n-1}\right|_v \!\max\!\left( |t_1|_v^n, \dots,\!|t_{n-1}|_v^n,\left| \frac{t_1}{t_{n+1}}+\cdots+\frac{t_{n-1}}{t_{2n-1}} \right|_v^n\!\!,|t_{n+1}|_v^n,\dots,\!|t_{2n-1}|_v^n \right)}.
$$
\label{archimedien}
\end{lemme}

\begin{lemme}
Avec $\underline{C}_n$ la variété torique de Coxeter de $\mathfrak{S}_n$, on a
$$
\left(1-\frac{1}{p}\right)^{2^n-n}\mu_p(X_{0,n}(\mathbb{Q}_p))=\left(1-\frac{1}{p}\right)^{2^n-n-1}\left(1-\frac{1}{p^n}\right)\frac{\left| \underline{C_n}(\mathbb{F}_p) \right|}{p^{n-1}}.
$$
\label{lemme19}
\end{lemme}
\noindent
\textit{Démonstration.--}
La démonstration est inspirée de la section 5 de \cite{Blomer2014} et de l'Annexe. On déduit de \cite[Lemma 19]{Blomer2014} l'égalité
$$
\left(1-\frac{1}{p}\right)^{2^n-n}\mu_p(X_{0,n}(\mathbb{Q}_p))=\frac{\left| \underline{O}(\mathbb{F}_p) \right|}{p^{\dim(O)}}.
$$
Pour déterminer $\left| \underline{O}(\mathbb{F}_p) \right|$, on utilise le fait que le $\underline{X_{0,n}}_{\mathbb{F}_p}$-torseur $\underline{O}_{\mathbb{F}_p}$ sous $\underline{T}_{\mathbb{F}_p}$ est trivial si bien que
$$
\left| \underline{O}(\mathbb{F}_p) \right|=\left| \underline{T}(\mathbb{F}_p) \right| \left| \underline{X_{0,n}}(\mathbb{F}_p) \right|.
$$
Les égalités
$$
\left| \underline{T}(\mathbb{F}_p) \right| =(p-1)^{2^n-n}
$$
et 
$$
\dim(O)=\dim(X_{0,n})+\dim(T)=\dim(X_{0,n})+\mbox{rg}\left( \mbox{Pic}(X_{0,n})\right)
$$
fournissent alors
$$
m_p\left(\underline{O}(\mathbb{Z}_p)\right)=\left(1-\frac{1}{p}\right)^{{\rm{rg}}\left( {\rm{Pic}}(X_{0,n})\right)}\frac{\left| \underline{X_{0,n}}(\mathbb{F}_p) \right|}{p^{\dim(X_{0,n})}}.
$$
On utilise alors pour finir le fait que $X_{0,n}$ soit un $\mathbb{P}^{n-1}$-fibré sur la variété torique $C_n$ pour aboutir à l'expression
$$
\left| \underline{X_{0,n}}(\mathbb{F}_p) \right|=\left| \underline{C_n}(\mathbb{F}_p) \right|\left| \underline{\mathbb{P}}_{\mathbb{F}_p}(\mathbb{F}_p) \right|=\frac{p^n-1}{p-1}\left| \underline{C_n}(\mathbb{F}_p) \right|.
$$
Il s'ensuit
$$
m_p\left(\underline{O}(\mathbb{Z}_p)\right)=\left(p-1\right)^{{\rm{rg}}\left( {\rm{Pic}}(X_{0,n})\right)-1}\left(1-\frac{1}{p^n}\right)\frac{\left| \underline{C_n}(\mathbb{F}_p) \right|}{p^{{\rm{rg}}\left( {\rm{Pic}}(X_{0,n})\right)+\dim(X_{0,n})-n}}.
$$
Comme on a
$$
\dim(X_{0,n})-n=2n-2-n=n-2=\dim(C_n)-1,
$$
on aboutit bien au résultat annoncé
$$
m_p\left(\underline{O}(\mathbb{Z}_p)\right)=\left(1-\frac{1}{p}\right)^{2^n-n-1}\left(1-\frac{1}{p^n}\right)\frac{\left| \underline{C_n}(\mathbb{F}_p) \right|}{p^{n-1}}.
$$
\hfill
$\square$\\

\subsection{Transformation de la constante obtenue par le Théorème \ref{theor1}}
L'objet de cette partie est de réécrire la constante $c_n$ fournie par le Théorème \ref{theor1} afin de vérifier que son expression coïncide avec la forme conjecturée par Peyre et explicitée en section 4.2.

\subsubsection{Mise sous forme de produit eulérien de la quantité $F(\mathbf{1})/\zeta(n)$}
Lorsque le $N-$uplet $\mathbf{z}$ est réduit, la variable $z_N$ n'est sujette à aucune condition de coprimalité tandis que toutes les autres variables sont sujettes à au moins une condition de coprimalité. On note alors $A$ l'ensemble des couples $(k,\ell)$ d'un $N-$uplet $\mathbf{z}$ réduit tels que $\mbox{pgcd}(z_k,z_{\ell})=1$ ainsi que $S=\{1,\dots,N-1\}$ de sorte que $G=(S,A)$ définisse un graphe pour lequel on peut appliquer \cite[theorem 5]{Blomer2014} afin d'écrire $F(\mathbf{1})$ défini en (\ref{F}) sous la forme d'un produit eulérien. Cette étape est capitale dans l'optique de démontrer la conjecture de Peyre puisque cette dernière prédit que le nombre de Tamagawa intervenant dans la constante est de la forme $$\underset{p}{\prod} \left(1-\frac{1}{p}\right)^{2^n-n}P_n\left(\frac{1}{p}\right)$$ avec $P_n$ un polynôme à coefficients entiers.
\begin{theor}[\cite{Blomer2014}]
Pour tout $U \subset A$, on définit ${\rm{ver}}(U) \subset A$ comme étant l'ensemble des sommets qui sont adjacents à au moins une arête de $U$. On a alors
$$
F(\mathbf{1})=\prod_p \sum_{k=0}^{2^n-2}\frac{b_k}{p^k}
$$
avec, pour tout $0 \leqslant k \leqslant 2^n-2$, 
$$
b_k=\sum_{U \subset A \atop {\rm{ver}}(U)=k} (-1)^{\#U}.
$$
De plus, on a $B_{0,n}=1$, $b_1=0$, $b_2=-\#A=-2^{n-1}\left(2^n+1\right)+3^n$ et la relation
$$
\sum_{k=0}^{2^n-2}b_k=0.
$$
\label{theor3}
\end{theor}
\noindent
\textit{Démonstration--} On renvoie au théorème 5 de la section 2 de \cite{Blomer2014} pour une preuve de ce résultat valable pour n'importe quel graphe $G=(S,A)$. La quantité $\#A$ s'obtient grâce à des arguments élémentaires de combinatoire.
\hfill
$\square$\\
\newline
\indent
Il n'apparaît pas évident \textit{a priori} de déterminer explicitement les valeurs de $b_k$ pour $k \geqslant 2$ pour des valeurs de $n$ quelconques ni de montrer directement que le polynôme $\overset{2^n-2}{\underset{k=0}{\sum}} b_k X^k$ admet 1 comme racine de multiplicité au moins $2^n-n-1$, comme cela est prévue par la conjecture de Peyre. On donne alors le lemme suivant qui permet d'expliciter un peu plus le produit eulérien définissant $F(\mathbf{1})$.
\begin{lemme}
On a
$$
F(\mathbf{1})=\prod_p \left(1-\frac{1}{p}\right)^{2^n-1} \left(1+ \sum_{k \geqslant 1} \frac{(k+1)^n-k^n}{p^k}\right).
$$
\label{lemmepdteulerien}
\end{lemme}
\noindent
\textit{Démonstration--}
Par définition de $F$ dans le Lemme \ref{holo} et au vu de (\ref{pe}), on a
$$
F(\mathbf{1})=\prod_p\left(1-\frac{1}{p}\right)^{2^n-1} G_p,
$$
où, pour tout nombre premier $p$, on a posé
 $$
 G_p=\sum_{(\nu_1,\dots,\nu_{N}) \in \mathbb{N}^N \atop \nu_i \nu_j=0 \hspace{0.5mm} {\rm{avec}} \hspace{0.5mm} (i,j) \in E_n} \frac{1}{p^{\nu_1+\cdots+\nu_{N}}}.
$$
En utilisant la bijection décrite en section 2 par le Lemme \ref{lemmefacto}, on obtient, pour tout nombre premier $p$, l'égalité
$$
G_p=\sum_{(k_1,\dots,k_{n}) \in \mathbb{N}^n } \frac{1}{p^{\max_{j} k_j}}.
$$
On s'inspire alors des calculs effectués en fin de section 4 de \cite{Br}. On regroupe la somme en fonction de la valeur $k:=\max_j k_j$. Posant $r_1=\#\{j \in \llbracket 1,n \rrbracket : k_j=k\}$ et $r_2=\#\{j \in \llbracket 1,n \rrbracket : 1 \leqslant k_j\leqslant k-1\}$, on aboutit à l'expression
$$
G_p=1+\sum_{\substack{k\geqslant1\\ 1 \leqslant r_1 \leqslant n \\ 0 \leqslant r_2 \leqslant n-r_1}} \frac{\rho(k,r_1,r_2)}{p^k},
$$
avec $\rho(k,r_1,r_2)$ le nombre de $\mathbf{k}=(k_1,\dots,k_n) \in \mathbb{N}^n$ correspondants. Or, on a
$$
\rho(k,r_1,r_2)=\binom{n}{r_1}\binom{n-r_1}{r_2}(k-1)^{r_2}
$$
si bien qu'on obtient
$$
G_p=1+ \sum_{k \geqslant 1} \frac{(k+1)^n-k^n}{p^k}.
$$
\hfill
$\square$\\
\newline
\indent
On dispose alors du lemme suivant qui montre que le rapport $F(\mathbf{1})/\zeta(n)$ est de la forme conjecturée par Peyre.
\begin{lemme}
Il existe un polynôme $P_n$ unitaire de degré $n-1$ tel que
$$
F(\mathbf{1})=\prod_p \left(1-\frac{1}{p}\right)^{2^n-n-1} P_n\left(\frac{1}{p}\right).
$$
Pour tout $n \in \mathbb{N}^{\ast}$, le polynôme est défini par la relation
\begin{eqnarray}
P_n(X)=(1-X)^{n+1} \sum_{k \geqslant 0} (k+1)^n X^k.
\label{pn}
\end{eqnarray}
\label{lemmepdteulerien2}
\end{lemme}
\noindent
\textit{Démonstration--}
Grâce à un changement d'indice, on déduit du Lemme \ref{lemmepdteulerien} l'expression
$$
F(\mathbf{1})=\prod_p \left(1-\frac{1}{p}\right)^{2^n} \sum_{k \geqslant 0} \frac{(k+1)^n}{p^k}.
$$
Une récurrence fournit alors l'existence de $P_n$ unitaire (palindromique) de degré $n-1$ tel que
$$
\sum_{k \geqslant 0} \frac{(k+1)^n}{p^k}=\frac{P_n\left(\frac{1}{p}\right)}{\left(1-\frac{1}{p}\right)^{n+1}},
$$
avec
$$
\forall n \in \mathbb{N}, \quad P_n(X)=(1-X)^{n+1} \sum_{k \geqslant 0} (k+1)^n X^k.
$$
Les polynômes $(P_n)_{n \geqslant 1}$ vérifient la relation de récurrence
$$
P_{n+1}(X)=(1+nX)P_n(X)+X(1-X)P'_n(X)
$$
et $P_1(X)=1$. Cela permet de conclure la preuve du lemme.
\hfill
$\square$\\
\newline
\indent

\subsubsection{Lien entre la quantité $F(\mathbf{1})/\zeta(n)$ et le nombre de Tamagawa associé à $X_{0,n}$}


On a prouvé lors du Lemme \ref{lemmepdteulerien2} que 
$$
F(\mathbf{1})=\prod_p \left(1-\frac{1}{p}\right)^{2^n-n-1} P_n\left(\frac{1}{p}\right),
$$
où 
$$
P_n(X)=\left(1-X\right)^{n+1}\sum_{k \geqslant 0} (k+1)^nX^k
$$
est un polynôme à coefficients entiers. D'après les remarques faisant suite au théorème 2 de Salberger présentes en Annexe et de manière plus rigoureuse d'après \cite[page 315]{Stem} et le Lemme \ref{lemme19}, on obtient que le facteur $p$-adique du nombre de Tamagawa est donné par
$$
\mu_p(X_{0,n}(\mathbb{A}_{\mathbb{Q}}))= \left(1-\frac{1}{p}\right)^{2^n-n-1} F\left(\mathfrak{S}_n,\frac{1}{p}\right)\left(1-\frac{1}{p^n}\right),
$$
où $F\left(\mathfrak{S}_n,t\right)$ est la fonction d'excédance de $\mathfrak{S}_n$ définie dans \cite{Stem}. D'après \cite[page 315]{Stem}, on a 
$$
F\left(\mathfrak{S}_n,X\right)=\left(1-X\right)^{n+1}\sum_{k \geqslant 0} (k+1)^nX^k
$$
si bien qu'on en déduit que $P_n=F\left(\mathfrak{S}_n,.\right)$. Ainsi, 
\begin{eqnarray}
\frac{F(\mathbf{1})}{\zeta(n)}=\prod_p \mu_p(X_{0,n}(\mathbb{A}_{\mathbb{Q}})).
\label{padic}
\end{eqnarray}
Pour conclure le traitement du nombre de Tamagawa
$$
\omega_H\left( X_{0,n}(\mathbb{A}_{\mathbb{Q}})^{{\rm Br}(X_{0,n})} \right)=\mu_{\infty}(X_{0,n}(\mathbb{A}_{\mathbb{Q}}))\prod_p \mu_p(X_{0,n}(\mathbb{A}_{\mathbb{Q}})),
$$
il reste à calculer la densité archimédienne. En utilisant le Lemme \ref{archimedien}, on obtient l'expression
\begin{eqnarray}
\mu_{\infty}(X_{0,n}(\mathbb{A}_{\mathbb{Q}}))=\!\!\int_{\mathbb{R}^{n-1}\times \mathbb{R}^{n-1}} \!\frac{\mbox{d}t_1\cdots\mbox{d}t_{n-1}\mbox{d}t_{n+1}\cdots\mbox{d}t_{2n-1}}{\left|t_{n+1}\cdots t_{2n-1}\right|\max\left(\left|t_1\right|^n,\dots,\left|t_{2n-1}\right|^n,\left| \frac{t_1}{t_{n+1}}+\cdots+\frac{t_{n-1}}{t_{2n-1}} \right|^n,1  \right)}.
\label{volumereel}
\end{eqnarray}
On démontre alors le résultat suivant.

\begin{lemme}
On a $\mu_{\infty}(X_{0,n}(\mathbb{A}_{\mathbb{Q}}))=n^2(n-1)!2^{n-1} \tilde{\beta}$, où $\tilde{\beta}$ a été défini en (\ref{betatilde}).
\label{omega}
\end{lemme}
\noindent
\textit{Démonstration.--} On part de l'expression (\ref{volumereel}) et on notera $I_{\infty}=\mu_{\infty}(X_{0,n}(\mathbb{A}_{\mathbb{Q}}))$ dans la suite de la démonstration. Il vient immédiatement l'égalité
$$
I_{\infty}\!=\!2^{n-1}\!\!\int_{\mathbb{R}^{n-1}\times \mathbb{R}^{n-1}\atop t_{n+1},\dots, t_{2n-1} \geqslant 0} \frac{\mbox{d}t_1\cdots\mbox{d}t_{n-1}\mbox{d}t_{n+1}\cdots\mbox{d}t_{2n-1}}{t_{n+1}t_{n+2}\cdots t_{2n-1}\max\left(\left|t_1\right|^n,\dots,t^n_{2n-1},\left| \frac{t_1}{t_{n+1}}+\cdots+\frac{t_{n-1}}{t_{2n-1}} \right|^n,1  \right)}.
$$
On utilise alors l'identité
$$
\forall s \geqslant 1, \qquad \frac{1}{s}=\int_{t \geqslant s} \frac{\mbox{d}t}{t^2}
$$
pour obtenir
$$
I_{\infty}=2^{n-1}\int_{\mathbb{R}^{n-1}\times \mathbb{R}^{n-1}\atop t_{n+1},\dots, t_{2n-1} \geqslant 0}\int_{t \geqslant \max\left(\left|t_1\right|^n,\dots,\left|t_{2n-1}\right|^n,\left| \frac{t_1}{t_{n+1}}+\cdots+\frac{t_{n-1}}{t_{2n-1}} \right|^n,1  \right)} \frac{\mbox{d}t_1\cdots\mbox{d}t_{n-1}\mbox{d}t_{n+1}\cdots\mbox{d}t_{2n-1}\mbox{d}t}{t_{n+1}t_{n+2}\cdots t_{2n-1}t^2}.
$$
Le changement de variables suivant
$$
\left\{
\begin{aligned}
&t=\frac{1}{u^n_{2n}}, \\
 &t_i=\frac{u_i}{u_{2n}} \quad \mbox{pour} \quad i \in \llbracket 1,2n-1\rrbracket \smallsetminus \{n\},\\
  \end{aligned}
  \right.
$$
de jacobien 
$$
\left|\det \begin{pmatrix}
-\frac{n}{u^{n+1}_{2n}}&0&\cdots & 0\\[3mm]
-\frac{u_1}{u^2_{2n}}&\frac{1}{u_{2n}}&\cdots & 0\\[3mm]
\vdots &\vdots &  \ddots & \vdots \\[3mm]
-\frac{u_{2n-1}}{u^2_{2n}}&0 & \cdots&   \frac{1}{u_{2n}}\\
 \end{pmatrix}\right|
=\frac{n}{u^{3n-1}_{2n}}
$$
fournit alors
$$
I_{\infty}=n2^{n-1}\int_{\substack{\mathbf{u} \in \mathbb{R}^{n-1}\times \mathbb{R}^{n}, \hspace{2mm} u_{n+1},\dots, u_{2n-1},u_{2n} \geqslant 0 \\ 1 \geqslant \max\left(\left|u_2\right|,\dots,u_{2n-1},u_{2n},u_{2n} \left| \frac{u_1}{u_{n+1}}+\cdots+\frac{u_{n-1}}{u_{2n-1}} \right|  \right)}} \!\frac{\mbox{d}u_1\cdots\mbox{d}u_{n-1}\mbox{d}u_{n+1}\cdots\mbox{d}u_{2n-1}\mbox{d}u_{2n}}{u_{n+1}u_{n+2}\cdots u_{2n-1}}.
$$
\'Etablissons à présent qu'on peut supposer que $u_{n+1} \leqslant u_{n+2} \leqslant \cdots \leqslant u_{2n}$ dans l'intégrale ci-dessus quitte à la mutliplier par un facteur $n!$. 
Pour toute permutation $\sigma$ de $\llbracket n+1,2n \rrbracket$, on notera
$$
I_{\infty, \sigma}=\int_{\substack{\mathbf{u} \in \mathbb{R}^{n-1}\times \mathbb{R}^{n}, \hspace{2mm} 0 \leqslant u_{\sigma(n+1)}\leqslant \dots\leqslant  u_{\sigma(2n)} \\ 1 \geqslant \max\left(\left|u_1\right|,\dots,u_{2n-1},u_{2n},u_{2n} \left| \frac{u_1}{u_{n+1}}+\cdots+\frac{u_{n-1}}{u_{2n-1}} \right|  \right)}} \frac{\mbox{d}u_1\cdots\mbox{d}u_{n-1}\mbox{d}u_{n+1}\cdots\mbox{d}u_{2n-1}\mbox{d}u_{2n}}{u_{n+1}u_{n+2}\cdots u_{2n-1}}.
$$
Si $\sigma(2n)=2n$, il est facile de voir que le changement de variables
\begin{eqnarray}
\left\{
\begin{aligned}
&v_i=u_{\sigma(n+i)-n} \quad \mbox{pour} \quad  1 \leqslant i \leqslant n-1, \\
&  v_i=u_{\sigma(n+i)}\quad \mbox{pour} \quad 1 \leqslant i \leqslant n-1 \\
& v_{2n}=u_{2n}
\end{aligned}
\right.
\label{chgt}
\end{eqnarray}
est de jacobien
$$
\left| \det\begin{pmatrix}
P_{\tilde{\sigma}} & 0&0\\
0& P_{\tilde{\sigma}}&0 \\
0&0&1\\
\end{pmatrix}
 \right|=\varepsilon(\tilde{\sigma})^2=1,
$$
où $\tilde{\sigma}$ est la permutation de $\llbracket n+1,2n-1\rrbracket$ induite par $\sigma$, $P_{\tilde{\sigma}}$ est la matrice de permutation associée à $\tilde{\sigma}$ et $\varepsilon(\tilde{\sigma})$ est la signature de la permutation $\tilde{\sigma}$. Ce changement de variables fournit alors la relation
$
I_{\infty, \sigma}=I_{\infty, {\rm Id}}.
$
Supposons à présent que $\sigma(2n) < 2n$. Dans ce cas, on considère $m=\sigma(2n)-n \in \llbracket 1,n-1 \rrbracket$. On effectue alors le changement de variables:
$$
\left\{
\begin{aligned}
&v_i=u_i \quad \mbox{pour} \quad i \in \llbracket 1,2n-1\rrbracket\smallsetminus\{n, m\}\\
&v_m=-u_{2n}\left( \frac{u_1}{u_{n+1}}+\frac{u_2}{u_{n+2}}+\cdots+\frac{u_{n-1}}{u_{2n-1}} \right)\\
\end{aligned}
\right.
$$
de jacobien $\frac{u_m}{u_{2n}}$. Il s'ensuit alors que $I_{\infty, \sigma}$ est égal à
$$
\int_{\substack{ \mathbf{v} \in\mathbb{R}^{n-1}\times \mathbb{R}^{n}, \hspace{2mm} 0  \leqslant v_{\sigma(n+1)}\leqslant \dots\leqslant  v_{\sigma(2n)} \\ 1 \geqslant \max\left(\left|v_1\right|, \dots ,v_{2n},v_{n+m}\left| \frac{v_1}{v_{n+1}}+\cdots+\frac{v_{m-1}}{v_{n+m-1}}+\frac{v_m}{v_{2n}}+\frac{v_{m+1}}{v_{n+m+1}}+\cdots+\frac{v_{n-1}}{v_{2n-1}} \right| \right)}} \!\!\!\!\frac{\mbox{d}v_1\cdots\mbox{d}v_{n-1}\mbox{d}v_{n+1}\cdots\mbox{d}v_{2n-1}\mbox{d}v_{2n}}{v_{n+1}\cdots v_{m-1}v_{m+1}\cdots v_{2n}}
$$
car
$$
\left|u_{m}\right|=u_{n+m}\left|  \frac{u_1}{u_{n+1}}+\cdots+\frac{u_{m-1}}{u_{n+m-1}}+\frac{v_m}{u_{2n}}+\frac{u_{m+1}}{u_{n+m+1}}+\cdots+\frac{v_{n-1}}{v_{2n-1}}\right|.
$$
Par choix de $m$, on se retrouve alors à nouveau dans un cas où $v_{n+m}=\max\left\{ v_{n+i} \hspace{1mm} : \hspace{1mm} 1 \leqslant i \leqslant n  \right\}$ et un changement de variables similaire à (\ref{chgt}) permet d'obtenir l'égalité
$
I_{\infty, \sigma}=I_{\infty, {\rm Id}}
$
également dans ce cas et ainsi de conclure que pour toute permutation  $\sigma$ de $\llbracket n+1,2n \rrbracket$, on a $
I_{\infty, \sigma}=I_{\infty, {\rm Id}}.
$ Par conséquent, il vient
\begin{eqnarray}
I_{\infty}=n2^{n-1}n!\int_{\substack{\mathbf{u} \in \mathbb{R}^{n-1}\times \mathbb{R}^{n}, \hspace{2mm} 0 \leqslant u_{n+1}\leqslant \dots\leqslant  u_{2n} \\ 1 \geqslant \max\left(\left|u_1\right|,\dots,u_{2n-1},u_{2n},u_{2n} \left| \frac{u_1}{u_{n+1}}+\cdots+\frac{u_{n-1}}{u_{2n-1}} \right|  \right)}} \frac{\mbox{d}u_1\cdots\mbox{d}u_{n-1}\mbox{d}u_{n+1}\cdots\mbox{d}u_{2n-1}\mbox{d}u_{2n}}{u_{n+1}u_{n+2}\cdots u_{2n-1}}.
\label{int}
\end{eqnarray}
On effectue alors le dernier changement de variables suivant, à $(u_1,\dots,u_{n-1},u_{2n})$ fixés,
$$
\left\{
\begin{aligned}
&v_1=\frac{u_{2n}}{u_{n+1}}\\[2mm]
&v_i=\frac{u_{n+i-1}}{u_{n+i}} \quad \mbox{pour} \quad i \in \llbracket 2,n-1\rrbracket\\
\end{aligned}
\right. \quad \mbox{soit} \quad \left\{
\begin{aligned}
&u_{n+1}=\frac{u_{2n}}{v_1}\\[2mm]
&u_{n+i}=\frac{u_{2n}}{\displaystyle \prod_{\scriptscriptstyle 1 \leqslant k \leqslant i}v_k} \quad \mbox{pour} \quad i \in \llbracket 2,n-1\rrbracket\\
\end{aligned}
\right.
$$
de jacobien
$$
\frac{u^{n-1}_{2n}}{\displaystyle\underset{\scriptscriptstyle 1 \leqslant k \leqslant n-1}{\prod}v^{n-k+1}_k}
$$
puisque la matrice jacobienne est une matrice triangulaire inférieure de coefficients diagonaux
$$
-\frac{u_{2n}}{\displaystyle v^2_i\underset{\scriptscriptstyle 1 \leqslant k \leqslant i-1}{\prod}v_k}
$$
pour tout $i \in \llbracket 1,n-1\rrbracket.$ Par conséquent, commme
$$
\frac{1}{u_{n+1}u_{n+2}\cdots u_{2n-1}}=\frac{u^{n-1}_{2n}}{\displaystyle\underset{\scriptscriptstyle 1 \leqslant k \leqslant n-1}{\prod}v^{n-k}_k},
$$
il s'ensuit que l'intégrale apparaissant dans (\ref{int}) est égale à
$$
\int_{\substack{0 \leqslant v_1,\dots,v_{n-1} \leqslant 1  \\ 1 \geqslant \max\left(\left|u_1\right|,\dots,|u_{n-1}|, \left| u_1v_1+\cdots+u_1v_1\cdots v_{n-1} \right|  \right)}} \int_{0 \leqslant u_{2n} \leqslant v_1v_2\cdots v_{n-1}}\frac{\mbox{d}u_1\cdots\mbox{d}u_{n-1}\mbox{d}v_{1}\cdots\mbox{d}v_{n-1}\mbox{d}u_{2n}}{v_{1}v_{2}\cdots v_{n-1}}=\tilde{\beta},
$$
où la dernière égalité est obtenue en intégrant par rapport à $u_{2n}$. Cela permet de conclure à l'égalité $\mu_{\infty}(X_{0,n}(\mathbb{A}_{\mathbb{Q}}))=n!2^{n-1}n \tilde{\beta}$ et achève la démonstration du lemme.
\hfill
$\square$

\subsection{Le dénouement}

Les Lemmes \ref{omega}, \ref{alpha}, la relation (\ref{padic}) ainsi que le Théorème \ref{theor1} impliquent alors que la constante $c_n=c_{\rm{Peyre}}$ est en accord avec la prédiction de Peyre et cela achève la démonstration des conjectures de Manin et Peyre pour les hypersurfaces projectives $W_n$.

%

\bibliographystyle{plain-fr}
\nocite{*}
\bibliography{bibliogr3}
\addcontentsline{toc}{section}{Annexe de Per Salberger}

\end{document}


\renewcommand{\contentsname}{Sommaire}
\renewcommand{\abstractname}{Abstract}
\setcounter{section}{0}
\rhead{}
\chead{}
\lhead{}
\fancyhead[CO]{\footnotesize{ \textsc{A crepant resolution}}}
\fancyhead[CE]{\footnotesize{\textsc{P. Salberger}}}
\subjclass{11G35, 14G05, 14J40.}
\keywords{toric variety, crepant resolution.}

\title{\textbf{Annexe: }\\ \textbf{A crepant resolution}}
\maketitle

\setcounter{tocdepth}{2}
\addtocounter{page}{43}
In a previous paper with Blomer and Brüdern \cite{Blomer2014}, we constructed a partial resolution $f_n:X_n \rightarrow W_n$ of the projective hypersurface $W_n \subseteq \mathbb{P}^{2n-1}$ ($n \geqslant 1$) with homogeneous coordinates $$(x_1,\dots,x_n,y_1,\dots,y_n)$$ defined by the equation
$$
x_1y_2\cdots y_{n-1}y_n+x_2y_1y_3\cdots y_n+\cdots=x_ny_1y_2\cdots y_{n-1}=0.
$$
\newline
The variety $X_n \subseteq \mathbb{P}^{2n-1} \times \mathbb{P}^{n-1} \times \mathbb{P}^{n-1}$ is the triprojective variety with trihomogeneous coordinates $(x_1,\dots,x_n,y_1,\dots,y_n;Y_1,\dots,Y_n;Z_1,\dots,Z_n)$ defined by the following equations
%
\begin{numcases}{}
 x_1Z_1+\cdots+x_nZ_n=0, \label{1}\\
 y_iY_j-y_jY_i=0 \quad \mbox{ for } \quad  1 \leqslant i<j\leqslant n, \label{2}\\
 Y_1Z_1=\cdots=Y_nZ_n. \label{3}
\end{numcases}
In \cite{Blomer2014}, we proved that the projection $\mbox{pr}_1: \mathbb{P}^{2n-1} \times \mathbb{P}^{n-1} \times \mathbb{P}^{n-1} \rightarrow \mathbb{P}^{2n-1}$ restricts to a birational projective morphism $f_n:X_n \rightarrow W_n$ which is \textit{crepant}. We have thus $f_n^{\ast} \omega_{W_n} \cong \omega_{X_n}$.
\newline
\par
We also showed for $n \geqslant 1$ in \cite{Blomer2014} that the restriction of $\mbox{pr}_2 \times \mbox{pr}_3: \mathbb{P}^{2n-1} \times \mathbb{P}^{n-1} \times \mathbb{P}^{n-1} \rightarrow \mathbb{P}^{n-1} \times \mathbb{P}^{n-1}$ to $X_n$ gives rise to a morphism $\lambda_n: X_n \rightarrow B_n$ such that $X_n$ is a $\mathbb{P}^{n-1}$-bundle over the biprojective variety $B_n \subseteq \mathbb{P}^{n-1} \times \mathbb{P}^{n-1}$ with bihomogeneous coordinates $(Y_1,\dots,Y_n;Z_1,\dots,Z_n)$ defined by the equations $(\ref{3})$. Note that for $n=1$, there is no equation and $B_1=\mathbb{P}^{0} \times \mathbb{P}^{0}$ is a point.
\newline
\par
Hence, $X_n$ is non-singular for $n \leqslant 3$ as $B_n$ is non-singular for $n \leqslant 3$. We have thus a crepant desingularisation $f_n:X_n \rightarrow W_n$ for $n \leqslant 3$. But for $n>3$, both $B_n$ and $X_n$ are singular. The aim of this note is to construct a crepant resolution $p_{0,n}:B_{0,n} \rightarrow B_n$ for all $n \geqslant 1$. This will give a crepant desingularisation $p_{X_n}: X_{0,n}=B_{0,n} \times_{B_n} X_n \rightarrow X_n$ by pulling back $\lambda_n:X_n \rightarrow B_n$ along $p_{0,n}$ and hence a crepant resolution $p_{X_n} \circ f_n: X_{0,n} \rightarrow W_n$ for all $n \geqslant 1$.
\newline
\par
For non-empty subsets $I=\{i_1,i_2,\dots,i_k\} \subseteq \{1,\dots,n\}$, we write $\mathbb{P}^{I}\times \mathbb{P}^{I}$ for the biprojective space $\mathbb{P}^{k-1}\times \mathbb{P}^{k-1}$ with bihomogeneous coordinates $\left(\mathbf{Y}^{I};\mathbf{Z}^{I}\right)=\left( Y^{I}_{i_1},\dots,Y^{I}_{i_k};Z^{I}_{i_1},\dots,Z^{I}_{i_k} \right)$ and $B^{I} \subseteq \mathbb{P}^{I}\times \mathbb{P}^{I}$ for the closed subvariety defined by
\begin{eqnarray}
Y^{I}_{i_1}Z^{I}_{i_1}=\cdots=Y^{I}_{i_k}Z^{I}_{i_k}. \label{4}
\end{eqnarray}
In particular, $B^{I}=B_n$ for $I=[n]:=\{1,\dots,n\}$. We will never allow $I=\varnothing$ in this note. For products $\prod_{I \subseteq [n]} V^{I}$ of varieties indexed by $I \subseteq [n]$, it will thus be understood that $I \neq \varnothing$.\\

\par
Now let $B_{0,n}$ be the closed subvariety of $\prod_{I \subseteq [n]} B^{I}$ defined by the equations
\begin{numcases}{}
Y^{I}_k Y^{J}_{\ell}=Y^{I}_{\ell} Y^{J}_{k} \quad \mbox{ for } \quad k, \ell \in J \subsetneq I \subseteq [n] \label{5}\\
Z^{I}_k Z^{J}_{\ell}=Z^{I}_{\ell} Z^{J}_{k} \quad \mbox{ for } \quad k, \ell \in J \subsetneq I \subseteq [n] \label{6}
\end{numcases} 
There is an obvious morphism $p_{0,n}: B_{0,n} \rightarrow B_n$ which sends $\prod_{I \subseteq [n]} \left( \mathbf{Y}^{I};\mathbf{Z}^{I} \right)$ to $\left( \mathbf{Y}^{[n]};\mathbf{Z}^{[n]} \right)$.\\
\par
We shall also consider the Coxeter toric variety $C_n \subseteq \prod_{I \subseteq [n]} \mathbb{P}^{I}$ of $\mathfrak{S}_n$ with multihomogeneous coordinates $\prod_{I \subseteq [n]} \left( \mathbf{Y}^{I} \right)$ defined by the equations $(\ref{5})$. It follows from a result of De Concini and Procesi (see \cite[Theorem 4.2]{CP95} and \cite[Prop. 2.9]{He}) that $C_n$ is a non-singular projective irreducible toric variety of dimension $n-1$.
\begin{prop}
Let $B^{\ast}_n \subsetneq B_n$ be the open subset where we do not have $Y_k=Z_k=0$ for any $k \in [n]$ and $B^{\ast}_{0,n}=p_{0,n}^{-1}\left( B^{\ast}_n\right)$. Then
\begin{enumerate}[(a)]
\item
The morphism $\mbox{pr}_1: B_{0,n} \rightarrow C_n$, which sends $\prod_{I \subseteq [n]} \left( \mathbf{Y}^{I};\mathbf{Z}^{I}\right)$ to $\prod_{I \subseteq [n]} \left( \mathbf{Y}^{I}\right)$ is an isomorphism.
\item
$B_{0,n}$ is a non-singular projective irreducible toric variety.
\item
$B_{0,n} \smallsetminus B^{\ast}_{0,n}$ is of codimension at least two in $B_{0,n}$.
\item
$p_{0,n}$ is surjective and maps $B^{\ast}_{0,n}$ isomorphically onto $B^{\ast}_n$.
\end{enumerate}
\end{prop}
\noindent
\textit{Proof.--} \textit{(a)} The morphism $\mbox{pr}_1: B_{0,n} \rightarrow C_n$ is dominant as any $\prod_{I \subseteq [n]} \left( \mathbf{Y}^{I}\right) \in C_n$ with $Y^{I}_k \neq 0$ for all $I \subseteq [n]$ and $k \in I$ is the image of $\prod_{I \subseteq [n]} \left( \mathbf{Y}^{I};\mathbf{Z}^{I}\right) \in B_{0,n}$ with $Z^{I}_k=1/Y^{I}_k$ for all $I \subseteq [n]$ and $k \in I$. As $C_n$ is integral, it is thus enough to show that $\mbox{pr}_1$ is a closed immersion and since $\mbox{pr}_1$ is proper and finitely presented that $\mbox{pr}_1$ is injective (see \cite[prop 8.11.5]{EGA}). That is, we have to show that $\prod_{I \subseteq [n]} \left( \mathbf{Z}^{I}\right)  \in \prod_{I \subseteq [n]} \mathbb{P}^{I}$ is uniquely determined by $\prod_{I \subseteq [n]} \left( \mathbf{Y}^{I}\right)  \in \prod_{I \subseteq [n]} \mathbb{P}^{I}$ for $\prod_{I \subseteq [n]} \left( \mathbf{Y}^{I};\mathbf{Z}^{I}\right) \in B_{0,n}$.
\par
To show this, we use induction with respect to $n$. If $n=1$, then $\mbox{pr}_1$ is just the projection of $\mathbb{P}^0 \times \mathbb{P}^0$ onto $\mathbb{P}^0$ and the assertion is trivial. So suppose $n>1$ and that  $\prod_{I \subseteq [n]} \left( \mathbf{Y}^{I};\mathbf{Z}^{I}\right) \in B_{0,n}$. Then $\mathbf{Z}^{I} \in \mathbb{P}^{I}$ is uniquely determined by $\prod_{I \subseteq [n]} \left( \mathbf{Y}^{I}\right)  \in \prod_{I \subseteq [n]} \mathbb{P}^{I}$ for any $I \neq [n]$ by the induction hypothesis. Also, if $Y^{[n]}_i=0$ for some $i \in [n]$, then $Z^{[n]}_j=0$ for some $j \in [n]$. If, say $Z^{[n]}_n=0$, then $\left( \mathbf{Z}^{[n]} \right)=\left(Z^{[n-1]}_1,\dots,Z^{[n-1]}_{n-1},0 \right)$ in $\mathbb{P}^{n-1}$ by (\ref{6}). Finally, if $Y^{[n]}_i\neq 0$ for all $i \in [n]$, then $\left(\mathbf{Z}^{[n]}\right)=\left(1/Y^{[n]}_1,\dots, 1/Y^{[n]}_n \right) \in \mathbb{P}^{[n]}$. Hence $\left(\mathbf{Z}^{[n]}\right)$ is uniquely determined by $\prod_{I \subseteq [n]} \left( \mathbf{Y}^{I}\right)$ also for points on $B_{0,n}$.\\
\par
\noindent
\textit{(b)} This follows from \textit{(a)} and the above result for the Coxeter variety $C_n$ due to De Concini and Procesi \cite{CP95}.\\
\par
\noindent
\textit{(c)} It suffices to show that the closed subset $F_k \subseteq B_{0,n}$ defined by the equations $Y^{[n]}_k=Z^{[n]}_k=0$ is of codimension at least 2 in $B_{0,n}$ for all $k \in [n].$ But for such points in $B_{0,n}$, we see from (\ref{5}) and (\ref{6}), that $Y^{I}_k=Z^{I}_k=0$ for all $I \subseteq [n]$ where $k \in I$ and $\# I \geqslant 2$. The projection from $F_k$ to $\prod_{I \subseteq [n] \smallsetminus \{k\}} B^{I}$ is thus a closed immersion and, by \textit{(a)}, isomorphic to a proper variety of the $(n-2)-$dimensional Coxeter subvariety of $\prod_{I \subseteq [n] \smallsetminus \{k\}} \mathbb{P}^{I}$ defined by (\ref{5}) and (\ref{6}) for $J \subsetneq I \subseteq [n] \smallsetminus \{k\}$. This proves that $\dim(F_k) \leqslant n-3$.\\
\par
\noindent
\textit{(d)} It is clear that the open subset of $B_{0,n}$ with $Y^{[n]}_1Z^{[n]}_1\neq 0$ is mapped isomorphically onto the open subset of $B_n$ with $Y^{[n]}_1Z^{[n]}_1 \neq 0$ as $\left( \mathbf{Y}^{I};\mathbf{Z}^{I}\right) \in \mathbb{P}^{I} \times \mathbb{P}^{I}$ is the projection of $\left( \mathbf{Y}^{[n]};\mathbf{Z}^{[n]}\right) \in \mathbb{P}^{[n]} \times \mathbb{P}^{[n]}$ for any $\prod_{I \subseteq [n]} \left( \mathbf{Y}^{I};\mathbf{Z}^{I}\right) \in B_{0,n}$ with  $Y^{[n]}_1Z^{[n]}_1=\cdots=Y^{[n]}_nZ^{[n]}_n\neq 0$. Hence $p_{0,n}\left( B^{\ast}_{0,n}\right)$ is dense in $B^{\ast}_n$. It is thus just as in \textit{(a)} enough to show that the map from $B^{\ast}_{0,n}$ to $B^{\ast}_n$ is injective. As we have already seen that $\prod_{I \subseteq [n]} \left( \mathbf{Y}^{I};\mathbf{Z}^{I}\right) \in B_{0,n}$ is uniquely determined by $\left( \mathbf{Y}^{[n]};\mathbf{Z}^{[n]}\right) \in \mathbb{P}^{[n]} \times \mathbb{P}^{[n]}$ if $Y^{[n]}_1Z^{[n]}_1\neq 0$, suppose instead that $Y^{[n]}_1Z^{[n]}_1= 0$. We may then, if $\prod_{I \subseteq [n]} \left( \mathbf{Y}^{I};\mathbf{Z}^{I}\right) \in B^{\ast}_{0,n}$, find a partition $[n]=J \cup K$ such that $Y^{[n]}_j=0$ if, and only if, $j \in J$ and $Z^{[n]}_k=0$ if, and only if, $k \in K$. It now follows from (\ref{5}) and (\ref{6}) that $\left( \mathbf{Y}^{I};\mathbf{Z}^{I}\right) \in \mathbb{P}^{I} \times \mathbb{P}^{I}$ is uniquely determined by $\left( \mathbf{Y}^{[n]};\mathbf{Z}^{[n]}\right)$. If $I$ is not contained in $J$, then $\mathbf{Y}^{I} \in  \mathbb{P}^{I}$ is just the projection of $\mathbf{Y}^{[n]}$ and if $I \subseteq J$, then $\mathbf{Y}^{I} \in  \mathbb{P}^{I}$ is represented by a tuple obtained by inverting the $I$-coordinates of $\mathbf{Z}^{[n]}$. Similarly, if $I$ is not contained in $K$, then  $\mathbf{Z}^{I} \in  \mathbb{P}^{I}$ is just the projection of $\mathbf{Y}^{[n]}$ and if $I \subseteq K$, then $\mathbf{Z}^{I} \in  \mathbb{P}^{I}$ is given by a tuple obtained by inverting the $I$-coordinates of $\mathbf{Y}^{[n]}$. This completes the proof.
\hfill
$\square$\\
\newline
\noindent
We now pull back the $\mathbb{P}^{n-1}$-bundle $\lambda_n:X_n \rightarrow B_n$ along $p_{0,n}:B_{0,n} \rightarrow B_n$. This gives us a $\mathbb{P}^{n-1}$-bundle $\lambda_{0,n}: B_{0,n} \times_{B_n} X_n \rightarrow B_{0,n}$. One can also describe $X_{0,n}$ more directly. It is the closed subvariety $X_{0,n} \subseteq \mathbb{P}^{2n-1} \times \prod_{I \subseteq [n]} \mathbb{P}^{I}\times \mathbb{P}^{I}$ with multihomogeneous coordinates $$\left(x_1,\dots,x_n,y_1,\dots,y_n;\prod_{I \subseteq [n]} \left( \mathbf{Y}^{I};\mathbf{Z}^{I}\right)\right)$$ defined by (\ref{4}), (\ref{5}), (\ref{6}) and (\ref{7}) and (\ref{8}) below
\begin{numcases}{}
x_1 Z^{[n]}_1+\cdots+x_n Z^{[n]}_n=0, \label{7}\\
 y_iY^{[n]}_j-y_jY^{[n]}_i=0 \quad \mbox{ for } \quad  1 \leqslant i<j\leqslant n. \label{8}
\end{numcases}

\begin{theor}
The projection $p_{X_n}:X_{0,n}= B_{0,n} \times_{B_n} X_n \rightarrow X_n$ which sends 
$$
\left(x_1,\dots,x_n,y_1,\dots,y_n;\prod_{I \subseteq [n]} \left( \mathbf{Y}^{I};\mathbf{Z}^{I}\right)\right) \longmapsto \left(x_1,\dots,x_n,y_1,\dots,y_n;\mathbf{Y}^{[n]};\mathbf{Z}^{[n]}\right)
$$
is a crepant desingularisation of $X_n$.
\end{theor}
\noindent
\textit{Proof.--} The variety $X_{0,n}$ is non-singular as it is a $\mathbb{P}^{n-1}$-bundle over a non-singular variety $B_{0,n}$. It is also clear from the Proposition 1 that the open subset $X^{\ast}_{0,n}=\lambda_{0,n}^{-1}\left(B^{\ast}_{0,n}\right) \subseteq X_{0,n}$ is mapped isomorphically onto the open subset $X^{\ast}_n=\lambda_n^{-1}\left(B^{\ast}_n\right) \subseteq X_n$. Therefore, the restrictions of $\omega_{X_{0,n}}$ and $p_{X_n}^{\ast} \omega_{X_n}$ to $X^{\ast}_{0,n}$ are isomorphic. By the Proposition 1, we have also that $X_{0,n} \smallsetminus X^{\ast}_{0,n}$ is of codimension at least two in $X_{0,n}$ and hence the restriction from $\mbox{Pic}(X_{0,n})$ to $\mbox{Pic}\left(X^{\ast}_{0,n}\right)$ is injective (see proposition 6.5b and corollary II.6.16 in \cite{Ha}). Hence $\omega_{X_{0,n}}$ and $p_{X_n}^{\ast} \omega_{X_n}$ are canonically isomorphic invertible sheaves on $X_{0,n}$.\hfill
$\square$\\
\newline
\noindent
\textbf{Remark.--} If $n \leqslant 3$, then $p_{X_n}$ is an isomorphism as $p_{0,n}$ is an isomorphism for such $n$.\\
\newline
\noindent
We can now derive the main theorem of this note.
\begin{theor}
Let $n \geqslant 1$ and $W_n \subseteq \mathbb{P}^{2n-1}$ be the projective hypersurface defined by the following equation
$$
x_1y_2\cdots y_{n-1}y_n+x_2y_1y_3\cdots y_n+\cdots=x_ny_1y_2\cdots y_{n-1}=0.
$$
Let $X_{0,n} \subseteq \mathbb{P}^{2n-1} \times \prod_{I \subseteq [n]} \mathbb{P}^{I}\times \mathbb{P}^{I}$ be the closed subvariety with multihomogeneous coordinates $$\left(x_1,\dots,x_n,y_1,\dots,y_n;\prod_{I \subseteq [n]} \left( \mathbf{Y}^{I};\mathbf{Z}^{I}\right)\right)$$ defined by the equations (\ref{4}), (\ref{5}), (\ref{6}), (\ref{7}) and (\ref{8}). Then the projection $\mbox{pr}_1: \mathbb{P}^{2n-1} \times \prod_{I \subseteq [n]} \mathbb{P}^{I}\times \mathbb{P}^{I} \rightarrow \mathbb{P}^{2n-1}$ restricts to a morphism $f_{0,n}:X_{0,n} \rightarrow W_n$ which is a crepant resolution of singularities of $W_n$. Moreover, $X_{0,n}$ is a $\mathbb{P}^{n-1}$-bundle over a variety $B_{0,n}$ isomorphic to the Coxeter toric variety $C_n$ of $\mathfrak{S}_n$.
\end{theor}
\noindent
\textit{Proof.--} The map $f_{0,n}$ is the composition of $p_{X_n}:X_{0,n} \rightarrow X_n$ and $f_{n}:X_n \rightarrow W_n$. We proved in \cite{Blomer2014} that $f_n$ is birational and crepant. We thus obtain the desired result by combining the results on $f_n$ and $X_n$ in \cite{Blomer2014} with the results on $p_{X_n}$ and $X_{0,n}$ from the Theorem 1. \hfill
$\square$\\
\newline
\noindent
\textbf{Remarks.--} By (a slight generalisation of) Manin's conjecture, we expect as a consequence of the Theorem 2 that we have $\mbox{rk}\left(\mbox{Pic}(X_{0,n})\right)-1=\mbox{rk}\left(\mbox{Pic}(B_{0,n})\right)=\mbox{rk}\left(\mbox{Pic}(C_n)\right)$ log factors in the asymptotic formula for the counting function on the open subset of $W_n$ with $y_1y_2\cdots y_{n-1}y_n\neq 0$.\\
\newline
\noindent
We also expect that the constant in the main term of this asymptotic formula is given by Peyre's Tamagawa constant for $X_{0,n}$, which is interpreted as an adelic volume of the universal torsor $\boldsymbol{\mathit{T}}_0$ over $X_{0,n}$ in \cite{Sal}. The $p$-adic factor of Peyre's Tamagawa constant should thus be
$$
\frac{\#\boldsymbol{\mathit{T}}_0\left(\mathbb{F}_p\right)}{\dim\left(\boldsymbol{\mathit{T}}_0\right)}=\left(1-\frac{1}{p}\right)^{{\rm{rk}}\left({\rm{Pic}}(X_{0,n})\right)}\frac{\# X_{0,n}\left(\mathbb{F}_p\right)}{p^{\dim(X_{0,n})}}=\left(1-\frac{1}{p}\right)^{{\rm{rk}}\left({\rm{Pic}}(C_n)\right)}\frac{\# C_n\left(\mathbb{F}_p\right)}{p^{\dim(C_n)}}\left(1-\frac{1}{p^n}\right).
$$
Here clearly, $\dim(C_n)=n-1$. The invariants ${\rm{Pic}}(C_n)$ and $\# C_n\left(\mathbb{F}_p\right)$ for the Coxeter toric variety $C_n$ are also known and may be found in \cite{Stem}. We have
$$
{\rm rk}\left({\rm{Pic}}(C_n)\right)=2^n-n-1
$$
as ${\rm{rk}}\left({\rm{Pic}}(C_n)\right)$ is equal to the second Betti number of $C_n$ and
$$
\# C_n\left(\mathbb{F}_p\right)=F\left(\mathfrak{S}_n,p\right)
$$
where $F\left(\mathfrak{S}_n,t\right)$ is the symmetric polynomial of degree $n-1$ in $t$ given by the associated excedance function \cite[pages 311 and 316]{Stem}. Recall here that the excedance number of any $w \in \mathfrak{S}_n$ is the quantity $e(w)=\#\{i \in [n] \hspace{1mm}:\hspace{1mm}w(i)>i\}$ (see \cite[page 309]{Stem}).\\
\newline
\noindent
\textbf{Examples.--} In $\mathfrak{S}_3$, we have $e(\mbox{Id})=0$, $e((12))=1$, $e((13))=1$, $e((23))=1$, $e((132))=1$ and $e((123))=2$. The generating excedance function $F\left(\mathfrak{S}_3,t\right)$ is thus $1 \times t^0+4 \times t^1+1\times t^2=1+4t+t^2$ and $\# C_3\left(\mathbb{F}_p\right)=1+4p+p^2$.\\
For $n=3$, we thus expect $\mbox{rk}\left(\mbox{Pic}(C_3)\right)=2^3-4=4$ log factors and the $p-$adic factor 
$$
\left(1-\frac{1}{p}\right)^4\left(1+\frac{4}{p}+\frac{1}{p^2}\right)\left(1-\frac{1}{p^3}\right)
$$
in the main term.\\
\newline
\noindent
In $\mathfrak{S}_4$, we have:
\begin{itemize}
\item[$\bullet$] $e(\mbox{Id})=0$;
\item[$\bullet$] For $\sigma$ one of the six transpositions or $\sigma=(132)$, $(142)$, $(143)$ or $(1432)$, one has $e(\sigma)=1$;
\item[$\bullet$] For $\sigma$ a product of two disjoint transpositions, a three cycle distinct from $(132)$, $(142)$ and $(143)$ and a four cycle distinct from $(1234)$, one hase $e(\sigma)=2$;
\item[$\bullet$] $e((1234))=3$.
\end{itemize}
So $F\left(\mathfrak{S}_4,t\right)=1+11t+11t^2+t^3$ and $\# C_4\left(\mathbb{F}_p\right)=1+11p+11p^2+p^3$.\\
For $n=4$, we thus expect $\mbox{rk}\left(\mbox{Pic}(C_4)\right)=2^4-5=11$ log factors and the $p-$adic factor 
$$
\left(1-\frac{1}{p}\right)^{11}\left(1+\frac{11}{p}+\frac{11}{p^2}+\frac{1}{p^3}\right)\left(1-\frac{1}{p^4}\right)
$$
in the main term.\\
\newline
\noindent
Similarly, for $n=5$, we expect $\mbox{rk}\left(\mbox{Pic}(C_5)\right)=2^5-6=26$ log factors and the $p-$adic factor 
$$
\left(1-\frac{1}{p}\right)^{26}\left(1+\frac{26}{p}+\frac{66}{p^2}+\frac{26}{p^3}=\frac{1}{p^4}\right)\left(1-\frac{1}{p^5}\right)
$$
in the main term.\\
\newline
\noindent
\textbf{Remarks.--} The odd Betti numbers vanish for any smooth projective toric variety. We have thus by the Weil conjectures (see \cite[Appendix]{Ha}) that
$$
\# C_n\left(\mathbb{F}_p\right)=\beta_{2(n-1)}p^{n-1}+\beta_{2(n-2)}p^{n-2}+\cdots+\beta_2 p+\beta_0 
$$
and that
$$
F\left(\mathfrak{S}_n,t\right)=\beta_{2(n-1)}t^{n-1}+\beta_{2(n-2)}t^{n-2}+\cdots+\beta_2 t+\beta_0 
$$
where $\beta_i=\dim\left( H^{i}\left( C_n(\mathbb{C})_{{\rm{an}}},\mathbb{Q} \right)\right)$ is the $i$-th Betti number of $C_n$. In \cite{Stem}, Stembridge gives thus a combinatorial interpretation of the Betti numbers of Coxeter toric varieties in terms of the excedance function of $\mathfrak{S}_n$.\\
For more background on the cohomology of toric varieties, see also \cite[section 4.5]{Fu}.
\bibliographystyle{alpha}
\bibliography{bibliogr3}